\documentstyle[12pt]{article}

\makeatletter

\newfam\gotfam

\font\fivgot = eufm5
\font\sixgot = eufm5	\@magscale1
\font\sevgot = eufm7
\font\egtgot = eufm7 \@magscale1
\font\tengot = eufm10
\font\elvgot = eufm10  	\@halfmag
\font\twlgot = eufm10  	\@magscale1
\font\frtngot = eufm10 	 \@magscale2
\font\svtngot = eufm10 	 \@magscale3
\font\twtygot = eufm10   \@magscale4
\newfam\bmfam

\font\fivmib = cmmib10  \@ptscale5
\font\sixmib = cmmib10  \@ptscale6
\font\sevmib = cmmib10  \@ptscale7
\font\egtmib = cmmib10  \@ptscale8
\font\tenmib = cmmib10
\font\elvmib = cmmib10   \@halfmag
\font\twlmib = cmmib10   \@magscale1
\font\frtnmib = cmmib10 	 \@magscale2
\font\svtnmib = cmmib10 	 \@magscale3
\font\twtymib = cmmib10   \@magscale4
\newdimen\indentationtitresection
\newdimen\indentationtitresubsection
\newdimen\indentationtitresubsubsection
\newdimen\indentationtitreNumero
\newdimen\indentationtitrenumero
\newdimen\indentationparag
\newdimen\indentationparagraphe

\def\@ptsize{0} \@namedef{ds@11pt}{\def\@ptsize{1}}
\@namedef{ds@12pt}{\def\@ptsize{2}}
\def\ds@twoside{\@twosidetrue \@mparswitchtrue} \def\ds@draft{\overfullrule
5pt}
\@options

\def\AA{{\bf A}}
\def\BB{{\bf B}}
\def\CC{{\bf C}}

\def\LL{{\bf L}}

\def\NN{{\bf N}}

\def\QQ{{\bf Q}}
\def\RR{{\bf R}}

\def\ZZ{{\bf Z}}

\def\cA{{\cal A}}

\def\cF{{\cal F}}

\def\cI{{\cal I}}
\def\cJ{{\cal J}}

\def\cL{{\cal L}}
\def\cM{{\cal M}}

\def\cO{{\cal O}}

\mathchardef\alphag="7C0B
\mathchardef\betag="7C0C
\mathchardef\gammag="7C0D
\mathchardef\deltag="7C0E
\mathchardef\varepsilong="7C22
\mathchardef\varphig="7C27
\mathchardef\psig="7C20
\mathchardef\zetag="7C10
\mathchardef\epsilong="7C0F
\mathchardef\rhog="7C1A
\mathchardef\taug="7C1C
\mathchardef\upsilong="7C1D
\mathchardef\iotag="7C13
\mathchardef\thetag="7C12
\mathchardef\pig="7C19
\mathchardef\sigmag="7C1B
\mathchardef\etag="7C11
\mathchardef\omegag="7C21
\mathchardef\kappag="7C14
\mathchardef\lambdag="7C15
\mathchardef\mug="7C16
\mathchardef\xig="7C18
\mathchardef\chig="7C1F
\mathchardef\nug="7C17
\mathchardef\varthetag="7C23
\mathchardef\varpig="7C24
\mathchardef\varrhog="7C25
\mathchardef\varsigmag="7C26
\mathchardef\Omegag="7C0A
\mathchardef\Thetag="7C02
\mathchardef\Sigmag="7C06
\mathchardef\Deltag="7C01
\mathchardef\Phig="7C08
\mathchardef\Gammag="7C00
\mathchardef\Psig="7C09
\mathchardef\Lambdag="7C03
\mathchardef\Xig="7C04
\mathchardef\Pig="7C05
\mathchardef\Upsilong="7C07

\def\theenumi{\arabic{enumi}}

\def\theenumii{\alph{enumii}}
\def\p@enumii{\theenumi}

\def\theenumiii{\roman{enumiii}}
\def\p@enumiii{\theenumi(\theenumii)}

\def\p@enumiv{\p@enumiii\theenumiii}

\def\verse{\let\\=\@centercr
 \list{}{\itemsep\z@ \itemindent -1.5em\listparindent \itemindent
 \rightmargin\leftmargin\advance\leftmargin 1.5em}\item[]\mbox{}\ignorespaces }

\def\quotation{\list{}{\listparindent 1.5em
 \itemindent\listparindent
 \rightmargin\leftmargin \parsep 0pt plus 1pt}\item[]\mbox{}\ignorespaces }
\let\endquotation=\endlist

\def\descriptionlabel#1{\hspace\labelsep \rm #1}
\def\description{\list{}{\labelwidth\z@ \itemindent-\leftmargin
 \let\makelabel\descriptionlabel}}

\def\titlepage{\@restonecolfalse\if@twocolumn\@restonecoltrue\onecolumn
 \else \newpage \fi \thispagestyle{empty}\c@page\z@}
\def\endtitlepage{\if@restonecol\twocolumn \else \newpage \fi}

\tabbingsep \labelsep

\skip\@mpfootins = \skip\footins
\newif\ifnumerauto

\def\numerauto{\numerautotrue\let\eqn=\@eqnauto\let\endeqn=\end@eqnauto\let\eqnn=\@eqnnauto\let\endeqnn=\end@eqnnauto\let\numeqn=\numeqnauto}
\def\numerman{\numerautofalse\let\eqn=\@eqnman \let\endeqn=\endeqnarray\let\eqnn=\eqn \let\endeqnn=\endeqn
																										\let\numeqn=\numeqnman}

\newif\ifnumeroinsection
\newif\ifnumeroinsubsection
\newif\ifnumeroinsubsubsection
\newif\ifnumerotheoreme

\def\numeroinsection{\@addtoreset{subparagraph}{section}%
																			\def\thesubparagraph{\thesection.\arabic{subparagraph}}%
																			\numeroinsectiontrue
																			\numeroinsubsectionfalse
																			\numeroinsubsubsectionfalse
																			\numerotheoremefalse}
\def\numeroinsubsection{\@addtoreset{subparagraph}{subsection}%
																			\def\thesubparagraph{\thesubsection.\arabic{subparagraph}}%
																			\numeroinsectionfalse
																			\numeroinsubsectiontrue
																			\numeroinsubsubsectionfalse
																			\numerotheoremefalse}
\def\numeroinsubsubsection{\@addtoreset{subparagraph}{subsubsection}%
																			\def\thesubparagraph{\thesubsubsection.\arabic{subparagraph}}%
																			\numeroinsectionfalse
																			\numeroinsubsectionfalse
																			\numeroinsubsubsectiontrue
																			\numerotheoremefalse}

\newif\ifcitationauto

\let\citationman=\citationautofalse

\newif\ifthmgauche

\let\thmdroite=\thmgauchefalse

\newif\ifnumthmvide

\newif\iftextefrancais

\newif\ifeqngauche
\let\eqngauche=\eqngauchetrue

\newif\ifeqninsection
\newif\ifeqninsubsection
\newif\ifeqninsubsubsection
\newif\ifeqntheoreme
\newif\ifeqnnumero

\def\@sect#1#2#3#4#5#6[#7]#8{\ifnum #2>\c@secnumdepth
     \def\@svsec{}\else
     \refstepcounter{#1}\edef\@svsec{\avantnumsections\csname the#1\endcsname\apresnumsections\hskip 1em }\fi
     \@tempskipa #5\relax
      \ifdim \@tempskipa>\z@
        \begingroup #6\relax
          \@hangfrom{\hskip #3\relax\@svsec}{\interlinepenalty \@M #8\par}
        \endgroup
       \csname #1mark\endcsname{#7}\addcontentsline
         {toc}{#1}{\ifnum #2>\c@secnumdepth \else
                      \protect\numberline{\csname the#1\endcsname}\fi
                    #7}\else
        \def\@svsechd{#6\hskip #3\@svsec #8\csname #1mark\endcsname
                      {#7}\addcontentsline
                           {toc}{#1}{\ifnum #2>\c@secnumdepth \else
                             \protect\numberline{\csname the#1\endcsname}\fi
                       #7}}\fi
     \@xsect{#5}}

\def\thepart{\Roman{part}}
\def\thesection {\arabic{section}}
\def\thesubsection {\thesection.\arabic{subsection}}
\def\thesubsubsection {\thesubsection.\arabic{subsubsection}}

\def\indentationtextechap{\let\indentationtextechappm=\@afterindenttrue}
\def\nonindentationtextechap{\let\indentationtextechappm=\@afterindentfalse}

\def\@chapappvide{}
\def\thepartvide{}

\long\def\chap#1#2{
\def\chaptest{#1}\let\newChapitre=\Chapitre
\def\thepartpleinauto{\Roman{part}}
\def\thepartpleinman{\uppercase\expandafter{\romannumeral \chaptest}}
\ifx\empty\chaptest\let\thepart=\thepartvide\let\Chapitre=\@chapappvide
\let\newespaceavanttitrechapsommaire=\chaptest
\else
	\ifnumerauto \let\thepart=\thepartpleinauto
	\else \let\thepart=\thepartpleinman
	\fi
\fi\part{#2}\let\Chapitre=\newChapitre
\ifx\empty\chaptest\addtocounter{part}{\m@ne}
\let\newespaceavanttitrechapsommaire=\espaceavanttitrechapsommaire
\fi}

\def\appendix{\par
\setcounter{part}{0}
 \setcounter{section}{0}
 \setcounter{subsection}{0}
 \def\thesection{\Alph{section}}
\let\Chapitre=\@chapappvide
\gdef\@prechapapp{\iftextefrancais Appendice\else Appendix\fi}
\gdef\thepartvide{\@prechapapp}
\long\def\chap##1##2{
\def\chaptest{##1}
\def\thepartpleinauto{\@prechapapp {} \Roman{part}}
\def\thepartpleinman{\@prechapapp {} \uppercase\expandafter{\romannumeral \chaptest}}
\ifx\empty\chaptest\let\thepart=\thepartvide
\else
	\ifnumerauto \let\thepart=\thepartpleinauto
	\else \let\thepart=\thepartpleinman
	\fi
\fi\part{##2}
\ifx\empty\chaptest\addtocounter{part}{\m@ne}\fi}
}

\long\def\sec#1#2{\def\sectiontest{#1}%
\ifnumerauto {}\else\def\thesection{\sectiontest}\fi%
\ifx\empty\sectiontest\let\sectionheadings=\sectionheadingsvide\section*{#2}\sectionmark{#2}\addcontentsline{toc}{section}{#2}%
\else\let\sectionheadings=\sectionheadingsplein\section{#2}\fi}

\long\def\seccentre#1#2{\def\sectiontest{#1}%
\ifnumerauto {}\else\def\thesection{\sectiontest}\fi%
\ifx\empty\sectiontest\let\sectionheadings=\sectionheadingsvide\sectioncentre*{#2}\sectionmark{#2}\addcontentsline{toc}{section}{#2}%
\else\let\sectionheadings=\sectionheadingsplein\sectioncentre{#2}\fi}

\long\def\subsec#1#2{\def\subsectiontest{#1}%
\ifnumerauto {}\else\def\thesubsection{\subsectiontest}\fi%
\ifx\empty\subsectiontest\let\subsectionheadings=\subsectionheadingsvide\subsection*{#2}\subsectionmark{#2}\addcontentsline{toc}{subsection}{#2}%
\else\let\subsectionheadings=\subsectionheadingsplein\subsection{#2}\fi}

\long\def\subsubsec#1#2{\def\ttest{#1}%
\ifnumerauto {}\else\def\thesubsubsection{\ttest}\fi%
\ifx\empty\ttest\subsubsection*{#2}%
\addcontentsline{toc}{subsubsection}{#2}%
\else\subsubsection{#2}\fi}

\long\def\numeroman#1#2{\def\ttest{#1}%
\def\thesubparagraph{#1}%
\def\Ttest{#2}%
\ifx\empty\ttest%
			\ifx\empty\Ttest\parag%
			\else\slsubparagraph*{#2\fintitrenumero}\mbox{}\ignorespaces \fi%
\else%
			\ifx\empty\Ttest%
\slsubparagraph*{\thesubparagraph\fintitrenumero}\mbox{}\ignorespaces %
			\else%
								\slsubparagraph{#2\fintitrenumero}\mbox{}\ignorespaces \fi
\fi}

\long\def\Numeroman#1#2{\def\ttest{#1}
\def\Ttest{#2}\def\thesubparagraph{#1}
\ifx\empty\ttest%
			\ifx\empty\Ttest\parag%
			\else\slsubsubsubsection*{\Ttest\fintitreNumero}\fi%
\else%
\ifx\empty\Ttest\slsubsubsubsection{}%
			\else\slsubsubsubsection{#2\fintitreNumero}\fi
\fi}

\long\def\numeroautomats#1#2{\def\ttest{#1}%
\def\Ttest{#2}%
\ifx\empty\ttest%
			\ifx\empty\Ttest\parag%
			\else\slsubparagraph*{#2\fintitrenumero}\mbox{}\ignorespaces \fi%
\else%
			\ifx\empty\Ttest\refstepcounter{subparagraph}%
																			\slsubparagraph*{\thesubparagraph\fintitrenumero}\mbox{}\ignorespaces %
			\else\slsubparagraph{#2\fintitrenumero}\mbox{}\ignorespaces \fi
\fi}

\long\def\Numeroautomats#1#2{\def\ttest{#1}\def\Ttest{#2}
\ifx\empty\ttest%
			\ifx\empty\Ttest\parag%
			\else\slsubsubsubsection*{\Ttest\fintitreNumero}\fi%
\else%
\ifx\empty\Ttest\slsubsubsubsection{}%
			\else\slsubsubsubsection{#2\fintitreNumero}\fi
\fi}

\long\def\numeroautomatss#1#2{\def\ttest{#1}%
\def\Ttest{#2}%
\ifx\empty\ttest%
			\ifx\empty\Ttest\parag%
			\else\slsubparagraph*{#2\fintitrenumero}\mbox{}\ignorespaces \fi%
\else%
			\ifx\empty\Ttest\refstepcounter{subparagraph}%
																			\slsubparagraph*{\thesubparagraph\fintitrenumero}\mbox{}\ignorespaces %
			\else\slsubparagraph{#2\fintitrenumero}\mbox{}\ignorespaces \fi
\fi}

\long\def\Numeroautomatss#1#2{\def\ttest{#1}\def\Ttest{#2}
\ifx\empty\ttest%
			\ifx\empty\Ttest\parag%
			\else\slsubsubsubsection*{\Ttest\fintitreNumero}\fi%
\else%
\ifx\empty\Ttest\slsubsubsubsection{}%
			\else\slsubsubsubsection{#2\fintitreNumero}\fi
\fi}

\long\def\numeroautomatsss#1#2{\def\ttest{#1}%
\def\Ttest{#2}%
\ifx\empty\ttest%
			\ifx\empty\Ttest\parag%
			\else\slsubparagraph*{#2\fintitrenumero}\mbox{}\ignorespaces \fi%
\else%
			\ifx\empty\Ttest\refstepcounter{subparagraph}%
																			\slsubparagraph*{\thesubparagraph\fintitrenumero}\mbox{}\ignorespaces %
			\else\slsubparagraph{#2\fintitrenumero}\mbox{}\ignorespaces \fi
\fi}

\long\def\Numeroautomatsss#1#2{\def\ttest{#1}\def\Ttest{#2}
\ifx\empty\ttest%
			\ifx\empty\Ttest\parag%
			\else\slsubsubsubsection*{\Ttest\fintitreNumero}\fi%
\else%
\ifx\empty\Ttest\slsubsubsubsection{}%
			\else\slsubsubsubsection{#2\fintitreNumero}\fi
\fi}

\long\def\numeroautomath#1#2{\def\ttest{#1}%
\def\Ttest{#2}%
\ifx\empty\ttest%
			\ifx\empty\Ttest\parag%
			\else\slsubparagraph*{#2\fintitrenumero}\mbox{}\ignorespaces \fi%
\else\refstepcounter{theoreme}%
			\ifx\empty\Ttest\slsubparagraph*{\thesubparagraph\fintitrenumero}\mbox{}\ignorespaces %
			\else\slsubparagraph{#2\fintitrenumero}\mbox{}\ignorespaces \fi
\fi}

\long\def\Numeroautomath#1#2{\def\ttest{#1}
\def\Ttest{#2}
\ifx\empty\ttest%
			\ifx\empty\Ttest\parag%
			\else\slsubsubsubsection*{\Ttest\fintitreNumero}\fi%
\else\refstepcounter{theoreme}%
\ifx\empty\Ttest\slsubsubsubsection{}%
			\else\slsubsubsubsection{#2\fintitreNumero}\fi
\fi}

\def\numero#1#2{%
  \ifnumerauto
     \ifnumeroinsection\numeroautomats{#1}{#2}\else
     \ifnumeroinsubsection\numeroautomatss{#1}{#2}\else
     \ifnumeroinsubsubsection\numeroautomatsss{#1}{#2}\else
     \ifnumerotheoreme\numeroautomath{#1}{#2}\else
     \numeroman{#1}{#2}\fi\fi\fi\fi
   \else
   \numeroman{#1}{#2}
   \fi}

\def\parag{\newparag{}\mbox{}\ignorespaces}

\def\rem{\iftextefrancais\numero{}{Remarque}\else
\numero{}{Remark}\fi}

\def\demo{\iftextefrancais\numero{}{D\'emonstration}\else
\numero{}{Proof}\fi}

\newdimen\indentationtitreenonce

\def\@xthm#1#2{
\ifnumerauto
\@begintheorem{#2}{\csname the#1\endcsname}
\else
\@begintheorem{#2}{\theoremetest}
\fi\ignorespaces}

\def\@ythm#1#2[#3]{
\ifnumerauto
\@opargbegintheorem{#2}{\csname the#1\endcsname}{#3}
\else
\@opargbegintheorem{#2}{\theoremetest}{#3}
\fi\ignorespaces}

\def\@thmcounter#1{\noexpand\arabic{#1}}
\def\@thmcountersep{\ifnumerauto .\else{}\fi}
\def\@makethmnumber#1#2{\bf #1 #2:}

\def\@begintheorem#1#2{%
\ifthmgauche%
						\ifnumthmvide%
									\styletexteenonce \trivlist \item[\hskip\labelsep{\mbox{\hspace{\indentationtitreenonce}}%
																								\styletitreenonce #1\fintitreenonce}]%
						\else%
									\styletexteenonce \trivlist \item[\hskip\labelsep{\mbox{\hspace{\indentationtitreenonce}}%
																								\styletitreenonce \avantnumenonce#2\apresnumenonce\ #1\fintitreenonce}]%
						\fi%
\else%
						\ifnumthmvide%
				\styletexteenonce \trivlist \item[\hskip\labelsep{\mbox{\hspace{\indentationtitreenonce}}%
																								\styletitreenonce #1\fintitreenonce}]%
						\else
				\styletexteenonce \trivlist \item[\hskip\labelsep{\mbox{\hspace{\indentationtitreenonce}}%
																								\styletitreenonce #1\ \avantnumenonce#2\apresnumenonce\fintitreenonce}]%
						\fi%
\fi}
\def\@opargbegintheorem#1#2#3{%
\ifthmgauche%
						\ifnumthmvide%
											\styletexteenonce \trivlist\item[\hskip \labelsep{\mbox{\hspace{\indentationtitreenonce}}%
																							\styletitreenonce #1\ {\styleprecisionenonce(#3)}\fintitreenonce}]%
						\else%
											\styletexteenonce \trivlist\item[\hskip \labelsep{\mbox{\hspace{\indentationtitreenonce}}%
																							\styletitreenonce \avantnumenonce#2\apresnumenonce\ #1\ {\styleprecisionenonce(#3)}\fintitreenonce}]%
						\fi
\else%
						\ifnumthmvide%
				\styletexteenonce \trivlist\item[\hskip \labelsep{\mbox{\hspace{\indentationtitreenonce}}%
																							\styletitreenonce #1\ {\styleprecisionenonce(#3)}\fintitreenonce}]%
						\else%
				\styletexteenonce \trivlist\item[\hskip \labelsep{\mbox{\hspace{\indentationtitreenonce}}%
																							\styletitreenonce #1\ \avantnumenonce#2\apresnumenonce\ {\styleprecisionenonce(#3)}\fintitreenonce}]%
						\fi
\fi}
\def\@endtheorem{\endtrivlist}

\def\newTheorem#1{\@ifnextchar[{\@oThm{#1}}{\@nThm{#1}}}

\def\@nThm#1#2{%
\@ifnextchar[{\@xnThm{#1}{#2}}{\@ynThm{#1}{#2}}}

\def\@xnThm#1#2[#3]{\expandafter\@ifdefinable\csname #1\endcsname
{\@definecounter{#1}\@addtoreset{#1}{#3}%
\expandafter\xdef\csname the#1\endcsname{\expandafter\noexpand
  \csname the#3\endcsname \@Thmcountersep \@Thmcounter{#1}}%
\global\@namedef{#1}{\@Thm{#1}{#2}}\global\@namedef{end#1}{\@endTheorem}}}

\def\@ynThm#1#2{\expandafter\@ifdefinable\csname #1\endcsname
{\@definecounter{#1}%
\expandafter\xdef\csname the#1\endcsname{\@Thmcounter{#1}}%
\global\@namedef{#1}{\@Thm{#1}{#2}}\global\@namedef{end#1}{\@endTheorem}}}

\def\@oThm#1[#2]#3{\expandafter\@ifdefinable\csname #1\endcsname
  {\global\@namedef{the#1}{\@nameuse{the#2}}%
\global\@namedef{#1}{\@Thm{#2}{#3}}%
\global\@namedef{end#1}{\@endTheorem}}}

\def\@Thm#1#2{\refstepcounter{#1}
							\@ifnextchar[{\@yThm{#1}{#2}}{\@xThm{#1}{#2}}}

\def\@xThm#1#2{
\ifnumerauto
\@beginTheorem{#2}{\csname the#1\endcsname}
\else
\@beginTheorem{#2}{\theoremetest}
\fi\ignorespaces}

\def\@yThm#1#2[#3]{
\ifnumerauto
\@opargbeginTheorem{#2}{\csname the#1\endcsname}{#3}
\else
\@opargbeginTheorem{#2}{\theoremetest}{#3}
\fi\ignorespaces}

\def\@Thmcounter#1{\noexpand\arabic{#1}}
\def\@Thmcountersep{\ifnumerauto .\else{}\fi}
\def\@makeThmnumber#1#2{\bf #1 #2:}

\def\@beginTheorem#1#2{%
\ifthmgauche%
						\ifnumthmvide%
									\styletexteenonce \trivlist \item[\hskip\labelsep{\mbox{\hspace{\indentationtitreenonce}}%
																								\styletitreenonce #1\fintitreEnonce}]\mbox{ }\samepage%
						\else
									\styletexteenonce \trivlist \item[\hskip\labelsep{\mbox{\hspace{\indentationtitreenonce}}%
																								\styletitreenonce \avantnumEnonce#2\apresnumEnonce\ #1\fintitreEnonce}]\mbox{ }\samepage%
						\fi%
\else%
						\ifnumthmvide%
				\styletexteenonce \trivlist \item[\hskip\labelsep{\mbox{\hspace{\indentationtitreenonce}}%
																								\styletitreenonce #1\fintitreEnonce}]\mbox{ }\samepage%
						\else
				\styletexteenonce \trivlist \item[\hskip\labelsep{\mbox{\hspace{\indentationtitreenonce}}%
																								\styletitreenonce #1\ \avantnumEnonce#2\apresnumEnonce\fintitreEnonce}]\mbox{ }\samepage%
						\fi%
\fi}

\def\@opargbeginTheorem#1#2#3{%
\ifthmgauche%
						\ifnumthmvide%
											\styletexteenonce \trivlist\item[\hskip \labelsep{\mbox{\hspace{\indentationtitreenonce}}%
																							\styletitreenonce #1\ {\styleprecisionenonce(#3)}\fintitreEnonce}]\mbox{ }\samepage%
						\else
											\styletexteenonce \trivlist\item[\hskip \labelsep{\mbox{\hspace{\indentationtitreenonce}}%
																							\styletitreenonce \avantnumEnonce#2\apresnumEnonce\ #1\ {\styleprecisionenonce(#3)}\fintitreEnonce}]\mbox{ }\samepage%
						\fi
\else%
						\ifnumthmvide%
				\styletexteenonce \trivlist\item[\hskip \labelsep{\mbox{\hspace{\indentationtitreenonce}}%
																							\styletitreenonce #1\ {\styleprecisionenonce(#3)}\fintitreEnonce}]\mbox{ }\samepage%
						\else
				\styletexteenonce \trivlist\item[\hskip \labelsep{\mbox{\hspace{\indentationtitreenonce}}%
																							\styletitreenonce #1\ \avantnumEnonce#2\apresnumEnonce\ {\styleprecisionenonce(#3)}\fintitreEnonce}]\mbox{ }\samepage%
						\fi
\fi}
\def\@endTheorem{\endtrivlist}

\numthmvidefalse

\newtheorem{tthm}{\ENONCE}
\newTheorem{Tthm}[tthm]{\ENONCE}

\def\enonce#1#2{\def\ENONCE{#2}\def\theoremetest{#1}%
\ifx\empty\theoremetest\numthmvidetrue\addtocounter{tthm}{-1}%
\else\numthmvidefalse\fi\begin{tthm}}
\def\endenonce{\end{tthm}\numthmvidefalse}

\def\Enonce#1#2{\def\ENONCE{#2}\def\theoremetest{#1}%
\ifx\empty\theoremetest\numthmvidetrue\addtocounter{tthm}{-1}%
\else\numthmvidefalse\fi\begin{Tthm}}
\def\endEnonce{\end{Tthm}\numthmvidefalse}

\def\th#1{\def\ENONCE{\iftextefrancais Th\'eor\`eme\else Theorem\fi}
\def\theoremetest{#1}%
\ifx\empty\theoremetest\numthmvidetrue\addtocounter{tthm}{-1}%
\else\numthmvidefalse\fi\begin{tthm}}

\def\prop#1{\def\ENONCE{Proposition}
\def\theoremetest{#1}%
\ifx\empty\theoremetest\numthmvidetrue\addtocounter{tthm}{-1}%
\else\numthmvidefalse\fi\begin{tthm}}

\def\lem#1{\def\ENONCE{\iftextefrancais Lemme\else Lemma \fi}
\def\theoremetest{#1}%
\ifx\empty\theoremetest\numthmvidetrue\addtocounter{tthm}{-1}%
\else\numthmvidefalse\fi\begin{tthm}}

\def\cor#1{\def\ENONCE{\iftextefrancais Corollaire\else Corollary\fi}
\def\theoremetest{#1}%
\ifx\empty\theoremetest\numthmvidetrue\addtocounter{tthm}{-1}%
\else\numthmvidefalse\fi\begin{tthm}}

\def\definit#1{\def\ENONCE{\iftextefrancais D\'efinition\else Definition\fi}
\def\theoremetest{#1}%
\ifx\empty\theoremetest\numthmvidetrue\addtocounter{tthm}{-1}%
\else\numthmvidefalse\fi\begin{tthm}}

\def\Th#1{\def\ENONCE{\iftextefrancais Th\'eor\`eme\else Theorem\fi}
\def\theoremetest{#1}%
\ifx\empty\theoremetest\numthmvidetrue\addtocounter{tthm}{-1}%
\else\numthmvidefalse\fi\begin{Tthm}}

\def\Prop#1{\def\ENONCE{Proposition}
\def\theoremetest{#1}%
\ifx\empty\theoremetest\numthmvidetrue\addtocounter{tthm}{-1}%
\else\numthmvidefalse\fi\begin{Tthm}}

\def\Lem#1{\def\ENONCE{\iftextefrancais Lemme\else Lemma \fi}
\def\theoremetest{#1}%
\ifx\empty\theoremetest\numthmvidetrue\addtocounter{tthm}{-1}%
\else\numthmvidefalse\fi\begin{Tthm}}

\def\Cor#1{\def\ENONCE{\iftextefrancais Corollaire\else Corollary\fi}
\def\theoremetest{#1}%
\ifx\empty\theoremetest\numthmvidetrue\addtocounter{tthm}{-1}%
\else\numthmvidefalse\fi\begin{Tthm}}

\def\Definit#1{\def\ENONCE{\iftextefrancais D\'efinition\else Definition\fi}
\def\theoremetest{#1}%
\ifx\empty\theoremetest\numthmvidetrue\addtocounter{tthm}{-1}%
\else\numthmvidefalse\fi\begin{Tthm}}

\def\@pnumwidth{1.55em} \def\@tocrmarg {2.55em}
\def\@dotsep{4.5} 

\def\tableofcontents{\section*{\titretableofcontents}\@mkboth{\titretableofcontents}{\titretableofcontents}\@starttoc{toc}}

\def\pagetableofcontents{\section*{\hfill\large{\bf
\titretableofcontents}\normalsize\hfill\mbox{}}\@mkboth{\titretableofcontents}{\titretableofcontents}\vspace{1truecm}
 \@starttoc{toc}\vspace*{\fill}\pagebreak}

\def\titresommaire{\iftextefrancais Sommaire\else Contents\fi}
\def\titretabledesmatieres{\iftextefrancais Table des mati\`eres\else Table of Contents\fi}

\def\sommaire{\let\titretableofcontents\titresommaire \tableofcontents}

\def\tabledesmatieres{\let\titretableofcontents\titretabledesmatieres \tableofcontents}

\def\pagesommaire{\let\titretableofcontents\titresommaire \pagetableofcontents}

\def\pagetabledesmatieres{\let\titretableofcontents\titretabledesmatieres \pagetableofcontents}

\def\l@part#1#2{\addpenalty{\@secpenalty}
 \addvspace{2.25em plus 1pt} \begingroup
 \@tempdima 3em \parindent \z@ \rightskip \@pnumwidth \parfillskip
-\@pnumwidth
 {\bf \leavevmode #1\hfil \hbox to\@pnumwidth{\hss #2}}\par
 \nobreak \endgroup}
\def\l@section#1#2{\addpenalty{\@secpenalty} \addvspace{1.0em plus 1pt}
\@tempdima 1.5em \begingroup
 \parindent \z@ \rightskip \@pnumwidth
 \parfillskip -\@pnumwidth
 \bf \leavevmode #1\hfil \hbox to\@pnumwidth{\hss #2}\par
 \endgroup}
\def\l@subsection{\@dottedtocline{2}{1.5em}{2.3em}}
\def\l@subsubsection{\@dottedtocline{3}{3.8em}{3.2em}}
\def\l@paragraph{\@dottedtocline{4}{7.0em}{4.1em}}
\def\l@subparagraph{\@dottedtocline{5}{10em}{5em}}

\def\listoffigures#1{\section*{\titrelistoffigures\@mkboth
 {\titrelistoffigures}{\titrelistoffigures}}\@starttoc{lof}}
\def\l@figure{\@dottedtocline{1}{1.5em}{2.3em}}

\def\listoftables#1{\section*{\titrelistoftables\@mkboth
 {\titrelistoftables}{\titrelistoftables}}\@starttoc{lot}}
\let\l@table\l@figure

\def\thebibliography#1{\section*{\titrebibliographie}%
\@mkboth{\titrebibliographie}{\titrebibliographie}%
\addcontentsline{toc}{section}{\titrebibliographie}
\list{[\arabic{enumi}]}{\settowidth\labelwidth{[#1]}%
\leftmargin\labelwidth
 \advance\leftmargin\labelsep
 \usecounter{enumi}}
 \def\newblock{\hskip .11em plus .33em minus -.07em}
 \sloppy
 \sfcode`\.=1000\relax}

\newbox\auteurbox
\newbox\titrebox
\newbox\titrelbox
\newbox\editeurbox
\newbox\anneebox
\newbox\anneelbox
\newbox\journalbox
\newbox\volumebox
\newbox\pagesbox
\newbox\diversbox
\newbox\collectionbox
\def\fabriquebox#1#2{\par\egroup
                     \setbox#1=\vbox\bgroup
                      \leftskip=0pt
                       \hsize=\maxdimen \noindent#2}
\def\bibref#1{\ifcitationauto\bibitem{#1}
  												\else\bibitem[#1]{#1}\fi\mbox{}\ignorespaces
														\setbox0=\vbox\bgroup}
\def\auteur{\fabriquebox\auteurbox\styleauteur}
	\def\titre{\fabriquebox\titrebox\styletitre}
		\def\titrelivre{\fabriquebox\titrelbox\styletitrelivre}
			\def\editeur{\fabriquebox\editeurbox\styleediteur}
				\def\journal{\fabriquebox\journalbox\stylejournal}
					\def\volume{\fabriquebox\volumebox\stylevolume}
						
{\catcode`\-=\active\gdef\annee{\fabriquebox\anneebox\catcode`\-=\active\def-{\hbox{\rm \string-\string-}}\styleannee\ignorespaces}}
{\catcode`\-=\active\gdef\anneelivre{\fabriquebox\anneelbox\catcode`\-=\active\def-{\hbox{\rm \string-\string-}}\styleanneelivre}}
{\catcode`\-=\active\gdef\pages{\fabriquebox\pagesbox\catcode`\-=\active\def-{\hbox{\rm\string-\string-}}\stylepages}}
{\catcode`\-=\active\gdef\divers{\fabriquebox\diversbox\catcode`\-=\active\def-{\hbox{\rm\string-\string-}}\rm}}
\def\ajoutref#1{\setbox0=\vbox{\unvbox#1\global\setbox1=\lastbox}\unhbox1 \unskip\unskip\unpenalty}
\newif\ifpreviousitem
\global\previousitemfalse
\def\separateur{\ifpreviousitem {,\ }\fi}
\def\voidallboxes
{\setbox0=\box\auteurbox
  \setbox0=\box\titrebox
   \setbox0=\box\titrelbox
  		\setbox0=\box\editeurbox
    \setbox0=\box\anneebox
     \setbox0=\box\anneelbox
   		\setbox0=\box\journalbox
      \setbox0=\box\volumebox
       \setbox0=\box\pagesbox
        \setbox0=\box\diversbox
          \setbox0=\box\collectionbox
           \setbox0=\null}
\def\fabriquelivre
    {\ifdim\ht\auteurbox>0pt   \ajoutref\auteurbox\global\previousitemtrue\fi
       \ifdim\ht\titrelbox>0pt  \separateur\ajoutref\titrelbox\global\previousitemtrue\fi
					\ifdim\ht\collectionbox>0pt   \separateur\ajoutref\collectionbox\global\previousitemtrue\fi
						\ifdim\ht\editeurbox>0pt  \separateur\ajoutref\editeurbox\global\previousitemtrue\fi
         			\ifdim\ht\anneelbox>0pt      \separateur \ajoutref\anneelbox  \fi\global\previousitemfalse}
\def\fabriquearticle
   {\ifdim\ht\auteurbox>0pt     \ajoutref\auteurbox \global\previousitemtrue\fi
     \ifdim\ht\titrebox>0pt  \separateur\ajoutref\titrebox\global\previousitemtrue\fi
      \ifdim\ht\titrelbox>0pt   \separateur{\rm in}\ \ajoutref\titrelbox\global\previousitemtrue\fi
				\ifdim\ht\journalbox>0pt  \separateur \ajoutref\journalbox\global\previousitemtrue\fi
       \ifdim\ht\volumebox>0pt    	 \ \ajoutref\volumebox\fi
        \ifdim\ht\anneebox>0pt      \ {\rm(}\ajoutref\anneebox \rm)\fi
         \ifdim\ht\pagesbox>0pt     \separateur\ajoutref\pagesbox\fi\global\previousitemfalse}
\def\fabriquedivers
   {\ifdim\ht\auteurbox>0pt     \ajoutref\auteurbox\global\previousitemtrue\fi
     \ifdim\ht\diversbox>0pt      \separateur\ajoutref\diversbox\fi\global\previousitemfalse}
\def\endbibref
{\egroup
    \ifdim\ht\journalbox>0pt \fabriquearticle
     \else\ifdim\ht\editeurbox>0pt \fabriquelivre
      \else\ifdim\ht\diversbox>0pt \fabriquedivers
       \fi\fi\fi
       .\voidallboxes}

\newif\if@restonecol
\def\theindex{\@restonecoltrue\if@twocolumn\@restonecolfalse\fi
\columnseprule \z@
\columnsep 35pt\twocolumn[\section*{\titreindex}]
 \@mkboth{\titreindex}{\titreindex}\thispagestyle{plain}\parindent\z@
 \parskip\z@ plus .3pt\relax\let\item\@idxitem}
\def\@idxitem{\par\hangindent 40pt}

\def\endtheindex{\if@restonecol\onecolumn\else\clearpage\fi}

\def\footnoterule{\kern-3\p@
 \hrule width .4\columnwidth
 \kern 2.6\p@}

\long\def\@makefntext#1{\parindent 1em\noindent
 \hbox to 1.8em{\hss$^{\@thefnmark}$}#1}

\long\def\@makecaption#1#2{
 \vskip 10pt
 \setbox\@tempboxa\hbox{#1 #2}
 \ifdim \wd\@tempboxa >\hsize \unhbox\@tempboxa\par \else \hbox
to\hsize{\hfil\box\@tempboxa\hfil}
 \fi}

\def\thefigure{\@arabic\c@figure}
\def\fps@figure{tbp}
\def\ftype@figure{1}
\def\ext@figure{lof}
\def\fnum@figure{Figure \thefigure\ignorespaces}
\def\figure{\@float{figure}}
\let\endfigure\end@float
\@namedef{figure*}{\@dblfloat{figure}}
\@namedef{endfigure*}{\end@dblfloat}

\def\thetable{\@arabic\c@table}
\def\fps@table{tbp}
\def\ftype@table{2}
\def\ext@table{lot}
\def\fnum@table{\iftextefrancais Tableau\else Table\fi\ \thetable\ignorespaces}
\def\table{\@float{table}}
\let\endtable\end@float
\@namedef{table*}{\@dblfloat{table}}
\@namedef{endtable*}{\end@dblfloat}

\newcounter{flotteur}
\def\theflotteur{\@arabic\c@flotteur}
\def\fps@flotteur{tbp}
\def\ftype@flotteur{1}
\def\ext@flotteur{lof}
\def\fnum@flotteur{}
\def\flotteur{\@float{flotteur}}
\let\endflotteur\end@float

\def\maketitle{\par
 \begingroup
 \def\thefootnote{\fnsymbol{footnote}}
 \def\@makefnmark{\hbox
 to 0pt{$^{\@thefnmark}$\hss}}
 \if@twocolumn
 \twocolumn[\@maketitle]
 \else \newpage
 \global\@topnum\z@ \@maketitle \fi\thispagestyle{plain}\@thanks
 \endgroup
 \setcounter{footnote}{0}
}
\def\@maketitle{\newpage
 \null
 \vskip 2em \begin{center}
 {\Large \@title \par} \vskip 1.5em {\lineskip .5em
\begin{tabular}[t]{c}\@author
 \end{tabular}\par}
 \vskip 1em {\@date} \end{center}
 \par
 \vskip 1.5em}
\def\abstract{\if@twocolumn
\section*{Abstract}
\else \small
\begin{center}
{\sc Abstract\vspace{-.5em}\vspace{0pt}}
\end{center}
\quotation
\fi}
\def\endabstract{\if@twocolumn\else\endquotation\fi}

\mark{{}{}}

\if@twoside \def\ps@headings{\let\@mkboth\markboth\def\@oddfoot{}\def\@evenfoot{}\def\@evenhead{\rm\thepage\hfil \sl \leftmark\hfil\hbox{}} \def\@oddhead{\hbox{}\hfil\sl \rightmark \hfil\rm\thepage} \def\sectionmark##1{\@mkboth{\ifnum \c@secnumdepth>\z@ \sectionheadings\relax \fi ##1}{}} \def\subsectionmark##1{\markright{\ifnum \c@secnumdepth >\@ne
 \subsectionheadings\relax \fi ##1}}}
\else \def\ps@headings{\let\@mkboth\markboth\def\@oddfoot{}\def\@evenfoot{} \def\@oddhead{\hbox{}\hfil\sl \rightmark \hfil \rm\thepage} \def\sectionmark##1{\markright{{\ifnum \c@secnumdepth >\z@
 \sectionheadings\relax \fi ##1}}}}
\fi
\def\ps@myheadings{\def\@oddhead{\hbox{}\sl\rightmark \hfil
\rm\thepage}\def\@oddfoot{}\def\@evenhead{\rm \thepage\hfil\sl\leftmark\hfil\hbox
{}}\def\@evenfoot{}\def\sectionmark##1{}\def\subsectionmark##1{}}

\def\today{\iftextefrancais\number\day\space\ifcase\month\or
janvier\or f\'evrier\or mars\or avril\or mai\or juin\or juillet\or
ao\^ut\or septembre\or octobre\or novembre\or d\'ecembre\fi
\space\number\year%
\else\ifcase\month\or
January\or February\or March\or April\or May\or June\or
July\or August\or September\or October\or November\or December\fi
\space\number\day, \number\year%
\fi}

\ps@plain  \onecolumn \if@twoside\else\fi

{\catcode`\;=\active
  \catcode`\:=\active
   \catcode`\!=\active
    \catcode`\?=\active
\gdef\frenchdactylography{\frenchspacing
      \catcode`\;=\active
       \catcode`\:=\active
        \catcode`\!=\active
         \catcode`\?=\active
\def;{\relax\ifhmode\ifdim\lastskip>\z@
       \unskip\fi\kern\fontdimen2 \font\kern -1.2 \fontdimen3 \font\fi\string;}%
\def:{\relax\ifhmode\ifdim\lastskip>\z@\unskip\fi\penalty\@M\ \fi\string:}%
\def!{\relax\ifhmode\ifdim\lastskip>\z@
       \unskip\fi\kern\fontdimen2 \font \kern -1.1 \fontdimen3 \font\fi\string!}%
\def?{\relax\ifhmode\ifdim\lastskip>\z@
       \unskip\fi\kern\fontdimen2 \font \kern -1.1 \fontdimen3
\font\fi\string?}%
}}

{\catcode`\;=\active
  \catcode`\:=\active
   \catcode`\!=\active
    \catcode`\?=\active
\gdef\englishdactylography{\nonfrenchspacing
      \catcode`\;=\active
       \catcode`\:=\active
        \catcode`\!=\active
         \catcode`\?=\active
\def;{\string;}%
\def:{\string:}%
\def!{\string!}%
\def?{\string?}%
}}

\def\textefrancais{\textefrancaistrue\frenchdactylography}
\def\texteanglais{\textefrancaisfalse\englishdactylography}

\def\pointir{\unskip . --- \ignorespaces}
\def\pointt{\unskip .\ignorespaces}

\def\ie{{\it i.e}}

\def\T{\S\kern .15em }
\def\qed{\mbox{$\Box$}}

\def\@eqnauto{\def\@eqncr{\nonumber\@seqncr}\eqnarray}
\def\end@eqnauto{\nonumber\endeqnarray}

\def\theequation{\arabic{equation}}

\def\@eqnnumd{\hbox{\rm (\theequation)}}%
\def\@eqnnumg{\hbox to .01pt{}\rlap{\rm \hskip -\displaywidth(\theequation)}}%

\def\@eqnnummand{\hbox{\rm \theequation}}%
\def\@eqnnummang{\hbox to .01pt{}\rlap{\rm \hskip -\displaywidth \theequation}}%

\def\@eqnman{\def\eqnstyle{}\eqnarray\global\def\@eqnnum{}%
\global\def\theequation{\eqnstyle}}

\def\@eqnnautos{\@addtoreset{equation}{section}%
																\def\eqnstyle{\thesection.\arabic{equation}}}

\def\@eqnnautoss{\@addtoreset{equation}{subsection}%
																\def\eqnstyle{\thesubsection.\arabic{equation}}}

\def\@eqnnautosss{\@addtoreset{equation}{subsubsection}%
																	\def\eqnstyle{\thesubsubsection.\arabic{equation}}}

\def\@eqnnautoth{\let\c@equation=\c@theoreme%
																\def\eqnstyle{\thetheoreme}}

\def\@eqnnauton{\let\c@equation=\c@subparagraph
															\def\eqnstyle{\thesubparagraph}}

\def\@eqnnautoe{%
															\def\eqnstyle{\arabic{equation}}}

\def\@eqnnauto{\ifeqninsection\@eqnnautos
														\else
																		\ifeqninsubsection\@eqnnautoss
																		\else
																					\ifeqninsubsubsection\@eqnnautosss
																					\else
																									\ifeqntheoreme\@eqnnautoth
																									\else
																													\ifeqnnumero\@eqnnauton
																													\else \@eqnnautoe
																													\fi
																									\fi
																					\fi
																		\fi
														\fi
\global\def\theequation{\eqnstyle}\eqnarray\global\def\@eqnnum{\ifeqngauche\@eqnnumg\else \@eqnnumd\fi}\global\def\theequation{\eqnstyle}}

\def\end@eqnnauto{\endeqnarray\global\def\theequation{\eqnstyle}}

\def\numeqnman#1{\global\def\@eqnnum{\ifeqngauche\@eqnnummang
																																		\else\@eqnnummand\fi}%
								\def\\{\@eqncr\global\def\theequation{\eqnstyle}
														\global\def\@eqnnum{}}
								\global\def\theequation{#1}
								\ifx\empty\theequation\nonumber
								\else\global\advance\c@equation\m@ne
								\fi}

\def\numeqnauto#1{\global\def\@eqnnum{\ifeqngauche\@eqnnummang
																																		\else\@eqnnummand\fi}%
																\def\\{\@eqncr\global\def\theequation{\eqnstyle}
																						\global\def\@eqnnum{\ifeqngauche\@eqnnumg
																																								\else\@eqnnumd
																																								\fi}}
																\global\def\theequation{#1}
																\ifx\empty\theequation\nonumber
																\else\global\advance\c@equation\m@ne
																\fi}

\newdimen\lengtharrow
 \newbox\exponantbox \newbox\indicebox

\def\dimmax#1#2{\ifdim#1<#2 #2\else #1\fi}

\def\arrowr#1#2%
    {\setbox\exponantbox=\hbox{$#2$}
      \setbox\indicebox=\hbox{$#1$}
       \lengtharrow=\dimmax{\wd\indicebox}{\wd\exponantbox}
        \ifdim \lengtharrow=0pt
															\mathrel{\smash{\mathop{\hbox to 10mm{\rightarrowfill}}}}
         \else \ifdim\lengtharrow<6mm
                \lengtharrow=10mm
                 \else \advance\lengtharrow by 4mm
                  \fi
					\mathrel{\smash{\mathop{\hbox to\lengtharrow{\rightarrowfill}}\limits%
                         ^{\styleflechehoriz#2}_{\styleflechehoriz#1}}}
									\fi}

\def\arrowl#1#2%
    {\setbox\exponantbox=\hbox{$#2$}
      \setbox\indicebox=\hbox{$#1$}
       \lengtharrow=\dimmax{\wd\indicebox}{\wd\exponantbox}
        \ifdim \lengtharrow=0pt
															\mathrel{\smash{\mathop{\hbox to 10mm{\leftarrowfill}}}}
         \else \ifdim\lengtharrow<6mm
                \lengtharrow=10mm
                 \else \advance\lengtharrow by 4mm
                  \fi
					\mathrel{\smash{\mathop{\hbox to\lengtharrow{\leftarrowfill}}\limits%
                         ^{\styleflechehoriz#2}_{\styleflechehoriz#1}}}
									\fi}

\newbox\resumebox
\newbox\abstractbox
\newbox\classbox
\newbox\motsclesbox
\newbox\keywordsbox
\newbox\Titrebox
\newbox\Titredeuxbox
\newbox\Titretroisbox
\newbox\Auteurbox
\newbox\Auteurdeuxbox
\newbox\Auteurtroisbox
\newbox\Auteurquatrebox
\newbox\Auteurcinqbox
\newbox\numeroseriebox

\newdimen\lengthresume
\newdimen\lengthclasmots
\newdimen\lengthboiteresume

\def\pagedetitre{\setcounter{page}{\m@ne}
														\setbox0=\vbox\bgroup}

\def\Date{\normalsize\bf\ifcase\month\or
janvier\or f\'evrier\or mars\or avril\or mai\or juin\or juillet\or
ao\^ut\or septembre\or octobre\or novembre\or d\'ecembre\fi
\space\number\year}

\def\endpagedetitre{\egroup
							\lengthresume=\dimmax{\ht\resumebox}{\ht\abstractbox}
							\lengthclasmots=\dimmax{\ht\classbox}{\ht\motsclesbox}
							\lengthboiteresume=\dimmax{\lengthresume}{\lengthclasmots}
																	\fabriquepagedetitre}

\def\fabriquepagedetitre{\par\noindent\mbox{}\hfill\vbox to 13truecm \bgroup
                      												 \hsize=13truecm \noindent\vfill\input psfig \pssilent\par\centerline{\psfig{file=\lieulogo prepub.ps,height=8.2truecm}}\par
																							\vspace*{4.5truecm}\vfill\par\egroup\hfill\mbox{}
																				\par\noindent\mbox{}\hfill\vbox to 7.5truecm \bgroup
                     												  \hsize=12truecm \noindent\boitefenetre
																							\vfill\par\egroup\hfill\mbox{}
																				\par\noindent\mbox{}\hfill\vbox to 1.7truecm \bgroup
         												              \hsize=15truecm \noindent\vfill\boiteadresse
																							\par\egroup\hfill\mbox{}\par
																				\thispagestyle{empty}
																				\vfill
																				\pagebreak
																				\ifdim\lengthboiteresume=0pt
																						\addtocounter{page}{1}
																				\else
																						\begin{quotation}\mbox{}\small
																						\boiteresume
																						\end{quotation}\normalsize\rm
																						\thispagestyle{empty}
																						\vfill
																						\pagebreak
																				\fi}

\def\boiteresume{\ifdim\ht\resumebox>0pt
																									\scparagraph{R\'esum\'e\pointir}
																									\ajoutref\resumebox\fi
														\ifdim\ht\abstractbox>0pt
																									\scparagraph{Abstract\pointir}
																									\ajoutref\abstractbox\fi
														\ifdim\ht\classbox>0pt
																									\scparagraph{Classification AMS\pointir}
																									\ajoutref\classbox\fi
														\ifdim\ht\motsclesbox>0pt
																									\scparagraph{Mots-cl\'es\pointir}
																									\ajoutref\motsclesbox\fi
														\ifdim\ht\keywordsbox>0pt
																									\scparagraph{Keywords\pointir}
																									\ajoutref\keywordsbox\fi
														\normalsize\rm}

\def\boitefenetre{{\baselineskip 18pt\vskip .5truecm
															\begin{center}
															\ifdim\ht\Titrebox>0pt
																		\ajoutref\Titrebox\fi
															\ifdim\ht\Titredeuxbox>0pt
																		\\ \ajoutref\Titredeuxbox\fi
															\ifdim\ht\Titretroisbox>0pt
																		\\ \ajoutref\Titretroisbox\fi
															\\ \vskip 1truecm
															\begin{tabular}[t]{c}
															\ifdim\ht\Auteurbox>0pt
																		\ajoutref\Auteurbox\fi
															\ifdim\ht\Auteurdeuxbox>0pt
																		\and \ajoutref\Auteurdeuxbox\fi
															\ifdim\ht\Auteurtroisbox>0pt
																		\and \ajoutref\Auteurtroisbox\fi
															\ifdim\ht\Auteurquatrebox>0pt
																		\and \ajoutref\Auteurquatrebox\fi
															\ifdim\ht\Auteurcinqbox>0pt
																		\and \ajoutref\Auteurcinqbox\fi
															\end{tabular}\end{center}
															\par\mbox{}\vfill\mbox{}
															\par\noindent\hspace*{1truecm}
															\ifdim\ht\numeroseriebox>0pt
																		{\bf N}${}^{\bf o}$\kern .15em\ajoutref\numeroseriebox\fi
															\hbox{}\hfill\Date\hspace{1truecm}\mbox{}
															}\normalsize\rm}

\def\boiteadresse{\begin{small}\begin{center}
F--91128 Palaiseau Cedex\\
T\'el.~: $(({\mit33}))\;({\mit1})\; {\mit 69}\; {\mit 33}\; {\mit 40}\; {\mit 88}\quad\bullet\quad$Fax~: $(({\mit33}))\;({\mit1})\; {\mit 69}\; {\mit 33}\; {\mit 30}\; {\mit 19}$\\
Internet~: secret@orphee.polytechnique.fr
\end{center}\end{small}}

\texteanglais
\long\def\ig#1{\relax}
\ig{Thanks to Roberto Minio for this def'n.  Compare the def'n of
\comment in AMSTeX.}

\newcount \c@coefa
\newcount \c@coefb
\newcount \c@coefc
\newdimen\tempdimen
\newdimen\xlen
\newdimen\ylen

\newcount\tempcounta
\newcount\tempcountb
\newcount\tempcountc
\newcount\tempcountd
\newcount\tempcounte
\newcount\tempcountf
\newcount\xext
\newcount\yext
\newcount\xoff
\newcount\yoff
\newcount\c@gap%
\newcount\c@arrowtypea
\newcount\c@arrowtypeb
\newcount\c@arrowtypec
\newcount\c@arrowtyped
\newcount\c@arrowtypee
\newcount\c@height
\newcount\c@width
\newcount\xpos
\newcount\ypos
\newcount\c@run
\newcount\c@rise
\newcount\c@arrowlength
\newcount\c@halflength
\newcount\c@arrowtype
\newsavebox{\tempboxa}%
\newsavebox{\tempboxb}%
\newsavebox{\tempboxc}%

\def\settoheight#1#2{\setbox\@tempboxa\hbox{#2}#1\ht\@tempboxa\relax}%
\def\settodepth#1#2{\setbox\@tempboxa\hbox{#2}#1\dp\@tempboxa\relax}%

\def\putbox(#1,#2)#3{\put(#1,#2){\makebox(0,0){#3}}}

\def\setsqparms[#1`#2`#3`#4;#5`#6]{%
\settripairparms[#1`#2`#3`#4`1;#6]%
\c@width #5
}

\def\settriparms[#1`#2`#3;#4]{\settripairparms[#1`#2`#3`1`1;#4]}%

\def\settripairparms[#1`#2`#3`#4`#5;#6]{%
\c@arrowtypea #1
\c@arrowtypeb #2
\c@arrowtypec #3
\c@arrowtyped #4
\c@arrowtypee #5
\c@height #6
\c@width #6
}

\def\resetparms{\settripairparms[1`1`1`1`1;500]\c@width 500}

\def\mvector(#1,#2)#3{
\put(0,0){\vector(#1,#2){#3}}%
\put(0,0){\vector(#1,#2){30}}%
}
\def\evector(#1,#2)#3{{
\c@arrowlength #3
\put(0,0){\vector(#1,#2){\c@arrowlength}}%
\advance \c@arrowlength by-30
\put(0,0){\vector(#1,#2){\c@arrowlength}}%
}}

\def\horsize#1#2{%
\settowidth{\tempdimen}{$#2$}%
#1=\tempdimen
\divide #1 by\unitlength
}

\def\vertsize#1#2{%
\settoheight{\tempdimen}{$#2$}%
#1=\tempdimen
\settodepth{\tempdimen}{$#2$}%
\advance #1 by\tempdimen
\divide #1 by\unitlength
}

\def\vertadjust[#1`#2`#3]{%
\vertsize{\tempcounta}{#1}%
\vertsize{\tempcountb}{#2}%
\ifnum \tempcounta<\tempcountb \tempcounta=\tempcountb \fi
\divide\tempcounta by2
\vertsize{\tempcountb}{#3}%
\ifnum \tempcountb>0 \advance \tempcountb by20 \fi
\ifnum \tempcounta<\tempcountb \tempcounta=\tempcountb \fi
}

\def\horadjust[#1`#2`#3]{%
\horsize{\tempcounta}{#1}%
\horsize{\tempcountb}{#2}%
\ifnum \tempcounta<\tempcountb \tempcounta=\tempcountb \fi
\divide\tempcounta by20
\horsize{\tempcountb}{#3}%
\ifnum \tempcountb>0 \advance \tempcountb by60 \fi
\ifnum \tempcounta<\tempcountb \tempcounta=\tempcountb \fi
}

\ig{ In this procedure, #1 is the paramater that sticks out all the way,
#2 sticks out the least and #3 is a label sticking out half way.  #4 is
the amount of the offset.}

\def\sladjust[#1`#2`#3]#4{%
\tempcountc=#4
\horsize{\tempcounta}{#1}%
\divide \tempcounta by2
\horsize{\tempcountb}{#2}%
\divide \tempcountb by2
\advance \tempcountb by-\tempcountc
\ifnum \tempcounta<\tempcountb \tempcounta=\tempcountb\fi
\divide \tempcountc by2
\horsize{\tempcountb}{#3}%
\advance \tempcountb by-\tempcountc
\ifnum \tempcountb>0 \advance \tempcountb by80\fi
\ifnum \tempcounta<\tempcountb \tempcounta=\tempcountb\fi
\advance\tempcounta by20
}

\def\putvector(#1,#2)(#3,#4)#5#6{{%
\xpos=#1
\ypos=#2
\c@run=#3
\c@rise=#4
\c@arrowlength=#5
\c@arrowtype=#6
\ifnum \c@arrowtype<0
    \ifnum \c@run=0
        \advance \ypos by-\c@arrowlength
    \else
        \tempcounta \c@arrowlength
        \multiply \tempcounta by\c@rise
        \divide \tempcounta by\c@run
        \ifnum\c@run>0
            \advance \xpos by\c@arrowlength
            \advance \ypos by\tempcounta
        \else
            \advance \xpos by-\c@arrowlength
            \advance \ypos by-\tempcounta
        \fi
    \fi
    \multiply \c@arrowtype by-1
    \multiply \c@rise by-1
    \multiply \c@run by-1
\fi
\ifnum \c@arrowtype=1
    \put(\xpos,\ypos){\vector(\c@run,\c@rise){\c@arrowlength}}%
\else\ifnum \c@arrowtype=2
    \put(\xpos,\ypos){\mvector(\c@run,\c@rise)\c@arrowlength}%
\else\ifnum\c@arrowtype=3
    \put(\xpos,\ypos){\evector(\c@run,\c@rise){\c@arrowlength}}%
\fi\fi\fi
}}

\def\bfig{\begin{picture}(\xext,\yext)(\xoff,\yoff)}
\def\efig{\end{picture}}

\def\putsplitvector(#1,#2)#3#4{
\xpos #1
\ypos #2
\c@arrowtype #4
\c@halflength #3
\c@arrowlength #3
\c@gap 140
\advance \c@halflength by-\c@gap
\divide \c@halflength by2
\ifnum \c@arrowtype=1
    \put(\xpos,\ypos){\line(0,-1){\c@halflength}}%
    \advance\ypos by-\c@halflength
    \advance\ypos by-\c@gap
    \put(\xpos,\ypos){\vector(0,-1){\c@halflength}}%
\else\ifnum \c@arrowtype=2
    \put(\xpos,\ypos){\line(0,-1)\c@halflength}%
    \put(\xpos,\ypos){\vector(0,-1)3}%
    \advance\ypos by-\c@halflength
    \advance\ypos by-\c@gap
    \put(\xpos,\ypos){\vector(0,-1){\c@halflength}}%
\else\ifnum\c@arrowtype=3
    \put(\xpos,\ypos){\line(0,-1)\c@halflength}%
    \advance\ypos by-\c@halflength
    \advance\ypos by-\c@gap
    \put(\xpos,\ypos){\evector(0,-1){\c@halflength}}%
\else\ifnum \c@arrowtype=-1
    \advance \ypos by-\c@arrowlength
    \put(\xpos,\ypos){\line(0,1){\c@halflength}}%
    \advance\ypos by\c@halflength
    \advance\ypos by\c@gap
    \put(\xpos,\ypos){\vector(0,1){\c@halflength}}%
\else\ifnum \c@arrowtype=-2
    \advance \ypos by-\c@arrowlength
    \put(\xpos,\ypos){\line(0,1)\c@halflength}%
    \put(\xpos,\ypos){\vector(0,1)3}%
    \advance\ypos by\c@halflength
    \advance\ypos by\c@gap
    \put(\xpos,\ypos){\vector(0,1){\c@halflength}}%
\else\ifnum\c@arrowtype=-3
    \advance \ypos by-\c@arrowlength
    \put(\xpos,\ypos){\line(0,1)\c@halflength}%
    \advance\ypos by\c@halflength
    \advance\ypos by\c@gap
    \put(\xpos,\ypos){\evector(0,1){\c@halflength}}%
\fi\fi\fi\fi\fi\fi
}

\def\setpos(#1,#2){\xpos=#1 \ypos#2}

\def\putmorphism(#1)(#2,#3)[#4`#5`#6]#7#8#9{{%
\c@run #2
\c@rise #3
\ifnum\c@rise=0
  \puthmorphism(#1)[#4`#5`#6]{#7}{#8}{#9}%
\else\ifnum\c@run=0
  \putvmorphism(#1)[#4`#5`#6]{#7}{#8}{#9}%
\else
\setpos(#1)%
\c@arrowlength #7
\c@arrowtype #8
\ifnum\c@run=0
\else\ifnum\c@rise=0
\else
\ifnum\c@run>0
    \c@coefa=1
\else
   \c@coefa=-1
\fi
\ifnum\c@arrowtype>0
   \c@coefb=0
   \c@coefc=-1
\else
   \c@coefb=\c@coefa
   \c@coefc=1
   \c@arrowtype=-\c@arrowtype
\fi
\c@width=2
\multiply \c@width by\c@run
\divide \c@width by\c@rise
\ifnum \c@width<0  \c@width=-\c@width\fi
\advance\c@width by60
\if l#9 \c@width=-\c@width\fi
\putbox(\xpos,\ypos){$#4$}
{\multiply \c@coefa by\c@arrowlength
\advance\xpos by\c@coefa
\multiply \c@coefa by\c@rise
\divide \c@coefa by\c@run
\advance \ypos by\c@coefa
\putbox(\xpos,\ypos){$#5$} }%
{\multiply \c@coefa by\c@arrowlength
\divide \c@coefa by2
\advance \xpos by\c@coefa
\advance \xpos by\c@width
\multiply \c@coefa by\c@rise
\divide \c@coefa by\c@run
\advance \ypos by\c@coefa
\if l#9%
   \put(\xpos,\ypos){\makebox(0,0)[r]{$#6$}}%
\else\if r#9%
   \put(\xpos,\ypos){\makebox(0,0)[l]{$#6$}}%
\fi\fi }%
{\multiply \c@rise by-\c@coefc
\multiply \c@run by-\c@coefc
\multiply \c@coefb by\c@arrowlength
\advance \xpos by\c@coefb
\multiply \c@coefb by\c@rise
\divide \c@coefb by\c@run
\advance \ypos by\c@coefb
\multiply \c@coefc by70
\advance \ypos by\c@coefc
\multiply \c@coefc by\c@run
\divide \c@coefc by\c@rise
\advance \xpos by\c@coefc
\multiply \c@coefa by140
\multiply \c@coefa by\c@run
\divide \c@coefa by\c@rise
\advance \c@arrowlength by\c@coefa
\ifnum \c@arrowtype=1
   \put(\xpos,\ypos){\vector(\c@run,\c@rise){\c@arrowlength}}%
\else\ifnum\c@arrowtype=2
   \put(\xpos,\ypos){\mvector(\c@run,\c@rise){\c@arrowlength}}%
\else\ifnum\c@arrowtype=3
   \put(\xpos,\ypos){\evector(\c@run,\c@rise){\c@arrowlength}}%
\fi\fi\fi}%
\fi\fi
\fi\fi}}

\def\puthmorphism(#1,#2)[#3`#4`#5]#6#7#8{{%
\xpos #1
\ypos #2
\c@width #6
\c@arrowlength #6
\putbox(\xpos,\ypos){$#3$\vphantom{$#4$}}%
{\advance \xpos by\c@arrowlength
\putbox(\xpos,\ypos){\vphantom{$#3$}$#4$}}%
\horsize{\tempcounta}{#3}%
\horsize{\tempcountb}{#4}%
\divide \tempcounta by2
\divide \tempcountb by2
\advance \tempcounta by30
\advance \tempcountb by30
\advance \xpos by\tempcounta
\advance \c@arrowlength by-\tempcounta
\advance \c@arrowlength by-\tempcountb
\putvector(\xpos,\ypos)(1,0){\c@arrowlength}{#7}%
\divide \c@arrowlength by2
\advance \xpos by\c@arrowlength
\vertsize{\tempcounta}{#5}%
\divide\tempcounta by2
\advance \tempcounta by20
\if a#8 %
   \advance \ypos by\tempcounta
   \put(\xpos,\ypos){\makebox(0,0){$#5$}}%
\else
   \advance \ypos by-\tempcounta
   \put(\xpos,\ypos){\makebox(0,0){$#5$}}%
\fi
}}

\def\putvmorphism(#1,#2)[#3`#4`#5]#6#7#8{{%
\xpos #1
\ypos #2
\c@arrowlength #6
\c@arrowtype #7
\settowidth{\xlen}{$#5$}%
\putbox(\xpos,\ypos){$#3$}%
{\advance \ypos by-\c@arrowlength
\putbox(\xpos,\ypos){$#4$}}%
{\advance\c@arrowlength by-140
\advance \ypos by-70
\ifdim\xlen>0pt
   \if m#8%
      \putsplitvector(\xpos,\ypos){\c@arrowlength}{\c@arrowtype}%
   \else
      \putvector(\xpos,\ypos)(0,-1){\c@arrowlength}{\c@arrowtype}%
   \fi
\else
   \putvector(\xpos,\ypos)(0,-1){\c@arrowlength}{\c@arrowtype}%
\fi}%
\ifdim\xlen>0pt
   \divide \c@arrowlength by2
   \advance\ypos by-\c@arrowlength
   \if l#8%
      \advance \xpos by-40
      \put(\xpos,\ypos){\makebox(0,0)[r]{$#5$}}%
   \else\if r#8%
      \advance \xpos by40
      \put(\xpos,\ypos){\makebox(0,0)[l]{$#5$}}%
   \else
      \putbox(\xpos,\ypos){$#5$}%
   \fi\fi
\fi
}}

\def\topadjust[#1`#2`#3]{%
\yoff=10
\vertadjust[#1`#2`{#3}]%
\advance \yext by\tempcounta
\advance \yext by 10
}

\def\botadjust[#1`#2`#3]{%
\vertadjust[#1`#2`{#3}]%
\advance \yext by\tempcounta
\advance \yoff by-\tempcounta
}

\def\leftadjust[#1`#2`#3]{%
\xoff=0
\horadjust[#1`#2`{#3}]%
\advance \xext by\tempcounta
\advance \xoff by-\tempcounta
}

\def\rightadjust[#1`#2`#3]{%
\horadjust[#1`#2`{#3}]%
\advance \xext by\tempcounta
}

\def\rightsladjust[#1`#2`#3]{%
\sladjust[#1`#2`{#3}]{\c@width}%
\advance \xext by\tempcounta
}

\def\leftsladjust[#1`#2`#3]{%
\xoff=0
\sladjust[#1`#2`{#3}]{\c@width}%
\advance \xext by\tempcounta
\advance \xoff by-\tempcounta
}

\def\adjust[#1`#2;#3`#4;#5`#6;#7`#8]{%
\topadjust[#1``{#2}]
\leftadjust[#3``{#4}]
\rightadjust[#5``{#6}]
\botadjust[#7``{#8}]}

\def\putsquare(#1)[#2`#3`#4`#5;#6`#7`#8`#9]{%
\setpos(#1)
\puthmorphism(\xpos,\ypos)[#4`#5`{#9}]{\c@width}{\c@arrowtyped}b%
\advance\ypos by \c@height
\puthmorphism(\xpos,\ypos)[#2`#3`{#6}]{\c@width}{\c@arrowtypea}a%
\putvmorphism(\xpos,\ypos)[``{#7}]{\c@height}{\c@arrowtypeb}l%
\advance\xpos by \c@width
\putvmorphism(\xpos,\ypos)[``{#8}]{\c@height}{\c@arrowtypec}r%
}

\def\square[#1`#2`#3`#4;#5`#6`#7`#8]{{
\xext=\c@width                              
\yext=\c@height                             
\topadjust[#1`#2`{#5}]
\botadjust[#3`#4`{#8}]
\leftadjust[#1`#3`{#6}]
\rightadjust[#2`#4`{#7}]
\begin{picture}(\xext,\yext)(\xoff,\yoff)
\putsquare(0,0)[#1`#2`#3`#4;#5`#6`#7`{#8}]
\end{picture}
}}

\def\putptriangle(#1,#2)[#3`#4`#5;#6`#7`#8]{%
\xpos=#1 \ypos=#2
\advance\ypos by \c@height
\puthmorphism(\xpos,\ypos)[#3`#4`{#6}]{\c@height}{\c@arrowtypea}a%
\putvmorphism(\xpos,\ypos)[`#5`{#7}]{\c@height}{\c@arrowtypeb}l%
\advance\xpos by\c@height
\putmorphism(\xpos,\ypos)(-1,-1)[``{#8}]{\c@height}{\c@arrowtypec}r%
}

\def\ptriangle[#1`#2`#3;#4`#5`#6]{{
\c@width=\c@height                         
\xext=\c@width                           
\yext=\c@width                           
\topadjust[#1`#2`{#4}]
\botadjust[#3``]
\leftadjust[#1`#3`{#5}]
\rightsladjust[#2`#3`{#6}]
\begin{picture}(\xext,\yext)(\xoff,\yoff)
\putptriangle(0,0)[#1`#2`#3;#4`#5`{#6}]%
\end{picture}%
}}

\def\putqtriangle(#1,#2)[#3`#4`#5;#6`#7`#8]{%
\xpos=#1 \ypos=#2
\advance\ypos by\c@height
\puthmorphism(\xpos,\ypos)[#3`#4`{#6}]{\c@height}{\c@arrowtypea}a%
\putmorphism(\xpos,\ypos)(1,-1)[``{#7}]{\c@height}{\c@arrowtypeb}l%
\advance\xpos by\c@height
\putvmorphism(\xpos,\ypos)[`#5`{#8}]{\c@height}{\c@arrowtypec}r%
}

\def\qtriangle[#1`#2`#3;#4`#5`#6]{{
\c@width=\c@height                         
\xext=\c@width                           
\yext=\c@height                          
\topadjust[#1`#2`{#4}]
\botadjust[#3``]
\leftsladjust[#1`#3`{#5}]
\rightadjust[#2`#3`{#6}]
\begin{picture}(\xext,\yext)(\xoff,\yoff)
\putqtriangle(0,0)[#1`#2`#3;#4`#5`{#6}]%
\end{picture}%
}}

\def\putdtriangle(#1,#2)[#3`#4`#5;#6`#7`#8]{%
\xpos=#1 \ypos=#2
\puthmorphism(\xpos,\ypos)[#4`#5`{#8}]{\c@height}{\c@arrowtypec}b%
\advance\xpos by \c@height \advance\ypos by\c@height
\putmorphism(\xpos,\ypos)(-1,-1)[``{#6}]{\c@height}{\c@arrowtypea}l%
\putvmorphism(\xpos,\ypos)[#3``{#7}]{\c@height}{\c@arrowtypeb}r%
}

\def\dtriangle[#1`#2`#3;#4`#5`#6]{{
\c@width=\c@height                         
\xext=\c@width                           
\yext=\c@height                          
\topadjust[#1``]
\botadjust[#2`#3`{#6}]
\leftsladjust[#2`#1`{#4}]
\rightadjust[#1`#3`{#5}]
\begin{picture}(\xext,\yext)(\xoff,\yoff)
\putdtriangle(0,0)[#1`#2`#3;#4`#5`{#6}]%
\end{picture}%
}}

\def\putbtriangle(#1,#2)[#3`#4`#5;#6`#7`#8]{%
\xpos=#1 \ypos=#2
\puthmorphism(\xpos,\ypos)[#4`#5`{#8}]{\c@height}{\c@arrowtypec}b%
\advance\ypos by\c@height
\putmorphism(\xpos,\ypos)(1,-1)[``{#7}]{\c@height}{\c@arrowtypeb}r%
\putvmorphism(\xpos,\ypos)[#3``{#6}]{\c@height}{\c@arrowtypea}l%
}

\def\btriangle[#1`#2`#3;#4`#5`#6]{{
\c@width=\c@height                         
\xext=\c@width                           
\yext=\c@height                          
\topadjust[#1``]
\botadjust[#2`#3`{#6}]
\leftadjust[#1`#2`{#4}]
\rightsladjust[#3`#1`{#5}]
\begin{picture}(\xext,\yext)(\xoff,\yoff)
\putbtriangle(0,0)[#1`#2`#3;#4`#5`{#6}]%
\end{picture}%
}}

\def\putAtriangle(#1,#2)[#3`#4`#5;#6`#7`#8]{%
\xpos=#1 \ypos=#2
{\multiply \c@height by2
\puthmorphism(\xpos,\ypos)[#4`#5`{#8}]{\c@height}{\c@arrowtypec}b}%
\advance\xpos by\c@height \advance\ypos by\c@height
\putmorphism(\xpos,\ypos)(-1,-1)[#3``{#6}]{\c@height}{\c@arrowtypea}l%
\putmorphism(\xpos,\ypos)(1,-1)[``{#7}]{\c@height}{\c@arrowtypeb}r%
}

\def\Atriangle[#1`#2`#3;#4`#5`#6]{{
\c@width=\c@height                         
\xext=\c@width                           
\yext=\c@height                          
\topadjust[#1``]
\botadjust[#2`#3`{#6}]
\multiply \xext by2 
\leftsladjust[#2`#1`{#4}]
\rightsladjust[#3`#1`{#5}]
\begin{picture}(\xext,\yext)(\xoff,\yoff)%
\putAtriangle(0,0)[#1`#2`#3;#4`#5`{#6}]%
\end{picture}%
}}

\def\putAtrianglepair(#1,#2)[#3]{\xpos=#1 \ypos=#2%
\putAtrianglepairx[#3]}
\def\putAtrianglepairx[#1`#2`#3`#4;#5`#6`#7`#8`#9]{%
\puthmorphism(\xpos,\ypos)[#2`#3`{#8}]{\c@height}{\c@arrowtyped}b%
\advance\xpos by\c@height
\puthmorphism(\xpos,\ypos)[\phantom{#3}`#4`{#9}]{\c@height}{\c@arrowtypee}b%
\advance\ypos by\c@height
\putmorphism(\xpos,\ypos)(-1,-1)[#1``{#5}]{\c@height}{\c@arrowtypea}l%
\putvmorphism(\xpos,\ypos)[``{#6}]{\c@height}{\c@arrowtypeb}m%
\putmorphism(\xpos,\ypos)(1,-1)[``{#7}]{\c@height}{\c@arrowtypec}r%
}

\def\Atrianglepair[#1`#2`#3`#4;#5`#6`#7`#8`#9]{{%
\c@width=\c@height
\xext=\c@width
\yext=\c@height
\topadjust[#1``]%
\vertadjust[#2`#3`{#8}]
\tempcountd=\tempcounta                       
\vertadjust[#3`#4`{#9}]
\ifnum\tempcounta<\tempcountd                 
\tempcounta=\tempcountd\fi                    
\advance \yext by\tempcounta                  
\advance \yoff by-\tempcounta                 
\multiply \xext by2 
\leftsladjust[#2`#1`{#5}]
\rightsladjust[#4`#1`{#7}]%
\begin{picture}(\xext,\yext)(\xoff,\yoff)%
\putAtrianglepair(0,0)[#1`#2`#3`#4;#5`#6`#7`#8`{#9}]%
\end{picture}%
}}

\def\putVtriangle(#1,#2)[#3`#4`#5;#6`#7`#8]{%
\xpos=#1 \ypos=#2
\advance\ypos by\c@height
{\multiply\c@height by2
\puthmorphism(\xpos,\ypos)[#3`#4`{#6}]{\c@height}{\c@arrowtypea}a}%
\putmorphism(\xpos,\ypos)(1,-1)[`#5`{#7}]{\c@height}{\c@arrowtypeb}l%
\advance\xpos by\c@height
\advance\xpos by\c@height
\putmorphism(\xpos,\ypos)(-1,-1)[``{#8}]{\c@height}{\c@arrowtypec}r%
}

\def\Vtriangle[#1`#2`#3;#4`#5`#6]{{
\c@width=\c@height                         
\xext=\c@width                           
\yext=\c@height                          
\topadjust[#1`#2`{#4}]
\botadjust[#3``]
\multiply \xext by2 
\leftsladjust[#1`#3`{#5}]
\rightsladjust[#2`#3`{#6}]
\begin{picture}(\xext,\yext)(\xoff,\yoff)%
\putVtriangle(0,0)[#1`#2`#3;#4`#5`{#6}]%
\end{picture}%
}}

\def\putVtrianglepair(#1,#2)[#3]{\xpos=#1 \ypos=#2%
\putVtrianglepairx[#3]}
\def\putVtrianglepairx[#1`#2`#3`#4;#5`#6`#7`#8`#9]{%
\advance\ypos by\c@height
\putmorphism(\xpos,\ypos)(1,-1)[`#4`{#7}]{\c@height}{\c@arrowtypec}l%
\puthmorphism(\xpos,\ypos)[#1`#2`{#5}]{\c@height}{\c@arrowtypea}a%
\advance\xpos by\c@height
\puthmorphism(\xpos,\ypos)[\phantom{#2}`#3`{#6}]{\c@height}{\c@arrowtypeb}a%
\putvmorphism(\xpos,\ypos)[``{#8}]{\c@height}{\c@arrowtyped}m%
\advance\xpos by\c@height
\putmorphism(\xpos,\ypos)(-1,-1)[``{#9}]{\c@height}{\c@arrowtypee}r%
}

\def\Vtrianglepair[#1`#2`#3`#4;#5`#6`#7`#8`#9]{{%
\xoff=0
\yoff=2 
\xext=\c@height                  
\c@width=\c@height                 
\yext=\c@height                  
\vertadjust[#1`#2`{#5}]
\tempcountd=\tempcounta        
\vertadjust[#2`#3`{#6}]
\ifnum\tempcounta<\tempcountd
\tempcounta=\tempcountd\fi
\advance \yext by\tempcounta
\botadjust[#4``]%
\multiply \xext by2
\leftsladjust[#1`#4`{#7}]%
\rightsladjust[#3`#4`{#9}]%
\begin{picture}(\xext,\yext)(\xoff,\yoff)%
\putVtrianglepair(0,0)[#1`#2`#3`#4;#5`#6`#7`#8`{#9}]%
\end{picture}%
}}

\def\putCtriangle(#1,#2)[#3`#4`#5;#6`#7`#8]{%
\xpos=#1 \ypos=#2
\advance\ypos by\c@height
\putmorphism(\xpos,\ypos)(1,-1)[``{#8}]{\c@height}{\c@arrowtypec}l%
\advance\xpos by\c@height
\advance\ypos by\c@height
\putmorphism(\xpos,\ypos)(-1,-1)[#3`#4`{#6}]{\c@height}{\c@arrowtypea}l%
{\multiply\c@height by 2
\putvmorphism(\xpos,\ypos)[`#5`{#7}]{\c@height}{\c@arrowtypeb}r}%
}

\def\Ctriangle[#1`#2`#3;#4`#5`#6]{{
\c@width=\c@height                          
\xext=\c@width                            
\yext=\c@height                           
\multiply \yext by2 
\topadjust[#1``]
\botadjust[#3``]
\sladjust[#2`#1`{#4}]{\c@width}
\tempcountd=\tempcounta                 
\sladjust[#2`#3`{#6}]{\c@width}
\ifnum \tempcounta<\tempcountd          
\tempcounta=\tempcountd\fi              
\advance \xext by\tempcounta            
\advance \xoff by-\tempcounta           
\rightadjust[#1`#3`{#5}]
\begin{picture}(\xext,\yext)(\xoff,\yoff)%
\putCtriangle(0,0)[#1`#2`#3;#4`#5`{#6}]%
\end{picture}%
}}

\def\putDtriangle(#1,#2)[#3`#4`#5;#6`#7`#8]{%
\xpos=#1 \ypos=#2
\advance\xpos by\c@height \advance\ypos by\c@height
\putmorphism(\xpos,\ypos)(-1,-1)[``{#8}]{\c@height}{\c@arrowtypec}r%
\advance\xpos by-\c@height \advance\ypos by\c@height
\putmorphism(\xpos,\ypos)(1,-1)[`#4`{#7}]{\c@height}{\c@arrowtypeb}r%
{\multiply\c@height by 2
\putvmorphism(\xpos,\ypos)[#3`#5`{#6}]{\c@height}{\c@arrowtypea}l}%
}

\def\Dtriangle[#1`#2`#3;#4`#5`#6]{{
\c@width=\c@height                         
\xext=\c@height                          
\yext=\c@height                          
\multiply \yext by2 
\topadjust[#1``]
\botadjust[#3``]
\leftadjust[#1`#3`{#4}]
\sladjust[#2`#1`{#4}]{\c@height}
\tempcountd=\tempcountd                
\sladjust[#2`#3`{#6}]{\c@height}
\ifnum \tempcounta<\tempcountd         
\tempcounta=\tempcountd\fi             
\advance \xext by\tempcounta           
\begin{picture}(\xext,\yext)(\xoff,\yoff)
\putDtriangle(0,0)[#1`#2`#3;#4`#5`{#6}]%
\end{picture}%
}}

\def\setrecparms[#1`#2]{\c@width=#1 \c@height=#2}%
\def\recurse[#1`#2`#3`#4;#5`#6`#7`#8`#9]{{%
\settowidth{\tempdimen}{#1}
\ifdim\tempdimen=0pt
  \savebox{\tempboxa}{\hbox{#2}}%
  \savebox{\tempboxb}{\hbox{#4}}%
  \savebox{\tempboxc}{\hbox{#7}}%
\else
  \savebox{\tempboxa}{\hbox{$\hbox{#1}\times\hbox{#2}$}}%
  \savebox{\tempboxb}{\hbox{$\hbox{#1}\times\hbox{#4}$}}%
  \savebox{\tempboxc}{\hbox{$\hbox{#1}\times\hbox{#7}$}}%
\fi
\tempcounte=\c@height
\divide\tempcounte by 2
\tempcountf=\tempcounte
\advance\tempcountf by \c@width
\xext=\tempcountf \yext=\c@height
\topadjust[#2`\usebox{\tempboxa}`{#5}]%
\botadjust[#4`\usebox{\tempboxb}`{#9}]%
\sladjust[#3`#2`{#6}]{\tempcounte}%
\tempcountd=\tempcounta
\sladjust[#3`#4`{#8}]{\tempcounte}%
\ifnum \tempcounta<\tempcountd
\tempcounta=\tempcountd\fi
\advance \xext by\tempcounta
\advance \xoff by-\tempcounta
\rightadjust[\usebox{\tempboxa}`\usebox{\tempboxb}`\usebox{\tempboxc}]%
\bfig
{\settriparms[-1`1`1;\tempcounte]%
\putCtriangle(0,0)[`#3`;#6`#7`{#8}]}%
\c@arrowtypea=-1 \c@arrowtypeb=0 \c@arrowtypec=1 \c@arrowtyped=-1
\putsquare(\tempcounte,0)[#2`\usebox{\tempboxa}`#4`\usebox{\tempboxb};%
#5``\usebox{\tempboxc}`#9]%
\efig
}}

\textefrancais

\indentationtitresection=\parindent
\indentationtitresubsection=\parindent
\indentationtitresubsubsection=\parindent
\indentationtitreNumero=\parindent
\indentationtitrenumero=\parindent
\indentationparagraphe=\parindent
\indentationtitreenonce=\parindent

\indentationtextechap
\def\indentationtextesection{+}
\def\indentationtextesubsection{+}
\def\indentationtextesubsubsection{+}
\def\indentationtexteNumero{+}

\def\avantnumsection{}
\def\apresnumsection{.}

\def\avantnumsubsection{}
\def\apresnumsubsection{.}

\def\avantnumsubsubsection{}
\def\apresnumsubsubsection{.}

\def\avantnumnumero{}
\def\apresnumnumero{.}
\def\fintitrenumero{{\rm\pointir}}

\def\avantnumNumero{}
\def\apresnumNumero{.}
\def\fintitreNumero{.}

\def\avantnumenonce{\ifthmgauche (\else{}\fi}
\def\apresnumenonce{\ifthmgauche)\else{}\fi}
\def\fintitreenonce{\pointir}

\def\avantnumEnonce{\ifthmgauche (\else{}\fi}
\def\apresnumEnonce{\ifthmgauche)\else{}\fi}
\def\fintitreEnonce{\pointt}

\let\styletitresection=\bf
\let\styletitresubsection=\bf
\let\styletitresubsubsection=\bf
\let\styletitrenumero=\sl
\let\styletitreNumero=\sl
\let\styletitreparagraphe=\bf
\let\styletitreenonce=\sc
\let\styleprecisionenonce=\rm
\let\styletexteenonce=\it

\def\espaceavanttitrechapsommaire{\hspace{1em}}
\let\newespaceavanttitrechapsommaire=\espaceavanttitrechapsommaire
\def\stylenumchap{\large\bf}
\def\avantchap{\centering}
\def\apreschap{}
\def\styletitrechap{\Large\bf}
\def\Chapitre{\iftextefrancais Chapitre\else Chapter\fi}

\def\titretableofcontents{\iftextefrancais Sommaire\else Contents\fi}
\def\titrelistoffigures{\iftextefrancais Liste des figures\else List of figures\fi}
\def\titrelistoftables{\iftextefrancais Liste des tableaux\else List of tables\fi}
\def\titrebibliographie{\iftextefrancais R\' ef\' erences\else References\fi}

\def\sectionheadingsplein{\T\thesection.\ }
\def\sectionheadingsvide{}
\def\subsectionheadingsplein{\T\thesubsection.\ }
\def\subsectionheadingsvide{}

\let\styleauteur=\sc
\let\styletitre=\it
\let\styletitrelivre=\sl
\let\stylejournal=\rm
\let\stylevolume=\bf
\let\styleannee=\rm
\let\stylepages=\rm

\let\styleediteur=\rm
\let\styleanneelivre=\rm

\let\styleflechehoriz=\textstyle

\def\lieulogo{/usr/local/logo/}

\ps@plain  \onecolumn \if@twoside\else\fi
\def\sectionmark#1{}

\numerman
\textefrancais
\citationman
\thmdroite
\eqngauche

\def\@normalsize{\@setsize\normalsize{15pt}\xiipt\@xiipt
\font\egtbf=cmbx8 \font\sixbf=cmbx6\scriptfont\bffam\egtbf \scriptscriptfont\bffam\sixbf 
\def\bm{\fam\bmfam\twlmib}\textfont\bmfam\twlmib
    \scriptfont\bmfam\egtmib \scriptscriptfont\bmfam\sixmib
\def\got{\fam\gotfam\twlgot}
\textfont\gotfam\twlgot \scriptfont\gotfam\egtgot \scriptscriptfont\gotfam\sixgot
\def\oldstyle{\fam\@ne\twlmi}
\abovedisplayskip 12pt plus3pt minus7pt\belowdisplayskip \abovedisplayskip
\abovedisplayshortskip \z@ plus3pt\belowdisplayshortskip 6.5pt plus3.5pt
minus3pt}

\def\small{\@setsize\small{13.6pt}\xipt\@xipt
\def\bm{\fam\bmfam\elvmib}\textfont\bmfam\elvmib
    \scriptfont\bmfam\egtmib \scriptscriptfont\bmfam\sixmib
\def\got{\fam\gotfam\elvgot}
\textfont\gotfam\elvgot \scriptfont\gotfam\egtgot \scriptscriptfont\gotfam\sixgot
\def\oldstyle{\fam\@ne\elvmi}
\abovedisplayskip 11pt plus3pt minus6pt\belowdisplayskip \abovedisplayskip
\abovedisplayshortskip \z@ plus3pt\belowdisplayshortskip 6.5pt plus3.5pt
minus3pt
\def\@listi{\parsep 4.5pt plus 2pt minus 1pt
 \itemsep \parsep
 \topsep 9pt plus 3pt minus 5pt}}

\def\footnotesize{\@setsize\footnotesize{12pt}\xpt\@xpt
\def\bm{\fam\bmfam\tenmib}\textfont\bmfam\tenmib
    \scriptfont\bmfam\sevmib \scriptscriptfont\bmfam\fivmib
\def\got{\fam\gotfam\tengot}
\textfont\gotfam\tengot \scriptfont\gotfam\sevgot \scriptscriptfont\gotfam\fivgot
\def\oldstyle{\fam\@ne\tenmi}
\abovedisplayskip 10pt plus2pt minus5pt\belowdisplayskip \abovedisplayskip
\abovedisplayshortskip \z@ plus3pt\belowdisplayshortskip 6pt plus3pt minus3pt
\def\@listi{\topsep 6pt plus 2pt minus 2pt\parsep 3pt plus 2pt minus 1pt
\itemsep \parsep}}

\def\scriptsize{\@setsize\scriptsize{9.5pt}\viiipt\@viiipt
\def\bm{\fam\bmfam\egtmib}\textfont\bmfam\egtmib
    \scriptfont\bmfam\sixmib \scriptscriptfont\bmfam\fivmib
\def\got{\fam\gotfam\egtgot}
\textfont\gotfam\egtgot \scriptfont\gotfam\sixgot \scriptscriptfont\gotfam\fivgot
}

\def\tiny{\@setsize\tiny{7pt}\vipt\@vipt
\def\bm{\fam\bmfam\sixmib}\textfont\bmfam\sixmib
    \scriptfont\bmfam\sixmib \scriptscriptfont\bmfam\sixmib
\def\got{\fam\gotfam\sixgot}
\textfont\gotfam\sixgot \scriptfont\gotfam\sixgot \scriptscriptfont\gotfam\sixgot
}

\def\large{\@setsize\large{18pt}\xivpt\@xivpt
\def\bm{\fam\bmfam\frtnmib}\textfont\bmfam\frtnmib
    \scriptfont\bmfam\tenmib \scriptscriptfont\bmfam\sevmib
\def\got{\fam\gotfam\frtngot}
\textfont\gotfam\frtngot \scriptfont\gotfam\tengot \scriptscriptfont\gotfam\sevgot
\def\oldstyle{\fam\@ne\frtnmi}}

\def\Large{\@setsize\Large{22pt}\xviipt\@xviipt
\def\bm{\fam\bmfam\svtnmib}\textfont\bmfam\svtnmib
    \scriptfont\bmfam\twlmib \scriptscriptfont\bmfam\tenmib
\def\got{\fam\gotfam\svtngot}
\textfont\gotfam\svtngot \scriptfont\gotfam\twlgot \scriptscriptfont\gotfam\tengot
\def\oldstyle{\fam\@ne\svtnmi}}

\def\LARGE{\@setsize\LARGE{25pt}\xxpt\@xxpt
\def\bm{\fam\bmfam\twtymib}\textfont\bmfam\twtymib
    \scriptfont\bmfam\frtnmib \scriptscriptfont\bmfam\twlmib
\def\got{\fam\gotfam\twtygot}
\textfont\gotfam\twtygot \scriptfont\gotfam\frtngot \scriptscriptfont\gotfam\twlgot
\def\oldstyle{\fam\@ne\twtymi}}

\def\huge{\@setsize\huge{30pt}\xxvpt\@xxvpt}

\normalsize

\if@twoside \oddsidemargin 21pt \evensidemargin 59pt \marginparwidth 85pt
\else \oddsidemargin 39.5pt \evensidemargin 39.5pt
\fi
\skip\footins 10.8pt plus 4pt minus 2pt
minus 4pt \dbltextfloatsep 20pt plus 2pt minus 4pt \@dblmaxsep 20pt
\@beginparpenalty -\@lowpenalty \@endparpenalty -\@lowpenalty \@itempenalty
-\@lowpenalty

\def\part{\par \addvspace{4ex} \indentationtextechappm \secdef\@part\@spart}
\def\@part[#1]#2{\ifnum \c@secnumdepth >\m@ne \refstepcounter{part}
\addcontentsline{toc}{part}{\thepart \newespaceavanttitrechapsommaire #1}\else
\addcontentsline{toc}{part}{#1}\fi { \avantchap
 \ifnum \c@secnumdepth >\m@ne \stylenumchap  \Chapitre {} \thepart
\par\nobreak \fi
\styletitrechap  #2 \markboth{}{}\apreschap\par } \nobreak
\vskip 3ex \@afterheading }
\def\@spart#1{{ \avantchap
\styletitrechap  #1\markboth{}{}\apreschap\par} \nobreak \vskip 3ex \@afterheading }

\def\sectioncentre{\par \addvspace{3.5ex plus .1ex minus .2ex}\indentationtextechappm \secdef\@sectioncentre\@ssectioncentre}
\def\@sectioncentre[#1]#2{\ifnum \c@secnumdepth >\z@ \refstepcounter{section}
\addcontentsline{toc}{section}{\thesection \hskip 1 em #1}\else
\addcontentsline{toc}{section}{#1}\fi
\sectionmark{#1}
\let\sectionheadings=\sectionheadingsplein{\centering \ifnum \c@secnumdepth >\z@
\styletitresection \avantnumsection\thesection\apresnumsection \hskip 1em #2\par\nobreak\fi}
\nobreak\vskip 2.3ex\@afterheading }
\def\@ssectioncentre#1{\ifnum \c@secnumdepth >\z@ \stepcounter{section}\addtocounter{section}{-1}\fi\let\sectionheadings=\sectionheadingsvide{\centering \styletitresection  #1\par\nobreak}
\nobreak\vskip 2.3ex\@afterheading}

\def\section{\@ifstar{\let\sectionheadings=\sectionheadingsvide}{\let\sectionheadings=\sectionheadingsplein}\let\avantnumsections=\avantnumsection\let\apresnumsections=\apresnumsection\@startsection{section}{1}{\indentationtitresection}{\indentationtextesection3.5ex plus\indentationtextesection1ex minus \indentationtextesection.2ex}{2.3ex plus .2ex}{\styletitresection}}

\def\subsection{\@ifstar{\let\subsectionheadings=\subsectionheadingsvide}{\let\subsectionheadings=\subsectionheadingsplein}\let\avantnumsections=\avantnumsubsection\let\apresnumsections=\apresnumsubsection\@startsection{subsection}{2}{\indentationtitresubsection}{\indentationtextesubsection3.5ex plus
\indentationtextesubsection1ex minus \indentationtextesubsection.2ex}{1.5ex plus .2ex}{\styletitresubsection}}

\def\subsubsection{\let\avantnumsections=\avantnumsubsubsection\let\apresnumsections=\apresnumsubsubsection\@startsection{subsubsection}{3}{\indentationtitresubsubsection}{\indentationtextesubsubsection3.5ex plus
\indentationtextesubsubsection1ex minus \indentationtextesubsubsection.2ex}{1.5ex plus .2ex}{\styletitresubsubsection}}

\def\slsubsubsubsection{\let\avantnumsections=\avantnumNumero\let\apresnumsections=\apresnumNumero\@startsection{subparagraph}{4}{\indentationtitreNumero}{\indentationtexteNumero3.5ex plus
\indentationtexteNumero1ex minus \indentationtexteNumero.2ex}{3pt plus 1pt minus 1pt}{\styletitreNumero}}

\def\newparag{\@startsection{paragraph}{5}{\indentationparag}{3.25ex
plus 1ex minus .2ex}{-1em}{\normalsize\bf}}

\def\newparagnoind{\@startsection{paragraph}{5}{\z@}{-3.25ex
plus -1ex minus -.2ex}{\z@}{\normalsize\bf}}

\def\paragraph{\@startsection{paragraph}{5}{\indentationparagraphe}{3.25ex plus
1ex minus .2ex}{-1em}{\styletitreparagraphe}}

\def\subparagraph{\@startsection{subparagraph}{4}{1.2em}{3.25ex
plus 1ex minus .2ex}{-1em}{\normalsize\bf}}

\def\slsubparagraph{\let\avantnumsections=\avantnumnumero\let\apresnumsections=\apresnumnumero\@startsection{subparagraph}{4}{\indentationtitrenumero}{3.25ex
plus 1ex minus .2ex}{-1em}{\styletitrenumero}}

\def\scparagraph{\@startsection{paragraph}{5}{0pt}{3.25ex
plus 1ex minus .2ex}{-1em}{\sc}}

\leftmargin\leftmargini
\labelwidth\leftmargini\advance\labelwidth-\labelsep
\def\@listi{\leftmargin\leftmargini}
\def\@listii{\leftmargin\leftmarginii
 \labelwidth\leftmarginii\advance\labelwidth-\labelsep
 \topsep 5pt plus 2.5pt minus 1pt
 \parsep 2.5pt plus 1pt minus 1pt
 \itemsep \parsep}
\def\@listiii{\leftmargin\leftmarginiii
 \labelwidth\leftmarginiii\advance\labelwidth-\labelsep
 \topsep 2.5pt plus 1pt minus 1pt
 \parsep \z@ \partopsep 1pt plus 0pt minus 1pt
 \itemsep \topsep}
\def\@listiv{\leftmargin\leftmarginiv
 \labelwidth\leftmarginiv\advance\labelwidth-\labelsep}
\def\@listv{\leftmargin\leftmarginv
 \labelwidth\leftmarginv\advance\labelwidth-\labelsep}
\def\@listvi{\leftmargin\leftmarginvi
 \labelwidth\leftmarginvi\advance\labelwidth-\labelsep}

\makeatother

\texteanglais
\numerman
\begin{document}

\indentationtitresection=0pt
\indentationtitresubsection=0pt
\indentationtitresubsubsection=0pt
\indentationtitreNumero=0pt
\indentationtitrenumero=0pt
\indentationparag=.5em
\indentationparagraphe=0pt
\indentationtitreenonce=0pt

\def\limind{\mathop{\oalign{lim\cr
\hidewidth$\longrightarrow$\hidewidth\cr}}}

\def\limproj{\mathop{\oalign{lim\cr
\hidewidth$\longleftarrow$\hidewidth\cr}}}

\def\boxit#1#2{\setbox1=\hbox{\kern#1{#2}\kern#1}%
\dimen1=\ht1 \advance\dimen1 by #1
\dimen2=\dp1 \advance\dimen2 by #1
\setbox1=\hbox{\vrule height\dimen1 depth\dimen2\box1\vrule}%
\setbox1=\vbox{\hrule\box1\hrule}%
\advance\dimen1 by .4pt \ht1=\dimen1
\advance\dimen2 by .4pt \dp1=\dimen2 \box1\relax}
\def \pext{\, \hbox{\boxit{0pt}{$\times$}}\,}

\def\barql{\bar \QQ_{\ell}}
\def\fd{F \!\hbox{-} \Delta}
\def\fdphi{F_{\varphi} \hbox{-} \Delta}
\def\ker{{\rm Ker} \,}
\def\st{{\rm St}}
\def\repr{{\rm Repr}}
\def\ind{{\rm Ind}}
\def\ch{{\rm ch} \,}
\def\chcl{{\rm chcl} \,}
\def\coker{{\rm Coker} \,}
\def\ord{{\rm ord}}
\def\jac{{\cal J\!{\em ac}}}
\title{Germs of arcs on singular algebraic varieties\\ and motivic integration}
\author{Jan Denef \and Fran\c cois Loeser}
\date{Revised Nov. 1997, to appear in Invent. Math.}
\maketitle

\bigskip



\sec{1}{Introduction}Let $k$ be a field of characteristic zero.
We denote by $\cM$ the Grothendieck ring
of algebraic varieties over $k$ ({\it i.e.}
reduced separated schemes of finite type over $k$). 
It is the ring generated by symbols $[S]$, for $S$ an algebraic variety over
$k$, 
with the relations
$[S] = [S']$ if $S$ is isomorphic to $S'$,
$[S] = [S \setminus S'] + [S']$ if $S'$ is closed in $S$
and
$[S \times S'] = [S] \, [S']$.
Note that, for $S$ an algebraic variety
over
$k$,
the mapping $S' \mapsto [S']$
from the set of closed subvarieties of  $S$ extends uniquely to a mapping
$W \mapsto [W]$
from the set of constructible subsets of $S$ to $\cM$,
satisfying
$[W \cup W'] = [W] + [W'] - [W \cap W']$.
We set $\LL := [\AA^1_k]$ and
$\cM_{\rm loc} := \cM [\LL^{-1}]$.
We denote by $\cM [T]_{\rm loc}$ 
the subring of $\cM_{\rm loc} [[T]]$
generated by 
$\cM_{\rm loc} [T]$  and the series
$(1 -\LL^a T^b)^{-1}$ with
$a$ in $\ZZ$ 
and $b$ in $\NN \setminus \{0\}$.

Let $X$ be an algebraic variety over
$k$.
We denote by $\cL (X)$ the scheme of germs of arcs on $X$. It is a
scheme over $k$ and for any field extension $k \subset K $
there is a natural bijection
$$
\cL (X) (K) \simeq
{\rm Mor}_{k-{\rm schemes}} ({\rm Spec} \, K[[t]], X),
$$
between the set of $K$-rational points of $\cL (X)$ and the set of $K
[[t]]$-rational points of $X$ (called the set of germs of arcs with
coefficients in $K$ on $X$). More precisely, the scheme $\cL (X)$
is defined as
the projective limit
$
\cL (X) := \limproj \cL_{n} (X),
$
in the category of $k$-schemes, of the 
schemes $\cL_{n} (X)$, $n \in \NN$, representing
the functor
$$
R \mapsto {\rm Mor}_{k-{\rm schemes}} ({\rm Spec}
\, R [t]/t^{n+1} R [t], X),
$$
defined on the category of $k$-algebras. (Thus, for any
$k$-algebra $R$, the set of $R$-rational points of $\cL_{n} (X)$
is naturally identified with the set of 
$R [t]/t^{n+1} R [t]$-rational points of $X$.)
The existence of $\cL_{n} (X)$ is well known, cf. [B-L-R] p.276, and the
projective limit exists
since the transition morphisms are affine. We shall denote 
by $\pi_n$ the canonical morphism $\cL (X) \rightarrow \cL_n (X)$
corresponding to truncation of arcs.
In the present paper, the schemes $\cL (X)$ and $\cL_{n} (X)$ will
always be considered with their reduced structure.
Note that the set-theoretical image
$\pi_n (\cL (X))$ is a constructible subset
of $\cL_{n} (X)$, as follows from a theorem of M. Greenberg [G],
see (4.4) below.
These constructible
sets $\pi_n (\cL (X))$ were first studied by J. Nash in [N], in relation
with Hironaka's resolution of singularities. They are also considered in the
papers [L-J], [H].

\medskip
The following result is the first
main result of the paper. It is an analogue
of the rationality of the Poincar\'e series associated to the
$p$-adic points on a variety
proved in [D1].

\begin{th}{1.1}Let $X$ be an algebraic variety over
$k$. The power series
$$P (T) := \sum_{n = 0}^{\infty} \, [\pi_n (\cL (X))] \, T^n,$$
considered as an element of 
$\cM_{\rm loc} [[T]]$, is rational and belongs to
$\cM [T]_{{\rm loc}}$.
\end{th}

The proof of the theorem is given in section 5 and
uses two main ingredients.
The first one is a result of 
J. Pas [P] on quantifier
elimination for semi-algebraic sets of power series
in characteristic zero.
The second one is 
M. Kontsevich's marvellous idea of 
motivic integration [K]. More precisely,
M. Kontsevich introduced the completion
$\widehat \cM$ of 
$\cM_{\rm loc}$ with respect to the filtration
$F^m \cM_{\rm loc}$, where $F^m \cM_{\rm loc}$ is
the subgroup of $\cM_{\rm loc}$ generated by
$\{ [S] \, \LL^{- i} \bigm \vert
i - \dim S \geq m\}$,
and defined, for smooth $X$, a motivic integration on $\cL (X)$ with
values
into $\widehat \cM$. This is an analogue of classical $p$-adic integration.
In the present paper we extend Kontsevich's
construction to semi-algebraic subsets of
$\cL (X)$,
with $X$ any pure dimensional algebraic variety over
$k$, not necessarily smooth.
For such an $X$, let $\BB$ be the set of all
semi-algebraic
subsets
of $\cL (X)$. We construct in section 3 a canonical measure
$\mu : \BB \rightarrow
\widehat \cM$. This allows us to define integrals
$$
\int_{A} \LL^{- \alpha} d \mu,
$$
for $A$ in $\BB$ and
$\alpha : A \rightarrow \ZZ \cup \{+ \infty\}$
a simple function
which is bounded from below. (Semi-algebraic
subsets
of $\cL (X)$ and simple functions are defined in section 2.)
The properties of this motivic integration, together
with resolution of singularities and the result of Pas, suffice
to prove the rationality of the image of $P (T)$ in $\widehat \cM
[[T]]$.
To prove the rationality of $P (T)$, considered as an element of
$\cM_{\rm loc } [[T]]$, one needs a more refined argument based on Lemma
2.8 and the use of an obvious lifting $\tilde \mu (A)$ in $\cM_{\rm loc}$
of $\mu (A)$, when $A$ is a stable semi-algebraic subset of $\cL (X)$
(a notion defined in section 2).

For an algebraic variety $X$, it is natural to consider
its motivic volume $\mu (\cL (X))$. In section 6, we give explicit
formulas for
$\mu (\cL (X))$ in terms of certain special
resolutions of singularities
of $X$ (which always exist). As a corollary we deduce that 
$\mu (\cL (X))$ always belongs to a certain localization of
the image of
$\cM_{\rm loc}$ in $\widehat \cM$
on which the Euler characteristic $\chi$
naturally extends with rational values. 
So we obtain a new invariant
of $X$,
the Euler characteristic
$\chi (\mu (\cL (X)))$, which is a rational number and coincides with 
the usual Euler characteristic of $X$ when
$X$ is smooth.
In section 7, we prove that, when $X$ is of pure
dimension
$d$,
the sequence $[\pi_n (\cL (X))] \, \LL^{- (n +1 )d}$ converges
to $\mu (\cL (X))$ in 
$\widehat \cM$. 
This result, which is an analogue of a $p$-adic
result
by J. Oesterl\'e [O],
gives, in some sense, a precise meaning to Nash's
guess one should consider the limit of the $\pi_n (\cL (X))$'s.
We conclude the paper by some
remarks on the Greenberg function in section 8.

For related results concerning motivic
Igusa functions, see [D-L].

\sec{2}{Semi-algebraic sets of power series}
\noindent (2.1) From now on we will denote by
$\bar k$ a fixed algebraic closure of $k$, and by $\bar k((t))$
the fraction field of $\bar k [[t]]$, where $t$ is one variable.
Let $x_1, \ldots, x_m$
be variables running over $\bar k ((t))$ and let $\ell_1, \ldots,
\ell_r$ be variables running over $\ZZ$.
A {\it semi-algebraic} condition
$\theta (x_1, \ldots, x_m; \ell_1, \ldots, \ell_r)$
is a finite boolean combination of conditions of the form
\begin{eqn}
\ord_t f_1 (x_1, \ldots, x_m) &\geq&
\ord_t f_2 (x_1, \ldots, x_m) + L (\ell_1, \ldots, \ell_r)
\numeqn{(i)}\\
\ord_t f_1 (x_1, \ldots, x_m) &\equiv&
L (\ell_1, \ldots, \ell_r) \; \; \; {\rm mod} \, d,
\numeqn{(ii)}
\end{eqn}
and
\begin{eqn}
h (\overline{ac} (f_1 (x_1, \ldots, x_m)),
\ldots,
\overline{ac} (f_{m'} (x_1, \ldots, x_m))) &=& 0,
\numeqn{(iii)}
\end{eqn}
where $f_i$ and $h$ are polynomials over $k$,
$L$ is a polynomial of degree $\leq 1$ over $\ZZ$, $d \in \NN$,
and $\overline{ac} (x)$ is the coefficient of lowest degree
of $x$ in $\bar k ((t))$ if $x \not= 0$, and is equal to 0 otherwise.
Here we use the convention that $(+ \infty) + \ell = + \infty$
and
$+ \infty \equiv \ell \; {\rm mod} \, d$, for all $\ell \in \ZZ$.
In particular the condition 
$f (x_1, \ldots, x_m) = 0$ is a semi-algebraic condition, for $f$ a
polynomial over $k$.
A subset of $\bar k ((t))^{m} \times \ZZ^{r}$
defined by a semi-algebraic condition is called semi-algebraic.
One defines similarly semi-algebraic subsets of 
$K ((t))^{m} \times \ZZ^{r}$ for $K$ an algebraically closed field
containing
$\bar k$.

A function $\alpha : \bar k ((t))^m \times \ZZ^n \rightarrow \ZZ$
is called {\it simple} if its graph is semi-algebraic. An easy result of
Presburger [Pr] implies that
$(\exists \, \ell_1 \in \ZZ) \, \theta (x_1, \ldots, x_m; \ell_1,
\ldots, \ell_r)$
is semi-algebraic when $\theta$ is a semi-algebraic condition.

We will use the following result on quantifier
elimination due to J. Pas [P].

\begin{th}{2.1}[{J. Pas [P]}]If $\theta$ is a semi-algebraic condition,
then $$(\exists \, x_1 \in \bar k ((t))) \,\, \theta (x_1, \ldots, x_m; \ell_1,
\ldots, \ell_r)$$
is semi-algebraic. Furthermore, for any
algebraically closed field $K$ 
containing
$\bar k$, 
$$(\exists \, x_1 \in K ((t))) \,\, \theta (x_1, \ldots, x_m; \ell_1,
\ldots, \ell_r)$$
is also semi-algebraic and may be defined by the same conditions (\ie.
independently of $K$).
\end{th}

Indeed, the first assertion follows from Theorem 4.1 in [P] together
with Chevalley's constructibility theorem and the above mentioned result
of Presburger; the second assertion follows directly from the remark at
the begining of \S\kern .15em 3 of [P].

The theorem of Pas is a refinement of older quantifier elimination
results of Ax and Kochen [A-K], and of Delon [De].

\bigskip
\noindent (2.2) Let $X$ be an
algebraic variety over $k$.
For $x \in \cL (X)$, we denote by $k_{x}$
the residue field of $x$ on $\cL (X)$,
and by $\tilde x$ the corresponding rational point
$\tilde x \in \cL (X) (k_{x}) = X (k_{x}[[t]])$.
When there is no danger of confusion we will often write $x$
instead of $\tilde x$.
A {\it semi-algebraic family of semi-algebraic subsets}
(for $n = 0$ a semi-algebraic subset)
$A_{\ell}$, $\ell \in \NN^n$, of $\cL (X)$ is a family of subsets
$A_{\ell}$ 
of $\cL (X)$ such that there exists a covering of $X$ by affine 
Zariski open
sets $U$ with
$$
A_{\ell} \cap \cL (U) =
\{ x \in \cL (U) \bigm \vert
\theta (h_1 (\tilde x), \ldots, h_{m} (\tilde x); \ell)\},
$$
where $h_1, \ldots, h_{m}$ are regular functions on 
$U$
and $\theta$ is a semi-algebraic condition. 
Here the $h_{i}$'s and
$\theta$ may depend on $U$ and $h_{i} (\tilde x)$ belongs to $k_{x} [[t]]$.

Let $A$ be a semi-algebraic subset
of $\cL (X)$.
A function $\alpha : A \times \ZZ^n
\rightarrow \ZZ \cup \{ + \infty \}$ is called {\it simple} if the 
family of subsets $\{x \in \cL (X) \bigm | \alpha (x, \ell_1, \ldots,
\ell_n) = \ell_{n + 1}\}$,
$(\ell_1, \ldots,
\ell_{n + 1}) \in \NN^{n + 1}$,
is a semi-algebraic family of semi-algebraic subsets
of $\cL (X)$.

A {\it Presburger subset} of $\ZZ^{r}$ is a subset defined by a
semi-algebraic condition 
$\theta (\ell_{1}, \dots, \ell_{r})$ as in (2.1) with $m= 0$. A
function
$\alpha : \ZZ^{r} \rightarrow \ZZ$ is called
a {\it Presburger function} if its graph is a Presburger subset of
$\ZZ^{r + 1}$.

If $f : X \rightarrow Y$ is a morphism of algebraic varieties over $k$
and $A$ is a semi-algebraic subset of $\cL (X)$,
then $f (A)$ is a semi-algebraic subset of $\cL (Y)$, by Pas's Theorem.

We denote by $\pi_{n}$ the canonical morphism
$\cL (X) \rightarrow \cL_{n} (X)$. If necessary we will use also the
notation $\pi_{n, X}$. If $X$ is smooth, then $\pi_{n}$ is
surjective by
Hensel's Lemma.

\bigskip
\noindent (2.3) The following basic result is a consequence of Pas's
Theorem.

\begin{prop}{2.3}Let $X$ be an algebraic variety over $k$, and let $A$
be
a semi-algebraic subset of $\cL (X)$. Then $\pi_{n} (A)$ is a
constructible subset of $\cL_{n} (X)$.
\end{prop}

\demo We may assume $X = \AA^{n}_{k}$. Let $s : \cL_{n} (X)
\rightarrow
\cL (X)$ be a section of  $\pi_{n} : \cL (X) \rightarrow \cL_{n} (X)$
which maps $\cL_{n} (X) (\bar k)$ into 
$X (\bar k [t]) \subset X (\bar k [[t]]) = \cL (X) (\bar k)$. We have $y
\in
\pi_{n} (A)$ if and only if $s (y) \in \pi_{n}^{-1} \pi_{n} (A)$.
It is easy to verify that this implies the proposition because
$\pi_{n}^{-1} \pi_{n} (A)$ is a semi-algebraic subset of $\cL (X)$
by Theorem 2.1.  \hfill $\qed$

\bigskip
\noindent (2.4) Let $A$ be a semi-algebraic subset of $\cL (X)$. We call
$A$ {\it weakly stable at level} $n \in \NN$ if $A$ is a union of fibers
of $\pi_{n} : \cL (X) \rightarrow \cL_{n} (X)$. We call $A$ {\it weakly
stable}
if it stable at some level $n$. Note that weakly stable
semi-algebraic subsets form a boolean algebra.

\begin{lem}{2.4}For each $i \in \NN$, let $A_{i}$ be a
weakly stable
semi-algebraic subset of $\cL (X)$. Suppose that
$A := \bigcup_{i \in \NN} A_{i}$ is semi-algebraic and weakly stable. Then
$A$ equals
the union of a {\it finite} number of the $A_{i}'s.$
\end{lem}

\demo We may assume $X$ is affine. By looking at the complements of the
$A_{i}$'s, it is enough to prove the following assertion.
For each $i$ in $\NN$, let $B_{i}$ be a weakly stable semi-algebraic
subset of $\cL (X)$, and suppose that the intersection of finitely many
of the
$B_{i}$'s  is always nonempty. Then 
$\bigcap_{i \in \NN} B_{i}$ is nonempty. To prove this assertion note that
each $B_{i}$ is a {\it finite} boolean combination
of closed subschemes of $\cL (X)$ whose ideals are finitely generated,
since
$B_{i}$ is weakly stable at some level $n_{i}$ and $\pi_{n_{i}} (B)$ is
constructible. For any finite subset $\Sigma$ of $\NN$ there exists a
field $K_{\Sigma}$ containing $k$ and a 
$K_{\Sigma}$-rational point on $\bigcap_{i \in \Sigma} B_{i}$.
Considering the ultraproduct $K^{\ast}$ of the 
$K_{\Sigma}$'s with respect to a {\it suitable} ultrafilter, we obtain a
$K^{\ast}$-rational point in $\bigcap_{i \in \NN} B_{i}$.
(This kind of argument is very classical in model theory,
see, {\it e.g.}, page 172 of [C-K].) This proves the assertion. \hfill $\qed$

\rem The above lemma may also be deduced from Corollaire 7.2.7 of [G-D].

\bigskip
\noindent (2.5) Let $X$, $Y$ and $F$ be algebraic varieties over $k$,
and let
$A$, {\it resp.} $B$, be a constructible subset of $X$,
{\it resp.} $Y$. We say that
a map
$\pi : A \rightarrow B$ is a
{\it piecewise trivial fibration with fiber}
$F$, if there exists a finite partition of $B$ in subsets $S$ which are
locally closed
in $Y$ such that $\pi^{- 1} (S)$ is locally closed in $X$ and
isomorphic, as
a variety over $k$, to $S \times F$, with $\pi$
corresponding under the isomorphism to the projection
$S \times F \rightarrow S$. We say that the map $\pi$ is
a
{\it piecewise trivial fibration over} some constructible subset $C$ of
$B$,
if the restriction of $\pi$ to $\pi^{- 1} (C)$ is a piecewise 
trivial fibration onto $C$.

\bigskip
\noindent (2.6) For $X$ an algebraic variety over $k$
and $e$ in $\NN$, we will use the notation
$$
\cL^{(e)} (X) := \cL (X) \setminus \pi^{-1}_{e, X} (\cL_e (X_{\rm sing})),
$$
where $X_{\rm sing}$ denotes the singular locus of $X$.

\bigskip
\noindent (2.7) Let $X$ be an algebraic variety over $k$ of pure
dimension
$d$ (in particular we assume that
$X$ is non empty)
and let $A$ be a semi-algebraic subset of $\cL (X)$. We call $A$
{\it stable at level} $n \in \NN$, if $A$ is weakly
stable at level $n$ and $\pi_{m + 1} (\cL (X)) \rightarrow \pi_{m} (\cL
(X))$
is a piecewise trivial fibration over $\pi_{m} (A)$ with fiber
$\AA^{d}_{k}$ for all $m \geq n$.

We call $A$ {\it stable} if it stable at some level $n$. Note that the
family of stable semi-algebraic subsets of $\cL (X)$ is closed under
taking finite intersections and finite unions. If $A$ is stable at level
$n$, then $[\pi_{m} (A)] = [\pi_{n} (A)] \LL^{(m - n) d}$ for all
$m \geq n$. If $A$ is weakly stable and
$A \cap \cL (X_{\rm sing}) = \emptyset$ (which is
for instance the case when $X$ is smooth), then $A$ is stable. Indeed this
follows from Lemma 2.4 and Lemma 4.1 below, because then
$A$ is the union of the weakly stable subsets $A \cap \cL^{(e)} (X)$.

\bigskip
\noindent (2.8) We say that a semi-algebraic
family $A_{\ell}$, $\ell \in \NN^{n}$, of semi-algebraic subsets of $\cL
(X)$,
{\it has a bounded representation} if there exists a covering of $X$ by
affine Zariski open sets $U$ such that on each $U$ the family is given
by
a semi-algebraic condition $\theta$ (cf. (2.2)) with, in the notation of
(2.1),
$\ord_{t}f_{i}$ bounded on $A_{\ell} \cap U$ for each fixed $\ell$.

Clearly, if the family $A_{\ell}$ has a bounded representation then 
each 
$A_{\ell}$ is weakly stable. The next lemma is essential for the proof 
of Theorem 5.1$'$ on which Theorem 1.1 is based.

\begin{lem}{2.8}Let $X$ be a quasi-projective algebraic variety over 
$k$ and let $A_{\ell}$, $\ell \in \NN^{n}$, be a semi-algebraic
family of semi-algebraic subsets of $\cL (X)$. Assume that 
$A_{\ell}$ is weakly stable for each $\ell$. Then the family 
$A_{\ell}$ is a finite boolean combination of semi-algebraic
families of semi-algebraic subsets of $\cL (X)$ which
have bounded representations.
\end{lem}

\demo Because $X$ is quasi-projective we can work with a covering of 
$X$ by affine open sets $U_{i}$, $i \in I$, such that $U_{i} \setminus 
U_{j}$ is the locus in $U_{i}$ of a single regular function on $U_{i}$.
Hence we may assume that $X$ is affine and that the family $A_{\ell}$,
$\ell \in \NN^n$, is given by a semi-algebraic condition as in
(2.2) with $U = X$.

Let $f_{1}, \ldots, f_{r}$ be the regular functions on $X$ appearing 
in the conditions of the form 2.1 (i), (ii) and (iii) in the 
description of the semi-algebraic family $A_{\ell}$.

Assume that
for $i = 1, 2, \ldots, e \leq r$, $\ord_{t} f_{i}$ is not 
bounded on $A_{\ell}$ for some $\ell$, possibly depending on $i$,
and that $\ord_{t} f_{e + 1}$, \ldots, $\ord_{t} f_{r}$ are bounded
on $A_{\ell}$ by a function
$\nu (\ell)$ of $\ell$. Our proof is by induction on $e$. We may 
assume
that $A_{\ell}$ is weakly stable at level $\nu (\ell)$. Taking for
$\nu (\ell)$ the smallest integer satisfying the above requirements, 
we see by Theorem 2.1 that $\nu : \NN^n \rightarrow \NN$ is
a Presburger function.

By Greenberg's theorem [G], cf. (4.4) below, there exists a Presburger
function $\alpha : \NN^n \rightarrow \NN$, with
$\alpha (\ell) \geq \nu (\ell)$ for all $\ell \in \NN^n$, such that, 
for all $x$ in $\cL (X)$, if 
\begin{eqn}
f_{1}  (x) \equiv f_{2} (x) &\equiv& \cdots \equiv f_{e} (x) \equiv 0
\; \; \; {\rm mod} \, t^{\alpha (\ell) + 1},
\numeqn{(1)}
\end{eqn}
then there exists $x'$ in $\cL (X)$ with $x \equiv x'
\; \; \; {\rm mod} \, t^{\nu (\ell) + 1}$
and
\begin{eqn}
f_{1}  (x') = f_{2} (x') = \cdots = f_{e} (x')& =& 0.\numeqn{(2)}
\end{eqn}
Note that $A_{\ell}$ is the union of the following semi-algebraic 
subsets of $\cL (X)$ which are weakly stable at level $\alpha 
(\ell)$~:
\begin{eqn}
A_{\ell, 1} &:= &A_{\ell} \cap \{x \in \cL (X) \bigm |
\ord_{t} f_{1} (x) \leq \alpha (\ell)\},\\
&\vdots&\\
A_{\ell, e} &:= &A_{\ell} \cap \{x \in \cL (X) \bigm |
\ord_{t} f_{e} (x) \leq \alpha (\ell)\},\\
B_{\ell} &:= &A_{\ell} \cap \{x \in \cL (X) \bigm |
x \quad {\rm satisfies} \quad (1)\}.
\end{eqn}
The semi-algebraic families $A_{\ell, 1}, \ldots A_{\ell, e}$, 
$ \ell \in \NN^n$, are finite boolean combinations of semi-algebraic 
families which have bounded representations, because of the induction 
hypothesis. Thus it only remains to prove the lemma for the family
$B_{\ell}$ instead of $A_{\ell}$. In the description of
$A_{\ell}$ as a semi-algebraic family of semi-algebraic subsets of 
$\cL (X)$, mentioned in the beginning of the proof, replace $f_{1}$, 
$f_{2}$,
\dots,
$f_{e}$ by 0 and add the conditions $\ord_{t} f_{j} (x) \leq \nu 
(\ell)$ for $j = e + 1, \ldots , r$. In this way we obtain a new 
semi-algebraic family of semi-algebraic subsets of $\cL (X)$ which we 
denote by $A'_{\ell}$, $\ell \in \NN^n$. Clearly this family $A'_{\ell}$
has a bounded representation and for each fixed $\ell$ the set
$A'_{\ell}$ is weakly stable at level $\nu (\ell)$.

If $x$ in $\cL (X)$ satisfies (1), then there exists $x'$ in $\cL (X)$,
with $x \equiv x'
\; \; \; {\rm mod} \, t^{\nu (\ell) + 1}$, satisfying (2), and we have 
that $x \in A_{\ell}$ if and only if $x' \in A_{\ell}$, since $A_{\ell}$
is weakly stable at level $\nu (\ell)$, which in turn is equivalent 
by (2)
to $x' \in A'_{\ell}$ which,  since $A'_{\ell}$
is weakly stable at level $\nu (\ell)$, is verified if and only 
if
$x \in A'_{\ell}$.
Thus
\begin{eqn}
B_{\ell} &=& A'_{\ell} \cap \{x \in \cL (X) \bigm \vert
x \quad {\rm satisfies} \quad (1)\}\\
&=& A'_{\ell} \setminus \bigcup_{i = 1, \ldots, e}
\{x \in \cL (X) \bigm \vert \ord_{t} f_{i} (x) \leq \alpha (\ell) \}.
\end{eqn}
This proves the lemma, since $A'_{\ell}$ and
$\{x \in \cL (X) \bigm \vert \ord_{t} f_{i} (x) \leq \alpha (\ell) \}$
have a bounded representation for each $i$. \hfill $\qed$

\sec{3}{Motivic integration}The basic idea
of motivic integration goes back to M. Kontsevich [K]. We generalize
here his idea to the much more general setting of semi-algebraic
subsets of $\cL (X)$ and simple functions on $\cL (X)$,
for $X$ an
algebraic variety over $k$ which is not necessarily smooth.

\bigskip
\noindent (3.1) Let $X$ be an algebraic variety over $k$ of pure
dimension $d$. Denote by $\BB$ the set of all
semi-algebraic subsets of $\cL (X)$, and by
$\BB_{0}$ the set of all $A$ in $\BB$ which are stable. Clearly 
there is a unique additive measure
$$\tilde \mu : \BB_{0} \longrightarrow \cM_{\rm loc}$$
satisfying 
$$
\tilde \mu (A) = [\pi_{n} (A)] \, \LL^{- (n+1)d}
$$
when $A$ is stable at level $n$.
Sometimes we shall denote $\tilde \mu$ by
$\tilde \mu_{\cL (X)}$. Let $A$ be in $\BB_{0}$ and let 
$\alpha : A \rightarrow \ZZ$ be a simple function  all whose fibers 
are stable. Then, by Lemma 2.4, $\vert \alpha \vert$ is bounded, and 
we may define
$$
\int_{A} \LL^{-\alpha} d \tilde \mu :=
\sum_{n \in \ZZ} \LL^{-n} \, \tilde \mu (\alpha^{-1} (n)),
$$
the sum at the right hand side being finite.

Next we want to extend the measure $\tilde \mu$ on $\BB_{0}$ to a 
measure $\mu$ on $\BB$. The key to achieving that is the following lemma which 
allows to partition any $A$ in $\BB$ into stable subsets $A_{i}$, $i 
\in \NN$, and a set of ``measure zero''. However, this leads to 
infinite sums, and for this reason $\mu$ will take values in the
completion $\widehat \cM$ of $\cM_{\rm loc}$, see (3.2) below.

\begin{lem}{3.1}Let $X$ be an algebraic variety over $k$
of pure
dimension $d$, and let $A$ be a semi-algebraic subset of $\cL (X)$. There
exists a closed subvariety $S$ of $X$, with
${\rm dim} \, S < {\rm dim} \, X$,
and a semi-algebraic family $A_i$, $i \in \NN$, of semi-algebraic
subsets of $A$ such that
$\cL (S) \cap A$ and the $A_i$'s form a partition of
$A$, each $A_i$ is stable at some level $n_i$,
and
\begin{eqn}
\lim_{i \rightarrow \infty} \, ({\rm dim}\, \pi_{n_i} (A_i)
- (n_i + 1) \, d) = -\infty.\numeqn{(3.1.1)}
\end{eqn}
Moreover, if $\alpha : \cL (X) \rightarrow \ZZ$ is a simple function,
we can choose the partition such that
$\alpha$ is constant on
each $A_i$.
\end{lem}

\demo We may assume that $X$ is affine and irreducible and that $A$ is given
by a semi-algebraic condition. Let $g$ be a nonzero regular function on
$X$ which vanishes on the singular locus of $X$.
Let $F$ be the product
of $g$ and all the functions $f_i$
(assumed to be regular and not identically zero on $X$) appearing in the conditions
of the form 2.1 (i), (ii) and (iii)  in the description
of the semi-algebraic set $A$. Then we can take  $S$
to be the locus of $F= 0$ and
$$
A_i = \{x \in A \setminus \cL (S) \bigm \vert
{\rm ord}_t F (x) = i\}.
$$
Lemma 4.1 implies that $A_{i}$ is stable and (3.1.1) follows from
Lemma 4.4.
This proves the first
assertion. The proof of the second assertion is quite
similar. \hfill $\qed$

\bigskip 
\noindent (3.2) Let $S$ be an algebraic variety over $k$.
We write $\dim S \leq n$
if all the irreducible components of
$S$ have dimension $\leq n$.
Similarly, for $M$ in $\cM$,
we write $\dim M \leq n$ if
$M$ may be  expressed as a linear combination
of algebraic varieties with $\dim \leq n$.
For $m$ in $\ZZ$, we denote
by  $F^m \cM_{\rm loc}$
the subgroup of $\cM_{\rm loc}$ generated by
$\{ [S] \, \LL^{- i} \bigm \vert
i - \dim S \geq m\}$. This
defines a decreasing filtration $F^m$
on $\cM_{\rm loc}$.
We denote by $\widehat \cM$ the completion of 
$\cM_{\rm loc}$ with respect to this filtration.
We do not know whether or not the natural morphism
$\cM_{\rm loc} \rightarrow \widehat \cM$ is injective, but what is 
important for the applications in section 6 is the fact that the Euler 
characteristic and the Hodge polynomial
of an 
algebraic variety $S$ only depend on the image of $[S]$ in 
$\widehat \cM$, see (6.1) below.

We denote by $\overline \cM_{\rm loc}$ the image of 
$\cM_{\rm loc}$ in $\widehat \cM$, thus
$$\overline \cM_{\rm loc} = \cM_{\rm loc} / \cap_{m} F^m \cM_{\rm 
loc}.$$

\begin{enonce}{3.2}{Definition-Proposition}Let $X$ be an algebraic
variety
over $k$ of pure dimension $d$. Let $\BB$ be the set of all semi-algebraic
subsets
of $\cL (X)$. There exists a unique map $\mu : \BB \rightarrow
\widehat \cM$ satisfying the following three properties.
\begin{enumerate}
\item[] (3.2.1) \, If $A \in
\BB$ is stable at level $n$, then
$\mu (A) = [\pi_{n} (A)] \LL^{- (n + 1) d}$.
\item[] (3.2.2) \, If $A \in
\BB$ is contained in $\cL (S)$ with $S$ a closed subvariety of $X$ with 
${\rm dim} \, S < {\rm dim} \, X$, then $\mu (A) = 0$.
\item[] (3.2.3) \, Let $A_{i}$ be in $\BB$ for each $i$ in $\NN$.
Assume that the
$A_{i}$'s are mutually disjoint and that
$A := \bigcup_{i \in \NN} A_{i}$ is semi-algebraic. Then
$\sum_{i \in \NN} \mu (A_{i})$ converges in $\widehat \cM$
to $\mu (A)$.
\end{enumerate}
We call this unique map the {\it motivic volume} on $\cL (X)$
and denote it by $\mu_{\cL (X)}$ or $\mu$. Moreover we have
\begin{enumerate}
\item[] (3.2.4) \, If $A$ and $B$ are in $\BB$, $A \subset B$,
and if $\mu (B)$ belongs to the closure $F^{m} (\widehat \cM)$
of $F^{m} \cM_{\rm loc}$ in $\widehat \cM$,
then $\mu (A) \in F^{m} (\widehat \cM)$.
\end{enumerate}
Hence, for $A$ in $\BB$ and $\alpha : A \rightarrow \ZZ \cup
\{+ \infty\}$ a simple
function,
we can define
$$
\int_{A} \LL^{- \alpha} d \mu := \sum_{n \in \ZZ} \mu (A \cap \alpha^{-1} (n))
\, \LL^{- n}
$$
in $\widehat \cM$, whenever the right hand side
converges in $\widehat \cM$, in which case we say that $\LL^{- \alpha}$
is integrable on $A$. If the function
$\alpha$ is bounded from below, then
$\LL^{- \alpha}$
is integrable on $A$, because of (3.2.4).
\end{enonce}

\demo The key ingredients in the proof are the lemmas 2.4, 3.1 and 4.3.
The uniqueness of $\mu$ follows directly from Lemma 3.1, so it only
remains to prove the existence of a map $\mu : \BB \rightarrow \widehat
\cM$
satisfying (3.2.1) up to (3.2.4).

Let $\BB_{0}$ denote the set of all $A$ in $\BB$ which are stable. Thus
$\BB_{0}$ is closed under finite unions and finite
intersections. Clearly, there exists a map 
$\mu_{0} : \BB_{0} \rightarrow \widehat
\cM$ satisfying (3.2.1) and (3.2.4) with $\mu$ and $\BB$
replaced by 
$\mu_{0}$ and $\BB_{0}$. Obviously 
$\mu_{0}$ is additive, hence Lemma 2.4 yields 
(3.2.3) with $\mu$ and $\BB$
replaced by 
$\mu_{0}$ and $\BB_{0}$.
Next let $\BB_{1}$ be the set of all $A$ in $\BB$ which can be written as
$A = \bigcup_{i \in \NN} A_{i}$ with the $A_{i}$'s in $\BB_{0}$ mutually
disjoint
and $\lim_{i \rightarrow \infty} \mu_{0} (A_{i}) = 0$.
For $A$ in $\BB_{1}$ we set $\mu_{1} (A) = \sum_{i = 0}^{\infty} \mu_{0}
(A_{i})$. This is independent of the choice of the $A_{i}$'s.
Indeed, suppose that also
$A = \bigcup_{i \in \NN} A'_{i}$ with the $A'_{i}$'s in $\BB_{0}$ mutually
disjoint
and $\lim_{i \rightarrow \infty} \mu_{0} (A'_{i}) = 0$.
Then
\begin{eqn}
\sum_{i = 0}^{\infty} \mu_{0} (A_{i}) &=&
\sum_{i = 0}^{\infty} \mu_{0} (\bigcup_{j \in \NN} (A_{i} \cap A'_{j}))
=
\sum_{i = 0}^{\infty} \sum_{j = 0}^{\infty} \mu_{0}(A_{i} \cap A'_{j})\\
&=&
\sum_{j = 0}^{\infty} \sum_{i = 0}^{\infty} \mu_{0}(A_{i} \cap A'_{j})
= \sum_{j = 0}^{\infty} \mu_{0} (A'_{j})
\end{eqn}
because (3.2.3) and (3.2.4) hold for 
$\mu$ and $\BB$
replaced by 
$\mu_{0}$ and $\BB_{0}$.
One verifies that (3.2.1) and (3.2.4) are true for 
$\mu$ and $\BB$
replaced by 
$\mu_{1}$ and $\BB_{1}$.
From Lemma 4.3 one easily deduces the following
\begin{enumerate}
\item[] (3.2.5) \, If $S$ is a closed subvariety of $X$ with
${\rm dim} \, S < {\rm dim} \, X$ and if $A$ belongs to $\BB_{1}$, then
$A \setminus \cL (S)$ belongs also to 
$\BB_{1}$ and $\mu_{1} (A \setminus \cL (S)) = \mu_{1} (A)$.
\end{enumerate}

Indeed, we may assume $A$ belongs to $\BB_{0}$ and consider
the following partitions by elements of $\BB_{0}$:
$$
A \setminus \cL (S) =
(A \setminus \pi^{-1}_{m} \pi_{m} (\cL (S))) \cup \bigcup_{n \geq
m}\Bigl(
(\pi^{-1}_{n} \pi_{n} (\cL (S))
\setminus \pi^{-1}_{n + 1} \pi_{n + 1} (\cL (S)))
\cap A
\Bigr),
$$
$$
A = (A \setminus \pi^{-1}_{m} \pi_{m} (\cL (S))) \cup
(\pi^{-1}_{m} \pi_{m} (\cL (S)) \cap A),
$$
for $m \in \NN$ large enough.

Next let $A$ be any element of $\BB$. Then, by Lemma 3.1,
there exists a closed subvariety $S$ of $X$ with 
${\rm dim} \, S < {\rm dim} \, X$ such that $A \setminus \cL (S)$
belongs to $\BB_{1}$.
Define $\mu$ by $\mu (A) = \mu_{1} (A \setminus \cL (S))$.
By (3.2.5), this definition is independent of the choice of $S$.
Indeed, if $S'$ is another 
closed subvariety of $X$ with 
${\rm dim} \, S' < {\rm dim} \, X$ such that $A \setminus \cL (S')$
belongs to $\BB_{1}$, then
\begin{eqn}
\mu_{1} (A \setminus \cL (S')) &=& 
\mu_{1} ((A \setminus \cL (S')) \setminus \cL (S))
=
\mu_{1} ((A \setminus \cL (S)) \setminus \cL (S'))\\
&=& \mu_{1} (A \setminus \cL (S)).
\end{eqn}
Clearly (3.2.1), (3.2.2) and (3.2.4) are satisfied and $\mu$
is additive on finite disjoint unions. It remains to prove
(3.2.3). Let $A$ and the $A_{i}$'s 
be elements of $\BB$ as in (3.2.3) and let $m$
be in
$\NN$. By an argument entirely similar to the proof of Lemma 3.1,
there exists {\it weakly stable} 
$A'$ and $A'_{i}$'s in $\BB$ such that
$A \subset A'$, $A_{i} \subset A'_{i}$ and
$\mu (A) - \mu (A')$ and 
$\mu (A_{i}) - \mu (A'_{i})$
belong to $F^{m} \widehat \cM$.
Moreover, replacing $A_{i}$ by $A_{i} \cup (A' \setminus A)$ and
$A'_{i}$ by $A'_{i} \cap A'$, we may assume
$A_{i} \cup (A' \setminus A) \subset A'_{i} \subset A'$.
Hence $A' = \bigcup_{i \in \NN} A'_{i}$, and by Lemma 2.4, $A'$ is the
union of a finite number of the sets $A'_{i}$, thus $A' = \bigcup_{i =
1, \ldots,  e} A'_{i}$ whenever $e$ is large enough.
Since 
$$
A' = \Bigl(\bigcup_{i =
1, \ldots, e} A_{i}\Bigr) \cup
\Bigl(\bigcup_{i =
1, \ldots, e} (A'_{i} \setminus A_{i})\Bigr),
$$
we get
$$
\mu (A) \equiv \mu (A') \equiv \sum_{i = 1}^{e} \mu (A_{i}) \quad
\hbox{mod}
\quad F^{m} \widehat \cM.
$$
Because this holds for all $m$ in $\NN$, we obtain (3.2.3).  \hfill  $\qed$

\bigskip 
\noindent (3.3) Let $X$ be an algebraic variety over $k$ of dimension
$d$, and let $\cI$ be a coherent sheaf of ideals on $X$.
We denote by 
${\rm ord}_t \cI$
the function
${\rm ord}_t \cI : \cL (X) \rightarrow \NN \cup \{+\infty\}$ given by
$\varphi \mapsto \min_{g} {\rm ord}_t g (\tilde \varphi)$,
where the minimum is taken over all $g$ in the stalk $\cI_{\pi_{0}
(\varphi)}$
of $\cI$ at $\pi_{0}
(\varphi)$. Note that ${\rm ord}_t \cI$ is a simple function.
Let $\Omega^{1}_{X}$ be the sheaf of differentials on $X$ and let
$\Omega^{d}_{X}$ be the
$d$-th exterior power of 
$\Omega^{1}_{X}$. If $X$ is smooth and $\cF$ is a coherent sheaf
on $X$ together with a natural morphism
$\iota : \cF \rightarrow \Omega^{d}_{X}$, we denote by $\cI (\cF)$ the sheaf of ideals on $X$
which is locally generated by functions $\iota (\omega) / dx$ with $\omega$ a
local section of $\cF$ and $dx$ a local volume form on $X$.
Denote by 
${\rm ord}_t \cF$ the simple function
${\rm ord}_t \cI (\cF)$.

We have the following change of variables
formula for birational morphisms, which generalizes the one in [K].

\begin{lem}{3.3}Let $X$ and $Y$ be algebraic varieties over $k$, of pure
dimension $d$. Assume that $Y$ is {\it smooth}.
Let $h : Y \rightarrow X$
be a proper birational morphism,
$A$ be a semi-algebraic
subset of $\cL (X)$ and let $\alpha : A \rightarrow
\NN$ be a simple function.
Then
$$
\int_A \LL^{- \alpha} d \mu = 
\int_{h^{ - 1} (A)}
\LL^{- \alpha \circ h - {\rm ord}_t h^{\ast} (\Omega_{X}^{d})} d \mu.
$$
Moreover assume that $A \cap \cL (h (E)) = \emptyset$, where $E$
is the exceptional locus of $h$, and that $A$ and the fibers of 
$\alpha$
are weakly stable (and hence stable). Then $h^{- 1} (A)$ and the fibers
of $\alpha \circ h + \ord_{t} h^{\ast} (\Omega^d_{X})$
on $h^{- 1} (A)$ are stable and the above formula also holds for $\mu$
replaced by $\tilde \mu$.
\end{lem}

\demo The first assertion follows directly from Lemma 3.1 and Lemma 
3.4 below. Now assume that 
$A \cap \cL (h (E)) = \emptyset$ and that $A$ and the fibers of 
$\alpha$ are weakly stable. Since $X_{\rm sing} \subset h (E)$,
we have $A \cap \cL (X_{\rm sing}) = \emptyset$. Hence $A$ and the 
fibers of 
$\alpha$ are stable by (2.7). Moreover $h^{-1} (A)$ is stable, because 
$Y$ is smooth. Note that $\ord_{t} h^{\ast} (\Omega^d_{X})$ does not 
take the value $+ \infty$ on $h^{-1} (A)$, since $h^{-1} (A) \cap \cL
(E)
= \emptyset$,
and that the fibers of $\ord_{t} h^{\ast} (\Omega^d_{X})$ on $h^{- 1} 
(A)$ are weakly stable, hence stable. Thus, by (3.1), $\alpha \circ h$
and $\ord_{t} h^{\ast} (\Omega^d_{X})$ are bounded on
$h^{- 1} 
(A)$. We conclude that the fibers of $\alpha \circ h
+ \ord_{t} h^{\ast} (\Omega^d_{X})$
on $h^{- 1} (A)$ are stable, because they are finite unions of 
intersections of a fiber of
$\alpha \circ h$ with a fiber of
$\ord_{t} h^{\ast} (\Omega^d_{X})$. Note also that $A \subset 
\cL^{(e')} (X)$ for some $e'$ in $\NN$, by Lemma 2.4. The last 
assertion of Lemma 3.3 follows now directly from Lemma 3.4 below.  \hfill $\qed$

\begin{lem}{3.4}Let $X$ and $Y$ be algebraic varieties over $k$, of pure
dimension $d$
and let
$h : Y \rightarrow X$
be a birational morphism.
Assume that $Y$ is {\it smooth}.
For $e$ and $e'$
in $\NN$, let $\Delta_{e, e'}$ be the semi-algebraic subset of $\cL
(Y)$
defined
by
$$
\Delta_{e, e'} := \{\varphi 
\in \cL (Y) \bigm \vert ({\rm ord}_t h^{\ast} (\Omega_{X}^{d}))
(\varphi) 
= e \, \, \hbox{\rm and} \, \, h (\varphi) \in \cL^{(e')} (X)\},
$$
where
$\cL^{(e')} (X)$ is defined as in 2.6.
For $n$ in $\NN$, let
$h_{n \ast} : \cL_n (Y) \rightarrow \cL_n (X)$ be the morphism
induced by $h$, and let $\Delta_{e, e', n}$ be the image of
$\Delta_{e, e'}$ in $\cL_n (Y)$.
There exists $c$ in $\NN \setminus \{0\}$,
such that, for all $e, e', n$ in $\NN$ with $n \geq 2e$,
$n \geq e + ce'$,
the following holds.
\begin{enumerate}
\item[(a)]The set $\Delta_{e, e', n}$ is  a union of fibers of
$h_{n \ast}$.
\item[(b)]The restriction of $h_{n \ast}$ to $\Delta_{e, e', n}$
is a piecewise trivial
fibration with fiber
$\AA^e_k$ onto its image.
\end{enumerate}
\end{lem}

\demo Because $Y$ is smooth, the canonical morphism $\cL (Y) \rightarrow
\cL_{n} (Y)$ is surjective.
Consider the following assertion~:
\begin{eqn}
& & h_{n \ast}^{- 1} 
(h_{n \ast} (\bar \varphi)) 
\subset
\{\bar y \in \cL_{n} (Y) 
\bigm \vert
\bar \varphi \equiv \bar y 
\;\;\; {\rm mod}\, \cL_{n - e} (Y)\} \quad
\hbox{for all} \quad \bar \varphi \in
\Delta_{e, e', n},
\numeqn{(a$'$)}
\end{eqn}
where $\bar \varphi \equiv \bar y 
\;\;\; {\rm mod}\, \cL_{n - e} (Y)$ means that
$\bar \varphi$ and $\bar y$ have the same image in 
$\cL_{n - e} (Y)$.
Since $ n - e \geq e, e'$, the right hand side of
(a$'$) is contained in $\Delta_{e, e', n}$, and assertion
(a$'$) implies assertion (a).
Because $h$ is birational, $Y$ is smooth, $X$ is smooth at each point
of $h (\Delta_{e, e'})$,
and the ideal sheaf
$\cI (h^{\ast} (\Omega^{d}_{X}))$ does not vanish at any point
of $\Delta_{e, e'}$, we have that $h (y) \not= h (y')$ whenever
$y \not= y'$, $y \in \Delta_{e, e'}$, $y' \in \cL (Y)$.
Hence assertion (a$'$) is implied by the following assertion~:
\begin{eqn}
& & 
\hbox{For all} \, \, \varphi \in \Delta_{e, e'}, \, x \in \cL (X),
\, \hbox{with} \, \, h (\varphi) \equiv x 
\;\;\; {\rm mod}\, \cL_{n} (X),\numeqn{(a$''$)} \\
& &\, \hbox{there exists} \, \, y \in \cL (Y)
\, \hbox{with} \,
h (y) = x 
\, \hbox{and} \, \, 
\varphi \equiv y 
\;\;\; {\rm mod}\, \cL_{n - e} (Y).
\end{eqn}

Thus, we have only to prove 
(a$''$) and (b). For this we may assume that $X$ and $Y$ are affine.
Moreover we may assume that $Y \subset \AA^{M}_{k}$ and that the first
$d$
coordinates $y_{1}$, $y_{2}$, \dots, $y_{d}$ on 
$\AA^{M}_{k}$ induce an \'{e}tale map $Y \rightarrow \AA_{k}^{d}$.
Then, by Lemma 4.2 (with $n = e = 0$), the natural map
$\cL_{n} (Y)  \rightarrow Y \times_{\AA^{d}_{k}} \cL_{n} (\AA_{k}^{d})$
is an isomorphism. For ease of notation 
we will assume that
$Y = \AA^{d}_{k}$. The general case can be proved in the
same way,  identifying $\cL (Y)$ with $\cL (\AA_{k}^{d})$ on each fixed
fiber of $\cL (Y) \rightarrow Y$, and taking for
$\cJ_{h}$ below the jacobian matrix of $h$ with respect to the system
of
local coordinates $y_{1}$, $y_{2}$, \dots, $y_{d}$ on $Y$.

Let us first prove (a$''$) and (b) in the special case where also 
$X = \AA^{d}_{k}$.
Let $\varphi$ be in $\Delta_{e, e'} (\bar k) \subset \cL (Y) (\bar k) = Y ( \bar
k [[t]])$. Denote by
$\cJ_{h}$ the jacobian matrix of $h$. For 
(a$''$) we have to prove that 
for all $v$ in $\bar k [[t]]^d$
there exists $u$ in $\bar k [[t]]^d$ such that
\begin{eqn}& &
h (\varphi + t^{n + 1 -e} u) 
=
h (\varphi) + t^{n + 1} v.
\numeqn{(1)}
\end{eqn}
By Taylor expansion, (1) is equivalent to
\begin{eqn}& &
t^{-e} \cJ_{h} (\varphi) u +
t \, \hbox{(higher order terms in $u$)} =
v.
\numeqn{(2)}
\end{eqn}
The existence of a solution $u$ follows from Hensel's Lemma,
since the inverse of $t^{-e} \cJ_{h} (\varphi)$
has entries in $\bar k [[t]]$, because ${\rm ord}_{t} {\rm det} \,
\cJ_{h} (\varphi) = e$. This proves 
(a$''$), hence (a$'$) and (a), when $X = \AA^{d}_{k}$.

Let us prove now assertion (b) assuming
$X = \AA^{d}_{k}$. Let the morphism
$s : \cL_n (X) \rightarrow \cL (X)$
be a section of the projection $\pi_{n} :
\cL (X) \rightarrow \cL_n (X)$.
For $\bar x$ in $h_{n \ast} (\Delta_{e, e', n} (\bar k))$, $s (\bar x)$
belongs to $h (\Delta_{e, e'})$, because  of (1), and hence $h^{- 1}$ is regular at
$s (\bar x)$. The mapping
$\theta : \bar x \mapsto h^{-1}  (s (\bar x))$
from 
$h_{n \ast} (\Delta_{e, e', n})$ to  $\Delta_{e, e'}$ is 
a piecewise morphism, meaning there exists a finite partition of
the domain of $\theta$ into locally closed subvarieties of $\cL_{n}
(X)$,
such that the restriction of $\theta$ to any subvariety is a morphism of
schemes. 
For $\bar x$ in $h_{n \ast} (\Delta_{e, e', n} (\bar k))$, we deduce from (a$'$)
and (2) that
\begin{eqn}& &
h_{n \ast}^{- 1} (\bar x) =
\{ \theta (\bar x) + t^{n + 1 - e} u \;\;\; {\rm mod}\, t^{n + 1}
\bigm \vert
u \in \bar k [[t]]^d
\, \hbox{and} \, \,
(\cJ_h (\theta (\bar x))) u \equiv 0 \;\;\; {\rm mod}\, t^{e}\}.\numeqn{(3)}
\end{eqn}
Thus we see that the fiber $h_{n \ast}^{- 1} (\bar x)$ can be identified with
a linear subspace of $\{u \; {\rm mod}\, t^e \bigm \vert u \in \bar k [[t]]^d\} \simeq
\AA^{de}_k$ given by linear equations whose coefficients 
in $\bar k [[t]]$
are functions of $\bar x$ which are
piecewise morphisms
on $h_{n \ast} (\Delta_{e, e', n}) \subset \cL_n (X)$.
Moreover, for a fixed $\bar x$ in
$h_{n \ast} (\Delta_{e, e', n} (\bar k))$,
the jacobian matrix $\cJ_h (\theta (\bar x))$ is equivalent over
$\bar k [[t]]$ to a diagonal matrix with diagonal elements
$t^{e_1}, \ldots, t^{e_d}$, with $e_1, \ldots, e_d$ in $\NN$,
and  $e_1 + \cdots + e_d = e$. Together with (3),  this gives an isomorphism
\begin{eqn}&&
h_{n \ast}^{- 1} (\bar x) \simeq \AA^e_{\bar k},\numeqn{(4)}
\end{eqn}
and it is now easy to verify that
$h_{n \ast \vert_{\Delta_{e, e', n}}}$ is a piecewise trivial fibration
onto its image.
This proves (b) when
$X = \AA^d_{k}$.

We now turn to the case where
$X = {\rm Spec}\, (k [x_{1}, \ldots, x_{N}] / I) \subset \AA^{N}_{k}$
and $Y = \AA^d_{k}$.
Denote again by $\cJ_{h}$ the jacobian matrix of $h : Y \rightarrow X
\hookrightarrow \AA_{k}^{N}$. 
By the argument at the beginning of the proof of Lemma 4.1, we see
that there exists $c$ in $\NN \setminus \{0\}$, such that, for any $e, e'$ in $\NN$,
the set $h (\Delta_{e, e'})$ is covered by a finite number of semi-algebraic
subsets
$A$ of $\cL (X)$, which are weakly stable at level $c \, e'$, such that, for each such
$A$, the following holds~: On $A$ the variety $X$ is a ``complete
intersection'', meaning that
$$
\cL (X) \cap A
=
\cL ({\rm Spec}\, (k [x_{1}, \ldots, x_{N}] / (f_{1}, \ldots, f_{N -
d})))
\cap A,
$$
for suitable $f_{1}, \ldots, f_{N -
d} \in I$ (which may depend on $A$), and moreover , for each $x$
in
$A$, there exists an $N - d$ by 
$N - d$ minor $\delta'$ of the matrix $\Delta := \frac{\partial (f_{1},
\ldots, f_{N - d})}{\partial (x_{1},
\ldots, x_{N})}$ satisfying ${\rm ord}_{t} \delta'(\tilde x) \leq
c \, e'$.
Denote by $\delta$ the minor of $\Delta$ formed by the first $N - d$
columns.
Up to shrinking $A$ and renumbering the coordinates, we may assume
there exists $e''$ in $\NN$, with $e'' \leq c \, e'$, such that, for all
$x$ in $A$,
\begin{eqn}
e''  = {\rm ord}_{t} \delta (\tilde x) &\leq& {\rm ord}_{t}
\delta' (\tilde x)
\numeqn{(5)}
\end{eqn}
for all $N - d$ by 
$N - d$ minors $\delta'$
of $\Delta$.
It is enough to prove (a$''$) and (b) for $\Delta_{e, e'}$ replaced by
$\Delta_{e, e'} \cap h^{-1} (A)$. From now on let $\varphi$ be
in
$\Delta_{e, e'} (\bar k) \cap h^{-1} (A)$
and assume $n \geq 2 e$, $n \geq e + ce'$. Since the product of the jacobian
matrices
$\Delta (h (\varphi))$ 
and $\cJ_{h} (\varphi)$
yields zero, one checks that the $d$ by $d$ minor
of 
$\cJ_{h} (\varphi)$ formed by the last $d$ rows has minimal
${\rm ord}_{t}$, equal to $e$, among all 
$d$ by $d$ minors
of 
$\cJ_{h} (\varphi)$. Indeed, the columns of 
$\cJ_{h} (\varphi)$ are solutions of the homogeneous linear system
of equations with matrix $\Delta (h (\varphi))$, and the first
$N - d$ components of any solution are fixed $\bar k [[t]]$-linear
combinations of the last $d$ components by Cramer's rule and (5).
Thus the first $N - d$ rows of 
$\cJ_{h} (\varphi)$ are $\bar k [[t]]$-linear
combinations of the last $d$ rows.

Let $p : X \rightarrow \AA^{d}_{k}$ denote the projection on the last
$d$ coordinates $(x_{1}, \ldots, x_{N}) \mapsto
(x_{N - d + 1}, \ldots, x_{N})$.
Denote by $\cJ_{p \circ h}$ the jacobian matrix of
$p \circ h$, {\it i.e.} $\cJ_{p \circ h}$ consists in the last $d$ rows
of the matrix $\cJ_{h}$. Thus
$$
e = ({\rm ord}_{t} h^{\ast} (\Omega_{X}^{d})) (\varphi)
=
{\rm ord}_{t} {\rm det} \, \cJ_{p \circ h} (\varphi).
$$

To prove
(a$''$) we have to show that 
for all $v$ in $\bar k [[t]]^N$, satisfying
$$
h (\varphi) + t^{n + 1} v \in \cL (X) (\bar k) \subset \bar k [[t]]^{N},
$$
there exists $u$ in $\bar k [[t]]^d$ such that
\begin{eqn}& &
h (\varphi + t^{n + 1 -e} u) 
=
h (\varphi) + t^{n + 1} v.
\numeqn{(1$'$)}
\end{eqn}
It follows from Lemma 4.2 (with $n$, $e$ replaced by $n - e$, $e''$)
(or alternatively by a direct argument using Taylor expansion)
that (1$'$) is equivalent to
\begin{eqn}& &
p \circ h (\varphi + t^{n + 1 -e} u) 
=
p \circ h (\varphi) + t^{n + 1} p (v).
\numeqn{(1$''$)}
\end{eqn}
The proof of (a$''$) and (b) proceeds now in the same way as
in the special case $X = \AA^{d}_{k}$ treated above from (1) to (4).
We only have to replace $\cJ_{h}$ by
$\cJ_{p \circ h}$, $v$ by $p (v)$, (1) by (1$'$), $\Delta_{e, e'}$ by
$\Delta_{e, e'} \cap h^{-1 } (A)$, $\Delta_{e, e', n}$ by
the image of $\Delta_{e, e'} \cap h^{-1 } (A)$ in $\cL_{n} (Y)$, and to take
for $s : \cL_{n} (X) \rightarrow \cL (X)$ a section of the projection
$\pi_{n} : \cL (X) \rightarrow \cL_{n} (X)$, whose restriction to
$\pi_{n} (A)$ is a piecewise morphism.
The existence of such a section $s$ is insured by Lemma 4.2. The
essential
point is that (3) remains valid when one replaces
$\cJ_{h}$ by $\cJ_{p \circ h}$, because
the first $N - d$ rows of
$\cJ_{h} (\varphi)$ are 
$\bar k [[t]]$-linear
combinations of the last $d$ rows.
With these modifications, the argument from (1) to (4) remains valid
in the present situation. This ends the proof of Lemma 3.4.  \hfill $\qed$

\rem When $X$ is smooth and $\cI$ is the ideal sheaf of an effective
divisor $D$ on $X$, the motivic integral
$\int_{\cL (X)} \LL^{- {\rm ord}_{t} \cI} d \mu$ was first introduced by
Kontsevich [K] and denoted by him $[\int_{X} e^{D}]$.

\bigskip

\noindent (3.5) In this subsection we consider a generalization
which will not be used in the present paper.
Let $X$ be an algebraic variety over $k$ of pure dimension $d$, and let
$\cF$
be a coherent sheaf on $X$
together with  a natural morphism $\iota~: \cF \rightarrow
\Omega^{d}_{X}$. We
denote by
${\rm ord}_{t} (\cF)$ the simple function from
$\cL (X) \setminus \cL (X_{\rm sing})$ to
$\NN \cup \{+ \infty\}$ defined by
\begin{eqn}
&& ({\rm ord}_{t} \cF) (\varphi)
=
\min_{\omega_{1} \in \cF_{\pi_{0} (\varphi)}} \,
\max_{\omega_{2} \in (\Omega^{d}_{X})_{\pi_{0} (\varphi)}}
{\rm ord}_{t} \, \frac{\iota (\omega_{1})}{\omega_{2}}
(\tilde \varphi),
\numeqn{(3.5.1)}
\end{eqn}
for any $\varphi$ in $\cL (X) \setminus \cL (X_{\rm sing})$.
(Here $X_{\rm sing}$ denotes the singular locus of $X$.)
This definition of ${\rm ord}_{t} (\cF)$ coincides with the one given in
(3.3) when $X$ is smooth. By (3.2), the map
$$
\BB \rightarrow \widehat \cM : \qquad
A \mapsto
\int_{A \setminus \cL (X_{\rm sing})}
\LL^{- {\rm ord}_{t} \cF} \, d\mu_{\cL (X)}
$$
is a $\sigma$-additive measure on $\cL (X)$ which we will
denote by $\mu_{\cF}$. Note that $\mu_{\Omega^{d}_{X}} = \mu_{\cL (X)}$.
Let $A$ be a semi-algebraic subset of $\cL (X)$ and
let $\alpha : A \rightarrow \NN \cup \{+ \infty \}$ be a simple
function.
We define the motivic integral
$\int_{A} \LL^{-\alpha} d \mu_{\cF}$ as in (3.2), but with $\mu$
replaced by
$\mu_{\cF}$. If $Y$ is any algebraic variety over $k$ of pure dimension
$d$ and $h : Y \rightarrow X$ is a proper birational morphism, then
\begin{eqn}
\int_{h^{-1} (A)} \LL^{- \alpha \circ h}
d \mu_{h^{\ast} (\cF)}
=
\int_{A} \LL^{- \alpha} d\mu_{\cF}.
\numeqn{(3.5.2)}
\end{eqn}
Indeed, this follows by considering a resolution of $Y$ and applying
Lemma 3.3 twice.

\sec{4}{Some lemmas}In this section we prove some lemmas which were
already used in sections 2 and 3 and which will be used again in section 7.

\medskip

\begin{lem}{4.1}Let $X$ be an algebraic variety over $k$ of pure dimension
$d$, and assume the notation of 2.6.
There exists $c$ in $\NN \setminus \{0\}$ such that,
for all $e$, $n$ in $\NN$ with $n \geq ce$, the following holds.
\begin{enumerate}
\item[(a)]The map
$$
\theta_n : \pi_{n + 1} (\cL (X)) \rightarrow \pi_{n} (\cL (X))
$$
is a piecewise trivial fibration over $\pi_n (\cL^{(e)} (X))$
with fiber $\AA_k^d$.
\item[(b)]Moreover
$$
[\pi_n (\cL^{(e)} (X))] = [\pi_{ce} (\cL^{(e)} (X))] \, \LL^{d (n- ce)}.
$$
\end{enumerate}
\end{lem}

\demo Since (b) is a direct consequence of (a), we only have to prove
(a). We may assume $X$ is affine, say
$X = {\rm Spec} \, (k[x_1, \ldots, x_N] / I) \subset \AA_k^N$.
Clearly $X_{\rm sing}$ is the intersection of hypersurfaces having
an equation of the form
$h \delta = 0$ where $h$ belongs to
$k[x_1, \ldots, x_N]$ and $\delta$ is some $N - d$ by $N - d$
minor of the matrix
$\Delta :=
\frac{\partial (f_1, \ldots, f_{N - d})}{\partial (x_1, \ldots, x_N)}$,
with $f_1, \ldots, f_{N - d}$ in $I$ and $h I \subset
(f_1, \ldots, f_{N - d})$. Hence, by Hilbert's Nullstellensatz,
there exists $c$ in $\NN \setminus \{0\}$ such that,
for any $e$ in $\NN$, the set $\cL^{(e)} (X)$ is
covered by a finite number of sets of the form
$$
A := \{\varphi \in \cL (\AA^N_k) \bigm \vert 
(h \delta)(\tilde \varphi) \not\equiv 0 \;\;\; \hbox{mod} \, t^{ce + 1}\}.
$$
Thus it is sufficient to prove that, for $n \geq ce$, the map
$\theta_n$ is a piecewise trivial fibration over
$\pi_n (\cL (X) \cap A)$ with fiber
$\AA^d_k$. Since $h (\tilde \varphi) \not= 0$ for all
$\varphi$ in $A$, 
$$
\cL (X) \cap A =
\cL ({\rm Spec} \, k [x_1, \ldots, x_N] / (f_1, \ldots, f_{N -d})) \cap A.
$$
Hence we may assume $I = (f_1, \ldots f_{N - d})$, {\it i.e}. $X$ is a complete
intersection.
Let $e' \in \NN$, $e' \leq ce$, and set
$$
A' :=
\{\varphi \in A \bigm \vert {\rm ord}_t \delta (\tilde \varphi) = e'
\, \hbox{and} \,
{\rm ord}_t \delta' (\tilde \varphi) \geq e'
\, \hbox{for all} \, N - d \, \hbox{by} \, N - d \, \hbox{minors} \,
\delta' \, \hbox{of} \Delta\}.
$$
It is sufficient to prove that the map $\theta_n$ is a piecewise
trivial fibration over
$\pi_n (\cL (X) \cap A')$ with fiber $\AA_k^d$.
We may assume $\delta$ is the minor of the first $N- d$ columns of $\Delta$.
Let $s : (\bar k [t] / t^{n + 1})^N \rightarrow \bar k [t]^N$
be the $\bar k$-linear map given by
$t^{\ell} \, \hbox{mod} \, t^{n + 1} \mapsto t^{\ell}$ for
$\ell = 0, 1, \ldots, n$. Let $\bar \varphi \in (\bar k [t] / t^{n + 1})^N$
be any $\bar k$-rational point of $\pi_n (\cL (X) \cap A')$. We have
$$
\theta_n^{- 1} (\bar \varphi)
=
\{s (\bar \varphi) + t^{n + 1}y \;\;\; \hbox{mod} \, t^{n + 2} \bigm \vert
y \in \bar k [[t]]^N, f (s (\bar \varphi) + t^{n  + 1} y) = 0
\},
$$
where $f$ is the column with entries $f_1, \ldots, f_{N - d}$.
By Taylor expansion, the condition
$f (s (\bar \varphi) + t^{n  + 1} y) = 0$ can be rewritten as
\begin{eqn}
f (s (\bar \varphi)) + t^{n + 1} \Delta (s (\bar \varphi)) y
+ t^{2 (n + 1)} (\cdots) = 0.
\numeqn{(1)}
\end{eqn}
There exists an $N - d$ by $N-d$ matrix $M$ over
$k [x_1, \ldots, x_N]$, independent of the choice
of
$\bar \varphi$, such that
$$
M \Delta = (\delta I_{N - d}, W),
$$
where $I_{N - d}$ is the identity matrix with $N - d$ columns
and $W$ is an $N - d$ by $d$ matrix
such that $W (s (\bar \varphi)) \equiv 0 \, \hbox{mod} \,t^{e'}$.
Indeed, to check this last congruence, one expresses the last $d$ columns
of $\Delta$ in terms of the first $N- d$ columns by Cramer's rule and then one
uses the definition of $A'$.

Condition (1) is equivalent to
\begin{eqn}
t^{-e' - n -1} (Mf) (s (\bar \varphi)) + t^{- e'} (M \Delta) (s (\bar \varphi)) y
+ t^{n + 1 -e'} (\cdots) = 0.
\numeqn{(1$'$)}
\end{eqn}
Note that $t^{- e'} (M \Delta) (s (\bar \varphi))$ is a matrix over
$\bar k [[t]]$, whose minor  determined by the first
$N - d$ columns is not divisible by $t$, because 
${\rm ord}_t \delta (\tilde \varphi) = e'$. Moreover, $n + 1 - e' \geq
1$.
Since $\bar \varphi$ is liftable to $\cL (X)$
({\it i.e.} belongs to $\pi_{n} (\cL (X))$, equation (1$'$) has a
solution
$y$ in $\bar k [[t ]]^{N}$, and thus
$t^{- e' -n - 1} (Mf) (s (\bar \varphi))$ is a column matrix over 
$\bar k [[t]]$. By Hensel's Lemma, we deduce that
$\theta_{n}^{-1} (\bar \varphi)$ is equal to the set of all
$s (\bar \varphi)  + t^{n + 1} y_{0}$, with $y_{0}$ in
${\bar k}^{N}$ such that
$$
t^{- e' - n - 1} (Mf) (s (\bar \varphi)) +
t^{-e'} (M \Delta) (s (\bar \varphi)) y_{0} \equiv 0 \;\;\; \hbox{mod} \, t.
$$
Thus the fiber $\theta_{n}^{-1} (\bar \varphi)$ is a $d$-dimensional
affine subspace of $\AA^{N}_{k}$, given by linear equations which
express the first $N - d$ coordinates in terms of linear combinations
of the last $d$ coordinates, with coefficients which are regular
functions on
each locally closed subset of $\cL_{n} (X)$ contained in
$\pi_{n} (\cL (X) \cap A')$. This proves that 
$\theta_{n}$ is a piecewise trivial fibration over
$\pi_{n} (\cL (X) \cap A')$ with fiber $\AA^{d}_{k}$. \hfill $\qed$

\medskip

Assume  now $X = {\rm Spec} \, (k [x_{1}, \ldots,
x_{N}] / (f_{1}, \ldots, f_{N - d}))$, with
$f_{1}$, \ldots, $f_{N - d}$ in $k [x_{1}, \ldots,
x_{N}]$. Fix $m, n, e$ in $\NN$ with $m > n \geq e$.
Set
$$\Delta := \frac{\partial (f_{1},
\ldots, f_{N - d})}{\partial (x_{1},
\ldots, x_{N})}, \quad \quad \delta := {\rm det} \,
\frac{\partial (f_{1},
\ldots, f_{N - d})}{\partial (x_{1},
\ldots, x_{N-d})},
$$ and
$$
A := \{\varphi \in \cL (X) \bigm \vert
{\rm ord}_{t} \delta (\tilde \varphi) = e \leq
{\rm ord}_{t} \delta' (\tilde \varphi)
\, 
\hbox{\rm for every}
\,
N- d
\, \hbox{\rm by} \, N - d
\, 
\hbox{\rm minor}
\,
\delta'
\,\hbox{\rm of}
\, \Delta
\}.
$$
Let $p : X \rightarrow \AA^{d}_{k}$  denote the projection
onto the last $d$ coordinates.
Denote by $\kappa$
the natural map $\kappa : \pi_{m} (A) \rightarrow
\pi_{n} (A) \times_{\cL_{n} (\AA^{d}_{k})}
\cL_{m} (\AA^{d}_{k})$
induced by
$\varphi \in \cL (X) \mapsto (\pi_{n} (\varphi), p (\varphi))$,
where the fiber product is with respect to the map
$\pi_{n} (A) \subset \cL_{n} (X)) \rightarrow \cL_{n} (\AA^{d}_{k})$
induced by $p$. 

\begin{lem}{4.2}The map $\kappa$ is a bijection. Moreover,
for any locally closed subvariety $Z$ of
$\cL_{n} (X) \times_{\cL_{n} (\AA^{d}_{k})}
\cL_{m} (\AA^{d}_{k})$ contained in the domain of
$\kappa^{-1}$, the restriction of
$\kappa^{-1}$ to $Z$ is a morphism from $Z$ to
$\cL_{m} (X)$.
\end{lem}

\demo By induction it is enough to treat the case where $m = n + 1$, but
this case is just a reformulation of material in the proof of Lemma
4.1.  \hfill $\qed$

\begin{lem}{4.3}Let $X$ be an algebraic variety over $k$ of dimension
$d$.
\begin{enumerate}
\item[(1)]For any $n$ in $\NN$,
$$
{\rm dim} \, \pi_{n} (\cL (X)) \leq (n + 1) d.
$$
\item[(2)]For any $n$, $m$ in $\NN$, with $m \geq n$,
the fibers of
$\pi_{m} (\cL (X)) \rightarrow \pi_{n} (\cL (X))$ are of dimension
$\leq (m - n ) d$.
\end{enumerate}
\end{lem}

\demo This lemma is probably  well known and
is implicit in
[O]. Assertion (1) follows from
assertion (2). Moreover, it suffices to prove (2) for $m = n + 1$,
and we may assume that $X$ is affine and that $k = \bar k$.
Substituting $a_{i} + t^{n + 1} x_{i}$ for the affine coordinates
$x_{i}$ in the equations defining $X$ in some affine space,
we see that each fiber of
$\pi_{n + 1} (\cL (X)) \rightarrow \pi_{n} (\cL (X))$ is contained in
the reduction mod $t$ of a scheme which is flat and of finite type over
${\rm Spec} \, (k [t])$ with generic fiber $X \otimes_{k} k(t)$.
But the reduction mod $t$ of such a scheme 
has dimension at most $d$. \hfill $\qed$

\bigskip 
\noindent (4.4) Let $X$ be an algebraic
variety over $k$. By Greenberg's theorem [G], for any $n$ in $\NN$
there exists $\gamma (n) \geq n$ in $\NN$ such that
$\pi_n (\cL (X))$ is the image of
$\cL_{\gamma (n)} (X)$ by the natural projection
and which is minimal for this property.
Furthermore the function $\gamma (n)$ is bounded by a real
linear function of $n$.
We call $\gamma$ the Greenberg function for $X$. It has been studied in
[L-J], [H].

\begin{lem}{4.4}Let $X$ be an algebraic
variety over $k$ of dimension $d$
and let $S$ be a closed subvariety of dimension $< d$. Let
$\gamma$ be the Greenberg function for $S$.
For any $n$, $i$, $e$ in $\NN$, with $n \geq i \geq \gamma (e)$,
$\pi_{n, X} (\pi_{i, X}^{-1} \cL_{i} (S))$ is of dimension $\leq (n + 1)
\, 
d - e - 1$.
\end{lem}

\demo We drop the subscript $X$ in $\pi_{n}$, etc. Clearly we may assume
$i = \gamma (e)$. By Lemma 4.3 (2) applied to the projection
$$
\pi_{n} (\pi_{\gamma (e)}^{-1} \cL_{\gamma (e)} (S))
\rightarrow
\pi_{e} (\pi_{\gamma (e)}^{-1} \cL_{\gamma (e)} (S))
$$
we obtain 
$$
{\rm dim} \, \pi_{n} (\pi_{\gamma (e)}^{-1} \cL_{\gamma (e)} (S))
\leq
(n - e) \, d
+
{\rm dim} \, \pi_{e} (\pi_{\gamma (e)}^{-1} \cL_{\gamma (e)} (S)).
$$
Since,
by definition of the Greenberg function,
$\pi_{e} (\pi_{\gamma (e)}^{-1} \cL_{\gamma (e)} (S)) = \pi_{e} (\cL
(S))$,
the result follows because, by Lemma 4.3 (1),
${\rm dim} \, \pi_{e} (\cL
(S)) \leq (e + 1) \, (d - 1)$. \hfill $\qed$

\sec{5}{Proof of Theorem 1.1 and rationality results}

\begin{th}{5.1}Let $X$ be an algebraic variety over $k$ of pure
dimension
$d$.
Let $A_n$, $n \in \ZZ^r$,
be a semi-algebraic family of semi-algebraic subsets of $\cL (X)$ and 
let
$\alpha : \cL (X)  \times \ZZ^r
\rightarrow \NN$ be a simple function. Then the power series
\begin{eqn}
\sum_{n \in \NN^r} T^n \, \int_{A_{n}} \LL^{- \alpha (\_, n)} d\mu
&&\numeqn{(1)}
\end{eqn}
in the variable $T = (T_{1}, \ldots, T_{r})$ belongs to the subring of
$\widehat \cM [[T]]$ generated by the image in 
$\widehat \cM [[T]]$ of $\cM_{\rm loc} [T]$, $(\LL^i
 - 1)^{-1}$ and $(1- \LL^{-a} T^b)^{- 1}$, with $i \in \NN \setminus 
 \{0\}$, $a \in \NN$, $b \in \NN^r \setminus \{0\}$.
\end{th}

\begin{cor}{}For any semi-algebraic subset $A$ of $\cL (X)$, the measure
$\mu (A)$ is in $\overline
\cM_{\rm loc} [((\LL^i - 1)^{-1})_{i \geq 1}]$, where
$\overline
\cM_{\rm loc}$ is the image of $\cM_{\rm loc}$ in $\widehat \cM$, 
cf. (3.2).  \hfill $\qed$
\end{cor}

\begin{th}{5.1$'$}Let $X$ be an algebraic variety over $k$ of pure
dimension
$d$.
Let $A_n$, $n \in \ZZ^r$,
be a semi-algebraic family of semi-algebraic subsets of $\cL (X)$ and 
let
$\alpha : \cL (X)  \times \ZZ^r
\rightarrow \NN$ be a simple function. 
Assume that $A_{n} \cap \cL (X_{\rm sing}) = \emptyset$ and that
$A_{n}$ and the fibers of
$\alpha (\_, n) : A_{n} \rightarrow \NN$
are weakly stable (and hence stable), for every $n \in \NN^r$.
Then the power series
\begin{eqn}
\sum_{n \in \NN^r} T^n \, \int_{A_{n}} \LL^{- \alpha (\_, n)} d \tilde \mu
&&\numeqn{(1)}
\end{eqn}
in the variable $T = (T_{1}, \ldots, T_{r})$ belongs to the subring of
$\cM_{\rm loc}[[T]]$ generated by $\cM_{\rm loc} [T]$
and the series
$(1- \LL^{-a} T^b)^{- 1}$, with  $a \in \NN$ and
$b \in \NN^r \setminus \{0\}$.
\end{th}

\noindent {\it Proof of Theorem 5.1$'$}. --- Using a resolution
of singularities $\pi : \tilde X \rightarrow X$, with exceptional
locus $\pi^{-1} (X_{\rm sing})$ and using Lemma 3.3, we may assume
that $X$ is smooth. Moreover we may also assume that $X$ is affine.
For $n$ in $\NN^r$, $m$ in $\NN$, set
$$
A_{n, m} := \{x \in A_{n} \bigm \vert \alpha (x, n) = m\}.
$$
For every $n$, the map $\alpha (\_, n) : A_{n} \rightarrow \NN$
is bounded, cf. (3.1). Hence by Theorem 2.1 there exists a Presburger
function $\theta : \ZZ^r \rightarrow \NN$ such that $\alpha (x, n) 
\leq
\theta (n)$ for every $x \in A_{n}$, $n \in \NN^r$. Thus the series 
(1) is equal to 
\begin{eqn}
\sum_{{n \in \NN^r, \, m \in \NN}\atop{m \leq \theta (n)}} \tilde \mu
(A_{n, m}) \, \LL^{- m} \, T^n.&&\numeqn{(2)}
\end{eqn}
By Lemma 2.8, the family $(A_{n, m})_{n \in \NN^r, \, m \in \NN}$ is a 
finite boolean combination of semi-algebraic families
which have bounded representations. Recall that each member of a family
with bounded representation is weakly stable and hence stable, 
because $X$ is smooth. Using that
$\tilde \mu (A' \cup A'') = \tilde \mu (A') + \tilde \mu (A'')
- \tilde \mu (A' \cap A'')$, we may assume that the family
$(A_{n, m})_{n, m}$ is a finite intersection of families each  of which
has the property that it or
its complement has bounded representation. Using that
$\tilde \mu (A' \setminus A'') = \tilde \mu (A') - \tilde \mu (A' 
\cap A'')$, we may assume that the family
$(A_{n, m})_{n, m}$ is a finite union of finite intersections of 
families with bounded representation. Since a finite intersection of 
families with bounded representation has bounded representation, we 
may finally assume that the family $(A_{n, m})_{n, m}$ has bounded 
representation.
We may further assume that in the bounded representation
(2.8) of $(A_{n, m})_{n, m}$, no covering by affine open subsets $U$
is needed.
Let $F$ be the
product of all the functions $f_i$
(assumed to be regular on $X$ and to have bounded
order on each $A_{n, m}$) appearing in the conditions of the form
2.1 (i), (ii) and (iii)  in the bounded representation of the family
$(A_{n, m})_{n, m}$.

Let $h : Y \rightarrow X$ be an embedded resolution of singularities of
the locus 
of $F= 0$ in $X$.
The exceptional locus of $h$ is contained in $h^{-1} (F^{- 1} (0))$.
We can cover $Y$ by affine open subsets $U$ on which
there exist regular functions
$z_1, \ldots, z_d$ inducing an {\'e}tale map
$U \rightarrow \AA^d_k$ such that on $U$
each $f_i \circ h$ is a monomial in
$z_1, \ldots, z_d$ multiplied by a regular function with no zeros on 
$U$. For such an open subset $U$ we may assume
that the $z_{i}$'s appearing in at least one of these monomials
are exactly $z_{1}$, $z_{2}$, \dots, $z_{d_{0}}$. Since
$\ord_{t} F$ is bounded on $A_{n, m}$, Lemma 3.3 yields that,
uniformly in $n$, $m$,
$\tilde \mu (A_{n, m})$ is a finite 
$\ZZ$-linear combination of
terms of the form
\begin{eqn}
\sum_{{\ell_1, \ldots, \ell_{d_{0}} \in \NN}\atop{\theta (\ell_1,
\ldots, \ell_{d_{0}},n, m)}} 
\LL^{- \beta (\ell_1,
\ldots, \ell_{d_{0}})} \, \tilde \mu
(W_{\ell_1, \ldots, \ell_{d_{0}}}),&&\numeqn{(3)}
\end{eqn}
where
$\theta (\ell_1,
\ldots, \ell_{d_{0}},n, m)$ is a semi-algebraic condition defining
a Presburger subset of
$\ZZ^{d_{0} + r + 1}$,
$\beta$
is a linear form with coefficients in $\NN$, and
$W_{\ell_1, \ldots, \ell_{d_{0}}}$ is of the form
$$\{
y \in \cL (U) \bigm \vert
{\rm ord}_t z_i = \ell_i \, {\rm for} \,
i = 1, \ldots, d_{0}, \quad{\rm and} \quad
(\overline{ac} (z_1), \ldots,
\overline{ac} (z_{d_{0}}), \pi_0 (y)) \in W\},$$
with $W$ a constructible subset of $(\AA^{1}_{k} \setminus \{0\})^{d_{0}}
\times U$ and $U$ as above.
(Recall that $\pi_0 : \cL (Y) \rightarrow \cL_0 (Y) = Y$
is the natural projection.) Note that the sum in (3) is finite for each 
fixed $n$, $m$, since $\ord_{t} F$ is bounded on $A_{n, m}$.

From Lemma 4.2 with $n = e = 0$ we obtain
\begin{eqn}
\tilde \mu
(W_{\ell_1, \ldots, \ell_{d_{0}}}) &=&
[W'] \, \LL^{-(\sum_{i = 1}^{d_{0}} \ell_{i}) -d},\numeqn{(4)}
\end{eqn}
where $W'$ is the set of
$(w_{1}, \ldots,w_{d_{0}}, y)$'s in $W \subset (\AA^1_{k} \setminus 
\{0\})^{d_{0}} \times U$ such that $z_{i} (y) = 0$
when $\ell_{i} > 0$ and $z_{i} (y) = w_{i}$
when $\ell_{i} = 0$.
We conclude by (2), (3) and (4), that the series (1) is a finite
$\cM_{\rm loc}$-linear combination of series
$f (\LL^{- 1}, T_{1}, \ldots, T_{r})$,
with  $f (X_{1}, X_{2}, \ldots, X_{r + 1})$ in
$\ZZ [X_{1}][[X_{2}, \cdots, X_{r + 1}]]$
given as in Lemma 5.2 below (with $r$ replaced by $r + 1$).
The theorem follows now directly from Lemma 5.2 and Lemma 5.3
below.  \hfill $\qed$

\bigskip

\noindent {\it Proof of Theorem 5.1}. ---  The argument is the same
as in the proof of Theorem 5.1$'$, but easier. Now we do not have
to work 
with bounded representations anymore, since $\widehat \cM$ is 
complete, and  we are allowed to work with infinite sums (which have no 
meaning in $\cM_{\rm loc}$). Finally we obtain that 
the series (1) is a finite
$\overline \cM_{\rm loc}$-linear combination of series
$f (\LL^{- 1}, T_{1}, \ldots, T_{r})$,
with  $f (X_{1}, X_{2}, \ldots, X_{r + 1})$ in
$\ZZ [[X_{1}, \cdots, X_{r + 1}]]$
given as in Lemma 5.2 below (with $r$ replaced by $r + 1$).
However the series 
$f (X_{1}, X_{2}, \ldots, X_{r + 1})$ might not be in
$\ZZ [X_{1}][[X_{2}, \cdots, X_{r + 1}]]$, so we cannot apply 
Lemma 5.3 (as in the proof of Theorem 5.1$'$). We leave the details to 
the reader.  \hfill $\qed$

\begin{lem}{5.2}Let $P$ be a Presburger
subset of $\ZZ^{m}$ and let $\varphi_{1} : \ZZ^m
\rightarrow \NN$, \dots, $\varphi_{r} :
\ZZ^m
\rightarrow \NN$ be Presburger functions. Assume that the fibers of 
the map
$\varphi : P \rightarrow \NN^r$ given by
$i \mapsto (\varphi_{1} (i), \ldots, \varphi_{r} (i))$
are finite. Then the power series
$f (X) := \sum_{i \in P} X^{\varphi (i)}$, in the variable
$X = (X_{1}, \ldots, X_{r})$, belongs to the subring of
$\ZZ [[X]]$ generated by $\ZZ [X]$ and the series
$(1 - X^c)^{- 1}$, with $c \in \NN^r \setminus \{0\}$.
\end{lem}

\demo We may assume that $P \subset \NN^m$. We first
consider
the special case where $r = m$ and $\varphi (i) = i$. By replacing $i$ 
by $di + a$ for suitable $d \in \NN \setminus \{0\}$, $a \in \NN^m$, 
we may assume that no congruence relations appear in the description
of $P$. Since $\sum_{P_{1} \cup P_{2}} = \sum_{P_{1}} +
\sum_{P_{2}} -\sum_{P_{1} \cap P_{2}}$, we may further assume that 
$P$ is the set of integral points in a rational convex polyhedron 
$\overline P \subset \RR_{+}^m$. (Here $\RR_{+} :=
\{x \in \RR \bigm \vert x \geq 0\}$.) A straightforward calculation 
yields the lemma when $\overline P$ is a polyhedral cone generated by 
part of a basis of $\ZZ^n$. Hence the lemma also holds in our special 
case when $\overline P$ is any rational convex polyhedral cone, since
any such cone can be decomposed in cones generated by part of a basis
of $\ZZ^n$ (see, {\it e.g.}, [Da] p.123-124).
Let ${\overline P}' \subset \RR_{+}^{m + 1}$
be a
rational convex polyhedral cone such that $\overline P$ is equal to
the intersection of ${\overline P}'$ with the hyperplane given by the 
equation $x_{m + 1} = 
1$ in $\RR^{m + 1}$. Set
$$g (X, T) := \sum_{{i \in \NN^m, \, \ell \in \NN}\atop{(i, \ell)
\in {\overline P}'}} X^i \, T^{\ell},$$ where $T$ is a new variable.
Since ${\overline P}'$ is a polyhedral cone, we already know that
$g (X, T)$ belongs to the subring of $\ZZ [[X, T]]$ generated by
$\ZZ [X, T]$ and the series
$(1 - X^c T^a)^{- 1}$, with $(c, a) \in \NN^{r + 1} \setminus \{0\}$.
Since $f = \frac{\partial g}{\partial T} (X, 0)$, this proves Lemma 
5.2 in the special case where $r = m$ and $\varphi (i) = i$.
For the general case, consider $$h (Y, X) := 
\sum_{i \in P, \, j = \varphi (i)} Y^i X^j,$$ where $Y = (Y_{1}, \ldots,
Y_{m})$.
The special case implies that $h (Y, X)$ belongs to
the subring of $\ZZ [[Y, X]]$ generated by
$\ZZ [Y, X]$ and the series $(1 - Y^b X^c)^{- 1}$ with
$(b, c) \in \NN^{m + r} \setminus \{0\}$. Note that
$h (Y, X)$ belongs to
$\ZZ [Y][[ X]]$, since the fibers of $\varphi$ are finite. By Lemma 
5.3 we can take in the above
$c \not= 0$, which finishes the proof of Lemma 5.2,
since $f (X) = h (1, X)$. \hfill $\qed$

\begin{lem}{5.3}Let $R$ be a commutative ring with unit. Let $Y = (Y_{1},
\ldots, Y_{m})$ and $X = (X_{1}, \ldots, X_{r})$ be variables. Assume 
that
$h (Y, X)$ belongs to
the subring of $R [[Y, X]]$ generated by
$R [Y, X]$ and the series $(1 - Y^a X^b)^{- 1}$ with
$(a, b) \in \NN^{m + r} \setminus \{0\}$. If, moreover, 
$h (Y, X)$ belongs to $R [Y][[ X]]$, then we can take in the above
$b \not= 0$.
\end{lem}

\demo Straightforward.  \hfill $\qed$

\bigskip

\noindent(5.4) {\it Proof of Theorem 1.1}.--- We will prove the
following theorem of which Theorem 1.1 is a special case.

\begin{th}{5.4}Let $X$ be an algebraic variety over $k$ and let
$A$ be a
semi-algebraic subset of $\cL (X)$. The power series
$$P_{A} (T) := \sum_{n = 0}^{\infty} [\pi_n (A)] \, T^n,$$
considered as an element of 
$\cM_{\rm loc} [[T]]$, is rational and belongs to
$\cM [T]_{{\rm loc}}$.
\end{th}

\demo We may assume $X$ is a closed subvariety
of a smooth connected  algebraic variety $W$ of dimension $d$.
That $P_{A} (T)$ is in  $\cM [T]_{{\rm loc}}$ 
follows directly from Theorem 5.1$'$
because
$$
[\pi_{n, X} (A)] \, \LL^{- (n + 1) d} = 
\int_{\pi_{n, W}^{- 1} (\pi_{n, X} (A))} \LL^0 d \tilde
\mu_{\cL (W)}
$$
and the family $\pi_{n, W}^{- 1} (\pi_{n, X} (A))$, $n \in \NN$,
is a semi-algebraic
family of stable
semi-alge\-braic subsets of $\cL (W)$
by Theorem 2.1.  \hfill $\qed$

\rem Note that the above proof of Theorem 5.4
only uses the material in \S\kern .15em 3 with
$X$ smooth. Moreover it does not use (3.2).

\sec{6}{Formula for $\mu (\cL (X))$ and applications}

\noindent (6.1) Let $S$ be an algebraic variety over $\CC$. Consider the Hodge
numbers
$$
e^{p, q} (S) := \sum_{i \geq 0} (- 1)^i h^{p, q} (H^i_c (S, \CC)),
$$
where $h^{p, q} (H^i_c (S, \CC))$ denotes the rank of the
$(p, q)$-Hodge component of the $i$-th cohomology group
with compact supports. Define the Hodge polynomial of $S$ as
$$
H (S; u, v) := \sum_{p, q} e^{p, q} (S) \, u^p v^q.
$$
Note that $H (\AA^1_{\CC}; u, v) = uv$ and that
${\rm deg} \, H (S; u, v) =  2 \, {\rm dim} \, S$.
The map $S \mapsto H (S; u, v)$ factors through
$\cM$ and induces a ring morphism
$H : \cM \rightarrow \ZZ [u, v]$,
which extends naturally to
a ring morphism
$H : \cM_{\rm loc} \rightarrow \ZZ [u, v][(uv)^{-1}]$.
Clearly ${\rm deg} \, H (Z; u, v) \leq  -2 m$
when $Z \in F^m \cM_{\rm loc}$, hence the kernel of the natural map
$\cM_{\rm loc} \rightarrow \widehat \cM$ is killed by $H$, and $H$ induces
a ring morphism
$$H : \overline \cM_{\rm loc} \longrightarrow \ZZ [u, v][(uv)^{-1}],$$
where $\overline \cM_{\rm loc}$ is the image of $\cM_{\rm loc}$ in
$\widehat \cM$.

We denote by
$\chi (S)$ the topological
Euler characteristic of $S$, {\it i.e.} the alternating sum of the rank
of
its Betti or de Rham cohomology groups. Clearly
$\chi (S) = H (S; 1, 1)$. Hence the map $S \mapsto \chi (S)$ factors through
$\overline \cM_{\rm loc}$ and induces a ring morphism
$\chi : \overline \cM_{\rm loc}\rightarrow \ZZ$.
Such
morphisms $H$ and $\chi$ can also be 
defined when $k \not=\CC$, using an embedding into $\CC$
of the field of definition of $S$. Indeed, the $e^{p,q} (S)$'s do not
depend on the embedding, since for a smooth projective variety $S$ they
are equal to $(- 1)^{p + q} {\rm dim} \, H^{q} (S, \Omega^{p}_{S})$.

More generally let ${\rm Mot}_k$ denote the category of Chow motives over $k$. It
follows from [G-S] and [G-N] that there exists a 
natural morphism
$\cM_{\rm loc} \rightarrow K_0 ({\rm Mot}_k)$.
The conjectural existence of a weight filtration
(with the expected properties) on the objects of the category
${\rm Mot}_k$
would imply that this morphism
factors through 
$\overline \cM_{\rm loc}$. In particular, without using any conjecture, this is true
when we replace Chow motives by their {\'e}tale or Hodge realizations.

\bigskip
\noindent (6.2) Let $X$ be an algebraic variety over $k$ of pure
dimension
$d$, and let 
$h : Y \rightarrow X$ be a resolution of $X$. By this we
mean
$Y$ is a smooth algebraic variety over $k$, $h$ is birational,
proper and defined over $k$, and the exceptional
locus $E$ of $h$ has normal crossings, meaning that the $k$-irreducible
components of $E$ are smooth and intersect transversally.
Let us denote the $k$-irreducible components of $E$ by $E_{i}$, $i \in
J$.
For $I \subset J$, set $E_{I} = \bigcap_{i \in I} E_{i}$ and
$E_{I}^{\circ} = E_{I} \setminus \bigcup_{j \not\in I} E_{j}$.
Note that the support of the subscheme defined by
the ideal sheaf  $\cI (h^{\ast} (\Omega^{d}_{X}))$
is contained in $E$.
(See (3.3) for the notation $\cI (\_)$.)

\bigskip
\noindent (6.3) Assume the hypothesises of (6.2) hold and moreover that
the ideal sheaf  $\cI (h^{\ast} (\Omega^{d}_{X}))$
is invertible. Such a resolution $h$ always
exists by Hironaka's Theorem. For $i$ in $I$, set
\begin{eqn}
\quad &&\nu_{i} = 1 +
(\hbox{multiplicity of } \, E_{i} \, \hbox{in the divisor of} \,\,
Y \, \hbox{determined by} \,  \cI (h^{\ast} (\Omega^{d}_{X}))).
\numeqn{(6.3.1)}
\end{eqn}

\begin{prop}{6.3.2}Let $W$ be a closed subvariety of $X$. The following
formula holds in $\widehat \cM$,
\begin{eqn}
\mu_{\cL (X)} (\pi_{0}^{-1} (W))
=
\LL^{- d} \sum_{I \subset J} [E_{I}^{\circ} \cap h^{-1} (W)] \,
\prod_{i \in I} \frac{\LL - 1}{\LL^{\nu_{i}} -1} \, .
\numeqn{(6.3.2)}
\end{eqn}
\end{prop}

\demo It follows directly from Lemma 3.3 (cf. the proof of 5.1). \hfill $\qed$

\begin{cor}{6.3.3}Let $X$ be a variety over $k$ of pure dimension $d$.
Then $\mu (\cL (X))$ belongs to
$\overline \cM_{\rm loc} [(\frac{\LL - 1}{\LL^{i} - 1})_{i \geq 1}]$. \hfill $\qed$
\end{cor}

The morphism 
$\chi : \overline \cM_{\rm loc} \rightarrow \ZZ$ defined in (6.1)
may be extended uniquely to a ring morphism
$\chi : \overline \cM_{\rm loc} [(\frac{\LL - 1}{\LL^{i} - 1})_{i \geq 1}]
\rightarrow \QQ$ by requiring
$\chi (\frac{\LL - 1}{\LL^{i} - 1}) = \frac{1}{i}$ for $i \geq 1$.

\begin{cor}{6.3.4}Let $X$ be a variety over $k$ of pure dimension $d$.
Then the Euler characteristic
$\chi (\mu (\cL (X)))$ is well defined as a rational number. \hfill $\qed$
\end{cor}

Similarly, the morphism
$H : \overline \cM_{\rm loc} \rightarrow
\ZZ [u, v][(uv)^{-1}]$ defined in (6.1)
may be extended naturally to
a morphism of rings
$$H : \overline \cM_{\rm loc}
\biggl[\Bigl(\frac{\LL - 1}{\LL^{i} - 1}\Bigr)_{i \geq 1}\biggr]
\longrightarrow \ZZ [u, v][(uv)^{-1}]\biggl[\Bigl(\frac{uv - 1}
{(uv)^{i} - 1}\Bigr)_{i \geq
1}\biggr].$$
So if $X$ is  a variety over $\CC$
of pure dimension $d$, we have
$H (\mu (\cL (X)); u, v)$ which is well defined in
$\ZZ [u, v][(uv)^{-1}][(\frac{uv - 1}{(uv)^{i} - 1})_{i \geq
1}]$.

\bigskip
\noindent (6.4) Assume the hypothesises of (6.2) hold.
Instead of requiring that $\cI (h^{\ast} (\Omega^{d}_{X}))$ is invertible we will
now assume 
$\cI (h^{\ast} (\Omega^{d}_{X}))$
is locally generated by elements whose zero
loci are contained in
$E$. This weaker condition is more practical for explicit
computations. To $x$ in $E_{I}^{\circ}$,
$I \subset J$,
we associate a polyhedron $\Delta_{x} \subset
\RR^{|I|}$
as follows. 
For $i$ in $I$ and $g$ a section of
$\cI (h^{\ast} (\Omega^{d}_{X}))$ on an open set of $Y$ intersecting $E_{i}$,
set
\begin{eqn}
\quad \quad \, \,\,  \, \nu_{i}  (g) = 1 +
(\hbox{multiplicity of $E_{i}$ in the divisor
of $g$}).
\numeqn{(6.4.1)}
\end{eqn}
Set
\begin{eqn}
\Delta_{x} := \hbox{convex hull of} \, \,  \bigcup_{g}
\Bigl((\nu_{i} (g))_{i \in I} + \RR_{+}^{|I|}\Bigr),
\numeqn{(6.4.2)}
\end{eqn}
where $g$ runs over all local sections at $x$
of
$\cI (h^{\ast} (\Omega^{d}_{X}))$  whose zero
loci are contained in
$E$,
and
$\RR_{+} = \{a \in \RR \bigm | a \geq 0\}$.
On each $k$-irreducible component $C$ of $E_{I}^{\circ}$,
$\Delta_{x}$ is constant with value say
$\Delta_{C}$.
The support function $\ell (\_, \Delta)$ of a polyhedron $\Delta$ in
the interior of
$\RR_{+}^{|I|}$ is given, for $\xi$ in $\RR_{+}^{|I|}$, by
\begin{eqn}
\ell (\xi, \Delta) := \min_{v \in \Delta} \, \xi \cdot v.
\numeqn{(6.4.3)}
\end{eqn}
We set
\begin{eqn}
Z (\Delta) := (\LL - 1)^{|I|} \sum_{\xi \in (\NN^{\times})^{|I|}}
\LL^{- \ell (\xi, \Delta)} \in \widehat \cM,
\numeqn{(6.4.4)}
\end{eqn}
where
$\NN^{\times} = \NN \setminus \{0\}$.
Note that $(\NN^{\times})^{|I|}$ may be partitioned into finitely many
sets $C$ of the form
$\NN^{\times} \eta_{C, 1} + \cdots +  \NN^{\times} \eta_{C, e_{C}}$,
with $\eta_{C, 1}, \ldots,  \eta_{C, e_{C}}$ in $\NN^{|I|}$
linearly independent, on which 
$\ell (\_, \Delta)$ is linear. Thus
\begin{eqn}
Z (\Delta) = (\LL - 1)^{|I|} \sum_{C} \prod_{i = 1, \ldots, e_{C}}
(\LL^{\ell (\eta_{C, i}, \Delta)} - 1)^{-1}.
\numeqn{(6.4.5)}
\end{eqn}

\begin{prop}{6.4.6}Let $W$ be a closed subvariety of $X$. The following
formula holds in $\widehat \cM$,
\begin{eqn}
\mu_{\cL (X)} (\pi_{0}^{-1} (W))
=
\LL^{- d} \sum_{C} \, [C \cap h^{-1} (W)] \,
Z (\Delta_{C}) \, ,
\numeqn{(6.4.6)}
\end{eqn}
where $C$ runs over all $k$-irreducible components of $E_{I}^{\circ}$,
for all $I \subset J$.
\end{prop}

\demo It follows directly from (6.4.4) and
Lemma 3.3 (cf. the proof of 5.1). \hfill $\qed$

\bigskip
\noindent (6.5) Let $\cI$ be an invertible ideal sheaf on $X$
and let $h : Y \rightarrow X$ be a
resolution of $X$ such that both
$E$ and the locus of $\cI \cO_{Y}$
are contained in a divisor with normal crossings
$E'$ in $Y$
(cf.  the notation in (6.2)).
The previous results may be extended
directly to  integrals of the form
$\int_{\pi_{0}^{-1} (W)} \LL^{- \ord_{t} \cI} d \mu$.
More precisely,
change notation from (6.2) so that the $E_{i}$'s, $i \in J$,
denote now the $k$-irreducible components of $E'$.
Let $N_{i}$ be the multiplicity of $E_{i}$ in the divisor
of $\cI \cO_{Y}$. Replace $\nu_{i}$ in (6.3), {\it resp.}
$\nu_{i} (g)$ in (6.4), by
$\nu_{i}+ N_{i}$, 
{\it resp.}
$\nu_{i} (g)+ N_{i}$. Then
$\int_{\pi_{0}^{-1} (W)} \LL^{- \ord_{t} \cI} d \mu$ equals the right hand side
of
(6.3.2), {\it resp.} (6.4.6), when the hypothesis of (6.3), {\it resp.}
(6.4), with $E$ replaced by $E'$,
is satisfied.

\bigskip
\noindent (6.6) {\sl Remark}. --- Assume $X$ is normal and let
$\omega_{X}$ be the
canonical sheaf of $X$.
By definition $\omega_{X} := j_{\ast} (\Omega^{d}_{X^{\circ}})$,
where $j : X^{\circ}
\hookrightarrow X$ denotes the inclusion of the smooth locus of $X$.
Assume moreover $X$ is Gorenstein, {\it i.e.} $\omega_{X}$ is
invertible,
and that all singularities of $X$ are canonical, meaning that there
exists a resolution
$h : Y \rightarrow X$ with $h^{\ast} (\omega_{X}) \subset
\Omega_{Y}^{d}$.
In this case Kontsevich [K] introduced the invariant
\begin{eqn}
\int_{\cL (Y)} \LL^{-{\rm ord}_{t} h^{\ast} (\omega_{X})}
d \mu_{Y},\numeqn{(6.6.1)}
\end{eqn}
which does not depend on the resolution $h$. Indeed the independence
follows from
Lemma 3.3 (with both varieties smooth) and the fact that two resolutions
are always dominated by a third one. If the resolution is crepant ({\it
i.e.}
$h^{\ast} (\omega_{X}) = \Omega^{d}_{Y}$), then
(6.6.1) equals $\LL^{- d} \sum_{I \subset J} [E_{I}^{\circ}] =
\LL^{-d} [Y]$ because of (6.3.2). In this way Kontsevich [K]
showed that $H (Y; u, v)$ is independent of the crepant resolution
$h$ (assuming such a resolution exists). This independence was
conjectured
in [B-D]. For some related constructions, see [V].

\sec{7}{Motivic volume as a limit}The following
result is an analogue of a result
by Oesterl\'e [O].

\begin{th}{7.1}Let $X$ be an algebraic variety over $k$ of pure dimension
$d$. Let $A$ be a semi-algebraic subset of $\cL (X)$.
The sequence $[\pi_n (A)] \, \LL^{- (n + 1)d}$ converges
in $\widehat \cM$ to $\mu (A)$. Furthermore when $A = \cL (X)$
this limit is nonzero.
\end{th}

\demo By 4.4 and Lemma 4.3 we
may assume that $X$ is affine and irreducible and that $A$ is given
by a semi-algebraic condition. As in the proof
of Lemma 3.1, let $g$ be a nonzero regular function on
$X$ which vanishes on the singular locus of $X$, let
$F$ be the product
of $g$ and all the functions $f_i$
(assumed to be regular and not identically zero on $X$)
appearing in the conditions
of the form 2.1 (i), (ii) and (iii)  in the description
of the semi-algebraic set $A$, and let $S$ 
be the locus of $F= 0$ in $X$.
For $e \geq 0$, set
$
A^{(e)}
:=  A \setminus \pi_{e, X}^{-1} \cL_{e} (S)
$.
Let $\gamma$ be the Greenberg function for $S$.
It follows from Lemma 4.4 that, for $n \geq \gamma (e)$,
\begin{eqn}
[\pi_n (A)] \, \LL^{- (n + 1)d} -
[\pi_n (A^{(\gamma (e))})] \, \LL^{- (n + 1)d}\numeqn{(1)}
\end{eqn}
belongs
to $F^{e + 1} \cM_{\rm loc}$. 
On $A^{(\gamma (e))}$, the function
$\ord_{t} F (x)$ is bounded by a multiple of $\gamma (e)$.
Hence, since $X_{\rm sing} \subset S$ and
$A^{(e)} \subset \cL^{(e)} (X)$, it follows from Lemma 4.1 that
$A^{(\gamma (e))}$ is stable at level 
$c \, \gamma (e)$, for some $c \geq 1$ independent of $e$.
Thus,
for $n$ and 
$n' \geq c \, \gamma (e)$,
\begin{eqn}
[\pi_n (A)^{(\gamma (e))}] \, \LL^{- (n + 1)d} =
[\pi_{n'} (A)^{(\gamma (e))}] \, \LL^{- (n' + 1)d}.\numeqn{(2)}
\end{eqn}
We deduce that, for $n$ and 
$n' \geq c \, \gamma (e)$,
$
[\pi_n (A)] \, \LL^{- (n + 1)d} -
[\pi_{n'} (A)] \, \LL^{- (n' + 1)d}
$ belongs
to $F^{e + 1} \cM_{\rm loc}$. This implies
that $[\pi_n (A)] \, \LL^{-(n + 1)d}$ is a Cauchy sequence,
hence converges to a limit
$\cA$ in $\widehat \cM$.
Since, by (3.2.1), $\mu (A^{(\gamma (e))}) =
[\pi_n (A^{(\gamma (e))})] \, \LL^{- (n + 1)d}$ when $n
\geq c \, \gamma (e)$,
it follows that
$\mu (A^{(\gamma (e))})$ converges to
$\cA$ in $\widehat \cM$, and we deduce that
$\cA = 
\mu (A)$, since, by (3.2.3), $\mu (A^{(\gamma (e))})$
converges to $\mu (A)$.
Because $\gamma (0) = 0$ and $\pi_{0} (\AA^{(0)}) = \pi_{0} (A) \setminus
S$,
we obtain from (1) and (2) (with $e = 0$, $n' = 0$, $n \rightarrow
\infty$)
that 
$\mu (A) - [\pi_0 (A) \setminus S] \, \LL^{- d}$
belongs to $F^{1} \widehat \cM$.
When $A = \cL (X)$,
we obtain  that $\mu (\cL (X))$ is nonzero, since
$[X \setminus S] \, \LL^{-d} \not\in F^{1} \widehat \cM$
(this last fact may be checked by considering
the degree of 
Hodge polynomials, cf. (6.1)). \hfill $\qed$

\rem It follows from Theorem 7.1 and Lemma 4.3 that, for any
algebraic variety $X$ over $k$ of dimension $d$, and any
semi-algebraic subset $A$ of $\cL (X)$,
$[\pi_n (A)] \, \LL^{- (n + 1)d}$ converges in
$\widehat \cM$ to $\mu_{\cL (X_{d})} (A \cap
\cL (X_{d}))$, where $X_{d}$ is the
union of the $d$-dimensional components of $X$.

\sec{8}{Remarks on the Greenberg function}Let $X$ be an algebraic
variety over $k$. The Greenberg function $n \mapsto \gamma (n)$
for $X$
has been defined in 4.4.
For $n \geq m$, we denote by $\pi_{n,m}$ the canonical morphism
$\cL_{n} (X) \rightarrow \cL_{m} (X)$. We set
$X_{n, j} := \pi_{n + j, n} \cL_{n + j} (X)$
and consider the Poincar\'e series
$$
P (T, U) = \sum_{n, j \in \NN} [X_{n, j}] \, T^n U^j
$$
as an element of  $\cM_{\rm loc} [[T, U]]$.

\medskip
The following result is an analogue of [D1] Theorem 5.1 and
Proposition 5.2.

\begin{th}{8.1}
\begin{enumerate}
\item [(1)]The series $P (T, U)$ is rational and can be written as
$$
P (T, U)  = q (T, U) (1 - U)^{-1}
\prod_{i = 1}^e (1 -\LL^{a_i} T^{b_i} U^{c_i})^{-1}
$$
with $q (T, U)$ in $\cM_{\rm loc} [T, U]$, $a_i$ in
$\ZZ$, $b_i$, $c_i$
in $\NN$, and $b_i \geq 1$.
\item[(2)]Suppose none of the factors
$1 -\LL^{a_i} T^{b_i} U^{c_i}$ divides 
$q (T, U)$ in $\cM_{\rm loc} [T, U]$.
Then $c := 1 + {\rm sup}_i c_i / b_i$ is the smallest real number such 
that there exists $d$ in $\RR$ with
$\gamma (n) \leq cn + d$ for all $n$ in $\NN$.
\end{enumerate}
\end{th}

\demo The proof is
very much the same as the one in [D1]. By a similar
proof as the one of Theorem 1.1 one obtains that
$P (T, U)$ has the form required in (1), except for the fact that
$b_i \geq 1$, which is proved in the same way as in [D1]. The proof
of (2) is also similar to the one of [D1]. However
we have to use
that the image in $\cM_{\rm loc}$ of a non-empty constructible
subset of an algebraic variety is nonzero, a fact that may be checked
by considering the degree of Hodge polynomials, cf. (6.1). \hfill $\qed$

\rem It follows from the result of Pas (Theorem 2.1) that
there exists a finite
partition of $\NN$ into congruences classes, such that the function
$\gamma$ is linear on each class for $n$ large enough.

\noindent J.D. University of Leuven, Department of Mathematics,

\noindent Celestijnenlaan 200B, 3001 Leuven, Belgium

\noindent{\tt Jan.Denef@wis.kuleuven.ac.be}
\bigskip

\noindent F.L. Centre de Math\'ematiques,
Ecole Polytechnique,

\noindent F-91128 Palaiseau
(URA 169 du CNRS)

\noindent et

\noindent Institut de Math\'{e}matiques,
Universit\'{e} P. et M. Curie, Case 82,
4 place Jussieu,

\noindent F-75252 Paris Cedex 05
(UMR 9994 du CNRS)

\noindent{\tt loeser@math.polytechnique.fr}

\end{document}